\documentclass[a4paper,11pt]{elsarticle}
\usepackage[utf8]{inputenc}
\usepackage{amsmath}
\usepackage{amsthm}
\usepackage{amsfonts}
\usepackage{amssymb}
\usepackage{empheq}
\usepackage{siunitx}
\usepackage{xcolor}
\graphicspath{{IMG/}}
\usepackage{tikz}
\usepackage{pgfplots}
\pgfplotsset{compat=1.17} 
\usepackage{enumitem}
\usepackage{subcaption}
\usepackage{hyperref}
\usepackage{cleveref}
\usepackage{float}
\usepackage[margin=0.8in]{geometry}
\usepackage{natbib}
\usepackage{comment}
\usepackage{autonum}

\makeatletter
\def\ps@pprintTitle{%
  \let\@oddhead\@empty
  \let\@evenhead\@empty
  \def\@oddfoot{}
  \let\@evenfoot\@oddfoot
}
\makeatother

\newtheorem{definition}{Definition}
\numberwithin{definition}{section}

\newtheorem{assumption}{Assumption}
\numberwithin{assumption}{section}

\newtheorem{theorem}{Theorem}
\numberwithin{theorem}{section}

\numberwithin{proposition}{section}
\theoremstyle{plain}
\newtheorem{remark}
{Remark}
\theoremstyle{plain}
\newtheorem{lemma}{Lemma}[section]
\usepackage{tikz}

\newcommand\T{\rule{0pt}{0pt}}
\newcommand\B{\rule{0pt}{0pt}}

\newcommand{\jump}[1]{[\mkern-1.5mu [#1] \mkern-1.5mu]}
\newcommand{\avg}[1]{\{ \mkern-5mu \{#1 \} \mkern-5mu \}}
\newcommand{\avgg}[1]{\left\{ \mkern-10mu \left\{#1 \right\} \mkern-10mu \right\}}

\renewcommand{\div}[1]{\nabla \mkern-2.5mu \cdot \mkern-2.5mu {#1}}
\newcommand{\divh}[1]{\nabla_h \mkern-2.5mu \cdot \mkern-2.5mu {#1}}

\definecolor{myred}{rgb}{0.9, 0.0, 0.0}
\definecolor{myblue}{rgb}{0.0, 0.28, 0.67}
\definecolor{mygreen}{rgb}{0.0, 0.7, 0.0}
\definecolor{myyellow}{rgb}{1.0, 0.55, 0.0}
\definecolor{mypurple}{rgb}{0.5, 0.1, 0.5}

\usepackage[thinlines]{easymat}

\begin{document}
\begin{frontmatter}
\title{Numerical modelling of wave propagation phenomena in thermo-poroelastic media via discontinuous Galerkin methods}

\author[mox]{Stefano Bonetti\corref{cor1}}
\author[mox]{Michele Botti}
\author[mox]{Ilario Mazzieri}
\author[mox]{Paola F. Antonietti}
\cortext[cor1]{\href{mailto:stefano.bonetti@polimi.it}{stefano.bonetti@polimi.it}}

\address[mox]{MOX Department of Mathematics - Politecnico di Milano, P.zza Leonardo da Vinci 32,Milano, Italy}

\begin{abstract}
We present and analyze a high-order discontinuous Galerkin method for the space discretization of the wave propagation model in thermo-poroelastic media. The proposed scheme supports general polytopal grids.  
Stability analysis and $hp$-version error estimates in suitable energy norms are derived for the semi-discrete problem. The fully-discrete scheme is then obtained based on employing an implicit Newmark-$\beta$ time integration scheme. A wide set of numerical simulations is reported, both for the verification of the theoretical estimates and for examples of physical interest. A comparison with the results of the poroelastic model is provided too, highlighting the differences between the predictive capabilities of the two models.
\end{abstract}

\end{frontmatter}

\section{Introduction}
\label{sec:introduction}
This paper deals with the numerical analysis of the fully-dynamic thermo-poroelastic model, that describes the wave propagation phenomena in thermo-poroelastic media. The study of these phenomena finds application in many fields, such as greenhouse gas reservoirs and geothermal energy extraction, that are crucial for environmental sustainability, or the study of thermo-elastic seismic energy release as a source mechanism for volcanic earthquakes, and induced earthquakes by human geological activities.

The theory of wave propagation in porous media has been first presented by Biot \cite{Biot1962} and then developed by Carcione \cite{Carcione2014}. Biot considered a fully-saturated porous media and investigated the presence of three kind of waves: two compressional ($P$) waves and a shear ($S$) wave. The two $P$-waves propagate in different ways, the first one -- denoted by ($E$) -- is a fast wave, while the second one -- denoted by (Biot) -- is a slow wave, that is diffusive at low frequencies and is slower than the fast wave $E$ at high frequencies. Biot has been the first to propose a model for poroelastic wave propagation taking into account the effect of temperature \cite{Biot1956}; this formulation was based only on the Fourier law for heat conduction. The sole presence of the diffusive operator - due to its nature - may lead to non-physical results, as it yields an infinite diffusion velocity of the temperature. To overcome this issue, Lord and Shuman \cite{Lord1967} introduced a relaxation term, based on the generalization of the Fourier law formulated by Cattaneo \cite{Cattaneo1948}. In general, the introduction of the coupling with the temperature induces the appearance of an additional $T$-wave (thermal), which is a slow diffusion wave. The analysis of the wave induced by temperature effects can be found also in \cite{Carcione2018, Yang2022}, where the thermo-elasticity problem is studied, and in \cite{Straughan2011}, where it is possible to find a comprehensive review on heat wave propagation. 

Our model - proposed in \cite{Carcione2019, Santos2021} -  is constituted by three equations: mass conservation, momentum conservation, and energy conservation, under the hypothesis of a linear thermo-poroelastic medium, and the inertial terms are included in all three equations. The poroelastic component of the model is written in the two-displacements formulation, namely where the unknowns are the solid and filtration displacements, while the pressure field is recovered post-processing. The relaxation terms in the energy conservation equation are included in order to get physically-consistent results. 

Other formulations for thermo-poroelastic coupling are available in the literature. As mentioned before, some works do not consider the relaxation terms \cite{Biot1956, Deresiewicz1957}. In \cite{Noda1990} the inertial terms in the temperature equation are not included, but a non-linear advection term is present in the model. Nield and Bejan in \cite{Nield2006} first included the heat transfer between the solid and the fluid, but they still need to include the third-order terms to correctly represent the physical behavior, cf. also \cite{Sharma2008} for a similar model. Among the literature on thermo-poroelasticity, we can also mention \cite{Iesan2014, Kumar2017}, even if their analysis is not based on Biot equations nor Darcy law, thus it is not directly comparable with our model. Last, we mention some recent works on the quasi-static thermo-poroelastic problem \cite{AntoniettiBonetti2022, Brun2019, Brun2020, Chen2022}, where both the linear and non-linear cases are treated. In particular, in \cite{AntoniettiBonetti2022} a discontinuous Galerkin approximation of the fully-coupled problem including a nonlinear advection term is considered.

For the spatial discretization of the problem we propose a discontinuous Galerkin finite element method on polytopal grids (PolyDG \cite{Cangiani2014}). Examples of PolyDG schemes can be found in  \cite{Antonietti2013,Bassi2012} for elliptic problems,  in \cite{Cangiani2017} for parabolic problems, and in   \cite{Antonietti2019, Botti2020_korn, Botti2021} for poroelasticity. Moreover, in \cite{Antonietti2021, Antonietti2022,  DeLaPuente2008} PolyDG methods for wave propagation problems in porous media are analyzed. The PolyDG methods fit well in this framework because they guarantee an arbitrary-order accuracy as well as a high-level of geometric flexibility, that allows to handle complex geometries, solutions with (local) low-regularity, and layered heterogeneous materials.

The major highlights of this paper are: \textit{(i)} a precise model derivation and mathematical formulation of the fully-dynamic thermo-poroelasticity problem; \textit{(ii)} a complete analysis of the PolyDG discretization establishing stability and error estimate in the $hp$-framework; and \textit{(iii)} an in-depth numerical investigation on wave propagation phenomena in thermo-poroelastic media including also a comparison with the poroelastic model which neglects the temperature effects. Through numerical experiments, we also demonstrated the applicability of this model for tests of physical interest in heterogeneous media. Indeed, we take advantage of the geometric flexibility of the PolyDG scheme by allowing all the model coefficients to be discontinuous in the computational domain.

The rest of the paper is organized as follows: in Section~\ref{sec:model_problem} we present the model problem, the assumptions on the model's coefficients, and its weak formulation. In Section~\ref{sec:discretization}, we derive the semi-discrete PolyDG formulation (cf. Section~\ref{sec:DG_formulation}) and then the fully-discrete formulation is obtained by the coupling with the Newmark-$\beta$ time integration scheme ((cf. Section~\ref{sec:fully_discrete}). Section~\ref{sec:convergence_test} and Section~\ref{sec:numerical_results} are devoted to the numerical investigation. In Section~\ref{sec:convergence_test} we validate the convergence of the method on regular manufactured solutions, while in Section~\ref{sec:numerical_results} we present a wide set of test cases aimed at assessing the numerical performances and comparing the results with benchmark configurations taken from the previous literature. Moreover, test cases of geophysical interest are considered. Last, in Section~\ref{sec:semidiscrete_analysis} we present the results about stability and error estimates for the semi-discrete problem. The proofs of these results are postponed to \ref{sec:staberr_proofs}.

\section{Model problem}
\label{sec:model_problem}
The goal of this section is to present the fully-dynamic thermo-poroelastic model and provide its mathematical formulation and physical derivation. We start by considering the linear thermo-poroelastic problem \cite{AntoniettiBonetti2022, Brun2019} with the additional contribution of the inertial terms. 
Let $\Omega\subset\mathbb{R}^d$, $d=2,3$, be an open, convex polygonal/polyhedral domain with Lipschitz boundary $\partial\Omega$.
Given a final time $T_f > 0$, the problem reads: \textit{find $(\mathbf{u}, \mathbf{w}, p, T)$ such that:}
\begin{subequations}
    \label{eq:TPE_system_1}
    \begin{empheq}[left=\empheqlbrace]{align}
    & \rho \ddot{\mathbf{u}} + \rho_f \ddot{\mathbf{w}} -\div{\boldsymbol{\sigma}} = \Tilde{\mathbf{f}}&\text{in }\Omega\times(0,T_f], \label{eq:momentum_cons_1} \\
    & \rho_f \ddot{\mathbf{u}} + \rho_w \ddot{\mathbf{w}} + \mathbf{K}^{-1} \dot{\mathbf{w}} + \nabla p = \Tilde{\mathbf{g}} &\text{in }\Omega\times(0,T_f], \label{eq:filtration_displ_1} \\
    & c_0 \dot{p} - b_0 \dot{T} + \alpha \div{\dot{\mathbf{u}}} + \div{\dot{\mathbf{w}}} = 0
    &\text{in }\Omega\times(0,T_f], 
    \label{eq:mass_cons_1} \\[3pt] 
    & a_0 \left( \dot{T} + \tau \ddot{T} \right) - b_0 \left( \dot{p} + \tau \ddot{p} \right) + \beta \left( \div{\dot{\mathbf{u}}} + \tau \div{\ddot{\mathbf{u}}} \right)  - \div{(\boldsymbol{\Theta} \nabla T )} = H &\text{in }\Omega\times(0,T_f], \label{eq:energy_cons_1}
    \end{empheq}
\end{subequations}
where the four unknowns $(\mathbf{u}, \mathbf{w}, p, T)$ represent the solid displacement, the filtration displacement, the pore pressure, and the temperature, respectively. 
The filtration displacement is a quantity that represents the relative displacement among the fluid and the pore matrix, scaled with respect to the porosity \cite{Matuszyk2014}.
Note that, in Problem~\eqref{eq:TPE_system_1}, the field $T$ expresses the variation of the temperature distribution with respect to a reference value $T_0$. The short-hand notation $\dot{\psi}$ and $\ddot{\psi}$ is used for denoting the first and second partial derivatives with respect to time of a function $\psi:\Omega\times(0,T_f]\to\mathbb{R}$, respectively.

Equations~\eqref{eq:momentum_cons_1}, \eqref{eq:mass_cons_1}, \eqref{eq:energy_cons_1} represent momentum conservation, mass conservation, and energy conservation, respectively. 
Equation \eqref{eq:filtration_displ_1} corresponds to Darcy's law in its dynamic form. 
The terms $\Tilde{\mathbf{f}}$, $\Tilde{\mathbf{g}}$, and H are source terms that represent a body force, a fluid mass source, and a heat source. The description of the model's coefficients appearing in \eqref{eq:TPE_system_1} is given in Section~\ref{sec:tpe_coeff} (cf. Table~\ref{table:TPE_coeff}). From now on, we assume that the hydraulic and thermal conductivities $\mathbf{K}$ and $\boldsymbol{\Theta}$ are isotropic, namely $\mathbf{K} = k \mathbf{I}$ and $\boldsymbol{\Theta} = \theta \mathbf{I}$. 
However, we remark that the general anysotropic case can be handled with minor modifications following the lines of \cite{AntoniettiBonetti2022}. 
The constitutive law for the total stress tensor $\boldsymbol{\sigma}$ is obtained as in \cite{Coussy2003} under the small deformations assumption taking also into account the hydraulic and thermal effects on the porous matrix:
\begin{equation}
    \label{eq:const_law_sigma}
    \boldsymbol{\sigma}(\mathbf{u},p,T) = 2 \mu \boldsymbol{\epsilon}(\mathbf{u}) + \lambda \nabla \cdot \mathbf{u} \mathbf{I} - \alpha p \mathbf{I} - \beta T \mathbf{I}, 
\end{equation}
where $\mathbf{I}$ is the identity tensor and $\boldsymbol{\epsilon}(\mathbf{u})= \frac{1}{2}(\nabla \mathbf{u} + \nabla \mathbf{u}^T)$ is the strain tensor. Finally, problem \eqref{eq:TPE_system_1} is closed by imposing suitable boundary and initial conditions.
\begin{remark}
As pointed out in \cite{Chiavassa2013}, the second equation in Problem~\ref{eq:TPE_system_1} is valid under a constraint on frequencies. 
Namely, the spectrum of the waves has to lie in the low-frequency range. Thus, in what follows, we consider the frequencies to be lower than the critical value $f_c = \phi / (2 \pi a k \rho_f )$.
\end{remark}

\subsection{Three-field formulation}

In the spirit of what is done to obtain the \textit{two-displacements formulation} for the poroelasticity problem (cf. \cite{Matuszyk2014}), we exploit \eqref{eq:mass_cons_1} to express $p$ and its partial derivatives in terms of the other three unknowns of the problem. Hence, we need to recast the expression for $\dot{p}$, $\ddot{p}$, and $\nabla p$ as follows:
\begin{equation}
\label{eq:const_law_press}
\dot{p} = - {c_0}^{-1} \left( \alpha \div{\dot{\mathbf{u}}} + \div{\dot{\mathbf{w}}} - b_0 \dot{T} \right), \quad \nabla p = - {c_0}^{-1} \left( \alpha \nabla (\div{\mathbf{u}}) + \nabla (\div{\mathbf{w}}) - b_0 \nabla T \right) + \Psi_0,
\end{equation}
where $\Psi_0:\Omega\to\mathbb{R}^d$ is the vector field taking into account the initial condition on the fluid content that appears due to the integration in time of \eqref{eq:mass_cons_1}.

Plugging \eqref{eq:const_law_sigma} and \eqref{eq:const_law_press} into problem \eqref{eq:TPE_system_1}, we obtain the following three-field thermo-poroelasticity system: \textit{find $(\mathbf{u}, \mathbf{w}, T):\Omega\times (0,T_f]\to \mathbb{R}^d\times\mathbb{R}^d\times\mathbb{R}$ satisfying}
\begin{subequations}
    \label{eq:TPE_system}
    \begin{empheq}[left=\empheqlbrace]{align}
    & \rho \ddot{\mathbf{u}} + \rho_f \ddot{\mathbf{w}} -
    \div{\bigg(2 \mu \boldsymbol{\epsilon}(\mathbf{u}) + \left( \lambda + \frac{\alpha^2}{c_0} \right) \div{\mathbf{u}} \, \mathbf{I} + \frac{\alpha}{c_0} \div{\mathbf{w}} \, \mathbf{I} - \left(\frac{b_0 \alpha}{c_0} + \beta \right) T \mathbf{I} \bigg)} = \mathbf{f}, \label{eq:momentum_cons} \\
    & \rho_f \ddot{\mathbf{u}} + \rho_w \ddot{\mathbf{w}} + \frac{1}{k} \dot{\mathbf{w}} - \frac{\alpha}{c_0} \nabla (\div{\mathbf{u}}) - \frac{1}{c_0} \nabla (\div{\mathbf{w}}) + \frac{b_0}{c_0} \nabla T = \mathbf{g}, \hspace{3cm} 
    \label{eq:mass_cons} \\[5pt] 
    & \left( a_0 - \frac{b_0^2}{c_0} \right) \left( \dot{T} + \tau \ddot{T} \right) + \left( \frac{b_0 \alpha}{c_0} + \beta \right) \left( \div{\dot{\mathbf{u}}} + \tau \div{\ddot{\mathbf{u}}} \right) + \frac{b_0}{c_0} \left( \div{\dot{\mathbf{w}}} + \tau \div{\ddot{\mathbf{w}}} \right) - \div{(\theta \nabla T )} = H. \label{eq:energy_cons} 
    \end{empheq}
\end{subequations}
Note that, the contribution of $\Psi_0$ in \eqref{eq:momentum_cons} and \eqref{eq:mass_cons} has been included in the forcing terms, that have been redefined as $\mathbf{f} = \Tilde{\mathbf{f}} - \Psi_0(\mathbf{x})$, $\mathbf{g} = \Tilde{\mathbf{g}} - \Psi_0(\mathbf{x})$.
For the sake of simplicity, we complete problem~\ref{eq:TPE_system} by imposing homogeneous Dirichlet boundary conditions on the whole boundary $\partial\Omega$ and by introducing suitable initial conditions, e.g., $$
 (\mathbf{u}, \mathbf{w}, T) (\cdot, t=0) =\left( \mathbf{u}_0, \mathbf{w}_0, T_0\right)
 \quad\text{and}\quad
 (\dot{\mathbf{u}}, \dot{\mathbf{w}}, \dot{T})(\cdot, t=0)=
 \left(\mathbf{u}_1, \mathbf{w}_1, T_1 \right)
 \quad\text{in }\Omega.
$$  
 
We observe that model \eqref{eq:TPE_system} is slightly different than the thermo-poroelastic problem investigated in \cite{Carcione2019}. Indeed, in \eqref{eq:energy_cons} we consider different multiplying factors for the divergence of the solid and filtration displacements. This difference comes from the fact that the conservation of thermal energy assumed as a starting point in \cite{Carcione2019} shows a dependence on the variation of the fluid content and not on the variation of the pore pressure (cf. the term in \eqref{eq:energy_cons_1} weighted by the coefficient $b_0$). Nevertheless, we point out that the two models lead to equivalent partial differential systems under proper choices of the thermo-poroelastic parameters. This is in agreement with the physical relation between the pore pressure and fluid content described in \cite{Coussy2003} and \cite{Carcione2014}. In this paper, we prefer to state the energy conservation as in \eqref{eq:energy_cons} because it reduces to the one considered in \cite{Brun2019, Brun2020, AntoniettiBonetti2022} in the case $\tau=0$.
 
\begin{remark}
The choice of the \textit{two-displacements} formulation for poroelasticity is not the only possibile option, but it turns out to be convenient for the coupling with the temperature. Additionally, it allows to write the problem in a purely second-order hyperbolic form.
\end{remark}

\subsection{On the thermo-poroelastic coefficients}
\label{sec:tpe_coeff}
The coefficients appearing in problems~\eqref{eq:TPE_system_1} and \eqref{eq:TPE_system}, along with their unit of measure and physical meaning are reported in Table~\ref{table:TPE_coeff}.
\begin{table}
    \centering 
    \begin{tabular}{ c | c | l }
    \textbf{Symbol} & \textbf{Unit} & \textbf{Quantity}  \T\B \\
    $a_0$ & \si[per-mode = symbol]{\pascal \per \kelvin \squared} & thermal capacity \T\B\\
    $b_0$ & \si{\per \kelvin} & thermal dilatation coefficient \T\B\\
    $c_0$ & \si{\per\pascal} & specific storage coefficient \T\B\\
    $\alpha$ & - & Biot--Willis constant \T\B\\
    $\beta$ & \si[per-mode = symbol]{\pascal \per \kelvin} & thermal stress coefficient \T\B\\ 
    $\mu, \lambda$ & \si{\pascal} & Lamé parameters  \T\B\\
    $k$ & \si[per-mode = symbol]
    {\metre\squared \per \pascal \per \second} & permeability divided by dynamic fluid viscosity \T\B\\
    $\theta$ & \si[per-mode = symbol]{\metre\squared \pascal \per \kelvin\squared \per \second} & effective thermal conductivity \T\B\\
    $\rho_f$ & \si[per-mode = symbol]
    {\kilogram \per \meter\cubed} & saturating fluid density  \T\B\\
    $\rho_s$ & \si[per-mode = symbol]
    {\kilogram \per \meter\cubed} & solid matrix density  \T\B\\
    $\phi$ & -  & porosity \T\B\\
    $a$ & - & tortuosity \T\B\\
    $\tau$ & \si{\second}  & Maxwell-Vernotte-Cattaneo relaxation time \\
    \end{tabular}
    \\[10pt]
    \caption{Thermo-poroelastic coefficients appearing in problem~\eqref{eq:TPE_system}}
    \label{table:TPE_coeff}
\end{table}
All the model parameters are intended as (possibly) heterogenous scalar fields. The densities $\rho$ and $\rho_w$ are given by
\begin{equation}
    \rho = \phi \rho_f + (1 - \phi) \rho_s > 0, \label{eq:density} \quad
    \rho_w = \frac{a}{\phi}\rho_f > 0, \label{eq:bulk_density}
\end{equation}
where the porosity $\phi$ (the ratio between the void space in a porous medium and its whole volume) and the tortuosity $a$ (the measure of the deviation of the fluid paths from straight streamlines) are such that $0 < \phi_0 \leq \phi \leq \phi_1 < 1$ and $a > 1$ \cite{Souzanchi2013}. 

We refer to \cite{Coussy2003, Carcione2019} for additional comments on the remaining physical coefficients presented in Table~\ref{table:TPE_coeff}. 
Moreover, in \cite{Brun2019, Brun2020, AntoniettiBonetti2022} the assumptions on the model parameters needed to ensure the well-posedness of the quasi-static termo-poroelastic problem are pointed out. However, by passing from the \textit{displacement-pressure} to the \textit{two-displacements} formulation, we need to slightly modify the assumptions on the thermal capacity, thermal dilatation and specific storage coefficients, i.e.
\begin{assumption}
\label{ass:c0}
The model parameters $a_0$, $b_0$, and $c_0$ are such that $a_0\ge b_0^2 c_0^{-1}$, $b_0 \geq 0$, and $c_0 > 0$.
\end{assumption}

The hypotheses on the coefficients multiplying the coupling and elliptic terms are rather standard. We assume that the Biot--Willis modulus and thermal stress coupling parameters satisfy $\phi<\alpha\le 1$ and $\beta>0$, respectively. The shear and dilatation moduli $\mu$ and $\lambda$ together with the conductivities $k$ and $\theta$ are all assumed to be strictly-positive. Finally, the relaxation parameter $\tau$ is a positive scalar field possibly equal to zero.

\subsection{Weak formulation}
\label{subsection:weak_formulation}

Before presenting the variational formulation of problem~\eqref{eq:TPE_system} we introduce the required notation. 
For $X\subseteq\Omega$, we denote by $L^p(X)$ the standard Lebesgue space of index $p\in [1, \infty]$ and by $H^q(X)$ the Sobolev space of index $q \geq 0$ of real-valued functions defined on $X$. 
The notation $\mathbf{L}^2(X)$ and $\mathbf{H}^q(X)$ is adopted in place of $\left[ L^2(X) \right]^d$ and $\left[ H^q(X) \right]^d$, respectively. 
In addition, we denote by $\mathbf{H}(\textrm{div},X)$ the space of $\mathbf{L}^2(X)$ vector fields whose divergence is square integrable. These spaces are equipped with natural inner products and norms denoted by $(\cdot, \cdot)_X = (\cdot, \cdot)_{L^2(X)}$ and $||\cdot||_X = ||\cdot||_{L^2(X)}$, with the convention that the subscript can be omitted in the case $X=\Omega$.

For $T_f>0$ and a Banach space $X$, we denote by $L^p((0,T_f]; X)$ the Bochner space of $L^p$- regular functions defined on $(0,T_f]$ with values in $X$. Being $||\cdot||_X$ a norm in $X$, we then define the norm in $L^p((0,T_f]; X)$ as
\begin{equation}
    ||u||_{L^p((0,T_f];X)} = \left( \int_0^{T_f} ||u(t)||_X^p dt \right)^{\frac{1}{p}}.
\end{equation}
Finally, given $k\in\mathbb{N}$, the usual notation $C^k((0,T];X)$ is used for the space of $X$-valued functions which are $k$-times continuously differentiable in $[0,T]$.
For the sake of brevity, in what follows, we make use of the symbol $x \lesssim y$ to denote $x \le C y$, where $C$ is a positive constant independent of the discretization parameters, but possibly dependent on physical coefficients and final time $T_f$.

To derive the weak formulation of problem \eqref{eq:TPE_system} we start by providing the definition of the functional spaces that take into account the essential boundary conditions, namely
\begin{equation}
\label{eq:func_spaces}
\begin{aligned}
V & = H_0^1(\Omega) = \left\{ \varphi \in H^1(\Omega) \ \text{s.t.} \  \varphi_{\vert \partial \Omega} = 0 \right\}, \\
\mathbf{W} & = \mathbf{H}_0(\textrm{div},X) = \left\{ \mathbf{w} \in \mathbf{H}(\textrm{div},X) \ \text{s.t.} \  \left( \mathbf{w} \mkern-2.5mu \cdot \mkern-2.5mu \mathbf{n} \right)_{\vert \partial \Omega} = 0 \right\}.\\
\end{aligned}
\end{equation}
We use the following notation through the article: $\mathbf{V} = \left[ V \right]^d$. Next, we multiply \eqref{eq:TPE_system} for suitable test functions and we sum up all the contributions to obtain: \textit{for any time $t \in (0, T_f]$, find $(\mathbf{u}, \mathbf{w}, T)(t) \in \mathbf{V} \times \mathbf{W} \times V$ such that\ $ \forall (\mathbf{v}, \mathbf{z}, S) \in \mathbf{V} \times \mathbf{W} \times V$:}
\begin{equation}
\label{eq:DYN_weak_form}
\begin{aligned}
& \mathcal{M}_{uw}((\ddot{\mathbf{u}}, \ddot{\mathbf{w}}), \left( \mathbf{v}, \mathbf{z} \right)) + \tau \mathcal{M}_T(\ddot{T},S) + \tau \ \mathcal{C}(\left(\ddot{\mathbf{u}}, \ddot{\mathbf{w}} \right), S) + \mathcal{B}( \dot{\mathbf{w}}, \mathbf{z}) +  \mathcal{M}_T(\dot{T},S) + \mathcal{C}(\left(\dot{\mathbf{u}}, \dot{\mathbf{w}} \right), S) \\
& + \mathcal{A}_{uw}((\mathbf{u}, \mathbf{w}), \left( \mathbf{v}, \mathbf{z} \right)) + \mathcal{A}_T(T,S) - \mathcal{C}(\left( \mathbf{v}, \mathbf{z} \right), T) = ((\mathbf{f}, \mathbf{g}, H), (\mathbf{v}, \mathbf{z}, S)),
\end{aligned}
\end{equation}
where for any $\left( \mathbf{u}, \mathbf{w}, T \right), \left( \mathbf{v}, \mathbf{z}, S \right) \in \mathbf{V} \times \mathbf{W} \times V$ we have set:
\begin{equation}
\label{eq:bil_form_cont}
\begin{aligned}
    & \mathcal{M}_{uw}(\left(\mathbf{u}, \mathbf{w} \right), \left(\mathbf{v}, \mathbf{z} \right)) = \left( \rho \mathbf{u} + \rho_f \mathbf{w}, \mathbf{v} \right) + \left( \rho_f \mathbf{u} + \rho_w \mathbf{w}, \mathbf{z} \right), \\[5pt]
    & \mathcal{M}_{T}(T,S) =  \left(\left( a_0 - b_0^2/c_0 \right) T, S \right), \\[5pt]
    & \mathcal{A}_{uw}(\left(\mathbf{u}, \mathbf{w} \right), \left(\mathbf{v}, \mathbf{z} \right)) =  (2\mu\, \boldsymbol{\epsilon}(\mathbf{u}), \boldsymbol{\epsilon}(\mathbf{v})) +  (\lambda\div{\mathbf{u}}, \div{\mathbf{v}}) +  ({c_0}^{-1}(\alpha \div{\mathbf{u}} + \div{\mathbf{w}}), \alpha \div{\mathbf{v}} + \div{\mathbf{z}} ), \\[5pt]
    & \mathcal{A}_{T}(T,S) = (\theta \, \nabla T, \nabla S), \\[5pt]
    & \mathcal{B}(\mathbf{w}, \mathbf{z}) = ({k}^{-1}\mathbf{w}, \mathbf{z}), \\[5pt]
    & \mathcal{C}(\left(\mathbf{u}, \mathbf{w} \right), S) = \left( \left[ \alpha b_0/c_0 + \beta \right] \div{\mathbf{u}} + b_0/c_0 \div{\mathbf{w}}, S \right). \\[5pt]
\end{aligned}
\end{equation}

\section{Discretization}
\label{sec:discretization}
The purpose of this section is to derive the fully-discrete scheme for problem~\eqref{eq:TPE_system}. We start by introducing the spatial discretization obtained via the PolyDG approximation, cf. Section~\ref{sec:DG_formulation}, and then we couple it with the implicit Newmark-$\beta$ time integration scheme, cf. Section~\ref{sec:fully_discrete}.

\subsection{Preliminaries}

First, we present the mesh assumptions, the discrete spaces, and some instrumental results for the design and analysis of PolyDG schemes.
We introduce a subdivision $\mathcal{T}_h$ of the computational domain $\Omega$, whose elements are polygons/polyhedrons in dimension $d = 2, 3$, respectively. Next, we define the interfaces (or internal faces) as subsets of the intersection of any two neighbouring elements of $\mathcal{T}_h$. If $d = 2$ an interface is a line segment, while if $d = 3$ an interface is a planar polygon, that we assume can be further decomposed into a set of triangles. The same holds for the boundary faces collected in the set $\mathcal{F}_B$ which yields a simplicial subdivision of $\partial\Omega$.  Accordingly, we define $\mathcal{F}_I$ to be the set of internal faces and the set of all the faces as $\mathcal{F}_h=\mathcal{F}_B\cup\mathcal{F}_I$. In what follows, we introduce the main assumptions on the mesh $\mathcal{T}_h$ (cf. \cite{CangianiDong:17, Cangiani2014}).
\begin{definition}[Polytopic-regular mesh]
\label{def:unif_regular}
A mesh $\mathcal{T}_h$ is polytopic-regular if for any $ \kappa \in \mathcal{T}_h$, there exist a set of non-overlapping simplices contained in $\kappa$, denoted by $\{S_{\kappa}^F\}_{F \subset \partial \kappa}$, such that, for any face $F \subset \partial \kappa$, the following condition holds: $h_{\kappa} \lesssim d \ |S_{\kappa}^F| \ |F|^{-1}$, with $h_{\kappa}$ denoting the diameter of the element $\kappa$ and with $| \cdot |$ denoting the Hausdorff measure.
\end{definition}

As a basis for the construction of the PolyDG approximation, we define fully-discontinuous polynomial spaces on the mesh $\mathcal{T}_h$. 
Given an element-wise constant polynomial degree $\ell:\mathcal{T}_h\to\mathbb{N}_{>0}
$ which determines the order of the approximation, the discrete spaces are defined such as
\begin{equation}
    \label{eq:discrete_spaces}
    \begin{aligned}
    V_h^{\ell} &= \left\{ v_h \in L^2(\Omega) : v_h |_{\kappa} \in \mathbb{P}^{\ell_{\kappa}}(\kappa) \ \ \forall \kappa \in \mathcal{T}_h \right\}, \quad \mathbf{V}_h^{\ell} = \left[ V_h^{\ell} \right]^d.\\
    \end{aligned}
\end{equation}
where, for each $\kappa\in\mathcal{T}_h$, the space $\mathbb{P}^{\ell_{\kappa}}(\kappa)$ is spanned by polynomials of maximum degree  $\ell_{\kappa}=\ell_{|\kappa}$. In order to analyze the convergence of the spatial discretization, we consider a mesh sequence $\{\mathcal{T}_h\}_{h\to0}$ satisfying the following properties:

\begin{assumption}
\label{ass:mesh_Th1}
The mesh sequence $\{\mathcal{T}_h\}_{h\to0}$ and the polynomial degree $\ell$ are such that
\begin{enumerate}[start=1,label={\bfseries A.\arabic*}]
    \item \label{ass:A1} $\{\mathcal{T}_h\}_{h\to0}$ is uniformly polytopic-regular;
    \item \label{ass:A2} For each $\mathcal{T}_h\in \{\mathcal{T}_h\}_{h\to0}$ there exists a shape-regular, simplicial covering $\mathcal{T}_h^*$ of $\mathcal{T}_h$ such that, for each pair $\kappa \in \mathcal{T}_h$ and $k \in \mathcal{T}_h^*$ with $\kappa \subset k$ it holds
    \begin{enumerate}[start=1,label={(\roman*)}]
    \item $h_{k} \lesssim h_{\kappa}$;
    \item $\underset{\kappa \in \mathcal{T}_h}{\mbox{max}} \ \mbox{card} \left\{\kappa' \in \mathcal{T}_h: \kappa' \cap k \neq 0, k \in \mathcal{T}_h^*, \kappa \subset k \right\} \lesssim 1$;
    \end{enumerate}
    \item \label{ass:A3} For each $\mathcal{T}_h\in \{\mathcal{T}_h\}_{h\to0}$ and for any pair of neighbouring elements $\kappa^+, \kappa^- \in \mathcal{T}_h$, the following $hp$-local bounded variation properties hold: $h_{\kappa^+} \lesssim h_{\kappa^-} \lesssim h_{\kappa^+}$ and $\ell_{\kappa^+} \lesssim \ell_{\kappa^-} \lesssim \ell_{\kappa^+}$.
\end{enumerate}
\end{assumption}
We remark that under \textit{\ref{ass:A1}} the following inequality (called \textit{discrete trace-inverse inequality}) holds (cf. \cite{Cangiani2017} for all the details):
\begin{equation}
    \label{eq:trace_inverse_ineq}
    ||v||_{L^2(\partial \kappa)} \lesssim \frac{\ell_{\kappa}}{h_{\kappa}^{1/2}} ||v||_{L^2(\kappa)} \quad \forall v \in \mathbb{P}^{\ell_{\kappa}}(\kappa),
\end{equation}
where the hidden constant is independent of $\ell_{\kappa}, h_{\kappa}$, and the number of faces per element. For deriving the discontinuous Galerkin formulation, we also need to introduce the average and jump operators. We start by defining them on each interface $F\in\mathcal{F}_I$ shared by the elements $\kappa^{\pm}$ as in \cite{Arnold2002}:
\begin{equation}
    \label{eq:avg_jump_operators}
    \begin{aligned}
    & \jump{a} = a^+ \mathbf{n^+} + a^- \mathbf{n^-}, \ 
    && \jump{\mathbf{a}} = \mathbf{a}^+ \odot \mathbf{n^+} + \mathbf{a}^- \odot \mathbf{n^-}, \ 
    &&\jump{\mathbf{a}}_n = \mathbf{a}^+ \cdot \mathbf{n^+} + \mathbf{a}^- \cdot \mathbf{n^-}, \\ 
    & \avg{a} = \frac{a^+ + a^-}{2}, \
    && \avg{\mathbf{a}} = \frac{\mathbf{a}^+ + \mathbf{a}^-}{2}, \ && \avg{\mathbf{A}} = \frac{\mathbf{A}^+ + \mathbf{A}^-}{2},
    \end{aligned}
\end{equation}
where $\mathbf{a} \odot \mathbf{n} = \mathbf{a}\mathbf{n}^T$, and $a, \ \mathbf{a}, \ \mathbf{A}$ are (regular enough) scalar-, vector-, and tensor-valued functions, respectively. The notation $(\cdot)^{\pm}$ is used for the trace on $F$ taken within the interior of $\kappa^\pm$ and $\mathbf{n}^\pm$ is the outer unit normal vector to $\partial \kappa^\pm$. Accordingly, on boundary faces $F\in\mathcal{F}_B$, we set
$$
 \jump{a} = a \mathbf{n},\ \
\avg{a} = a,\ \
\jump{\mathbf{a}} = \mathbf{a} \odot \mathbf{n},\ \
\avg{\mathbf{a}} = \mathbf{a},\ \
\jump{\mathbf{a}}_n = \mathbf{a} \cdot \mathbf{n},\ \
\avg{\mathbf{A}} = \mathbf{A}.
$$ 
From now on, for the sake of simplicity, we assume that the model parameters are element-wise constant. Moreover, for later use, we can introduce the quantities
\begin{equation}
     c_{0,{\kappa}} = c_0 \rvert_{\kappa}, \ \ \theta_{\kappa} = \theta \rvert_{\kappa}, \ \ \lambda_{\kappa} = \lambda \rvert_{\kappa},
    \ \ \text{and } \mu_{\kappa} = \mu \rvert_{\kappa}.
\end{equation}

\subsection{Discontinuous Galerkin semi-discrete problem}
\label{sec:DG_formulation}

We are now ready to derive the semi-discrete PolyDG approximation of the fully-dynamic thermo-poroelastic problem~\eqref{eq:TPE_system}.
In the following discussion, we choose the Interior Penalty formulation \cite{Wheeler1978, Arnold1982, Ern2021}. Thus, the PolyDG semi-discretization of problem \eqref{eq:DYN_weak_form} reads:
\newline
\textit{for any $t \in (0, T_f]$, find $(\mathbf{u}_h, \mathbf{w}_h, T_h)(t) \in \mathbf{V}^{\ell}_h \times \mathbf{V}^{\ell}_h \times V^{\ell}_h $ such that:}
\begin{equation}
\label{eq:semidiscr_form_tau}
\begin{aligned}
& \mathcal{M}_{uw}((\ddot{\mathbf{u}}_h, \ddot{\mathbf{w}}_h), \left( \mathbf{v}_h, \mathbf{z}_h \right)) + \tau \mathcal{M}_T(\ddot{T}_h,S_h) + \tau \ \mathcal{C}_h(\left(\ddot{\mathbf{u}}_h, \ddot{\mathbf{v}}_h \right), S_h) + \mathcal{B}( \dot{\mathbf{w}}_h, \mathbf{z}_h) +  \mathcal{M}_T(\dot{T}_h,S_h) \\
& + \mathcal{C}_h(\left(\dot{\mathbf{u}}_h, \dot{\mathbf{v}}_h \right), S_h) + \mathcal{A}_{uw,h}((\mathbf{u}_h, \mathbf{w}_h), \left( \mathbf{v}_h, \mathbf{z}_h \right)) + \mathcal{A}_{T,h}(T_h,S_h) - \mathcal{C}_h(\left( \mathbf{v}_h, \mathbf{w}_h \right), T_h) \\
& \qquad\qquad\qquad\qquad\qquad = ((\mathbf{f}, \mathbf{g}, H), (\mathbf{v}_h, \mathbf{z}_h, S_h)) \qquad\qquad \forall (\mathbf{v}_h, \mathbf{z}_h, S_h) \in \mathbf{V}^{\ell}_h \times \mathbf{V}^{\ell}_h \times V^{\ell}_h,
\end{aligned}
\end{equation}
\textit{supplemented by initial conditions $(\mathbf{u}_{h,0}, \mathbf{w}_{h,0}, T_{h,0},\dot{\mathbf{u}}_{h,0}, \dot{\mathbf{w}}_{h,0}, \dot{T}_{h,0})$ that are fitting approximations of the initial conditions of problem~\eqref{eq:TPE_system}}. 
The bilinear forms labelled with the subscript ${}_h$  appearing in \eqref{eq:semidiscr_form_tau} are given by
\begin{equation}
\label{eq:bil_form_discr}
\begin{aligned}
    & \mathcal{A}_{uw,h}((\mathbf{u}, \mathbf{w}), ( \mathbf{v}, \mathbf{z})) = \mathcal{A}_{e,h}(\mathbf{u}, \mathbf{v}) + \mathcal{A}_{p,h}(\alpha \mathbf{u} + \mathbf{w}, \alpha \mathbf{v} + \mathbf{z}) \\[5pt]
    & \mathcal{A}_{T,h}(T,S) = \left(\theta \, \nabla_h T, \nabla_h S\right) - \sum_{F \in \mathcal{F}_h} \int_F \big(\avg{\theta \, \nabla_h T} \mkern-2.5mu \cdot \mkern-2.5mu \jump{S} + \jump{T} \mkern-2.5mu \cdot \mkern-2.5mu \avg{\theta \, \nabla_h S} - \varrho \jump{T} \mkern-2.5mu \cdot \mkern-2.5mu \jump{S} \big), \\[5pt] 
    & \mathcal{C}_h(\left(\mathbf{u}, \mathbf{w} \right), S) = \left( \frac{\alpha b_0+ \beta c_0}{c_0} \divh{\mathbf{u}} + \frac{b_0}{c_0} \divh{\mathbf{w}}, S \right) - \hspace{-1mm}
    \sum_{F \in \mathcal{F}_h} \int_F \hspace{-1mm}\left(\avgg{\frac{\alpha b_0+ \beta c_0}{c_0} S} \jump{\mathbf{u}}_n + \avgg{\frac{b_0}{c_0} S} \jump{\mathbf{w}}_n \right)
\end{aligned}
\end{equation}
with
\begin{equation}
\label{eq:bil_form_discr_detailed}
\begin{aligned}
    \mathcal{A}_{e,h}(\mathbf{u}, \mathbf{v}) = \ & (2 \mu\boldsymbol{\epsilon}_h(\mathbf{u}),\boldsymbol{\epsilon}_h(\mathbf{v}))
    - \sum_{F \in \mathcal{F}_h} \int_F \big( \avg{2 \mu\boldsymbol{\epsilon}_h(\mathbf{u})} \mkern-2.5mu : \mkern-2.5mu \jump{\mathbf{v}} + \jump{\mathbf{u}} \mkern-2.5mu : \mkern-2.5mu \avg{2 \mu\boldsymbol{\epsilon}_h(\mathbf{v})} - \sigma \jump{\mathbf{u}} \mkern-2.5mu : \mkern-2.5mu \jump{\mathbf{v}} \big)\\
    & + (\lambda \divh{\mathbf{u}},  \divh{\mathbf{v}}) -  \sum_{F \in \mathcal{F}_h} \int_F \big(\avg{\lambda \divh{\mathbf{u}}} \jump{\mathbf{v}}_\mathbf{n} + \jump{\mathbf{u}}_\mathbf{n} \avg{ \lambda\divh{\mathbf{v}}} - \xi \jump{\mathbf{u}}_\mathbf{n} \jump{\mathbf{v}}_\mathbf{n} \big), \\[5pt]
    \mathcal{A}_{p,h}( \mathbf{w}, \mathbf{z}) = \ & (c_0^{-1} \divh{\mathbf{w}},  \divh{\mathbf{z}}) - \sum_{F \in \mathcal{F}_h} \int_F \big(\avg{c_0^{-1} \divh{\mathbf{w}}} \jump{\mathbf{z}}_\mathbf{n} + \jump{\mathbf{w}}_\mathbf{n} \avg{c_0^{-1}  \divh{\mathbf{z}}} - \zeta \jump{\mathbf{w}}_\mathbf{n} \jump{\mathbf{z}}_\mathbf{n}\big).
\end{aligned}
\end{equation}
Here, for all $a \in V_h^{\ell}$ and $\mathbf{a}\in \mathbf{V}_h^{\ell}$, $\nabla_h a$ and $\divh{\mathbf{a}}$ denote the broken differential operators whose restrictions to each element $\kappa \in \mathcal{T}_h$ are defined as $\nabla w_{|\kappa}$ and $\nabla \cdot w_{|\kappa}$, respectively. Then, the broken version of the strain tensor is defined as $\boldsymbol{\epsilon}_h(\mathbf{u}) = \left(\nabla_h \mathbf{u} + \nabla_h \mathbf{u}^T\right)/2$. Last, we are left to define the stabilization functions $\sigma, \xi, \zeta$ and  $\varrho \in L^{\infty}(\mathcal{F}_h)$. Following \cite{Cangiani2017} we select
\begin{equation}
    \label{eq:stabilization_func}
    \begin{aligned}
    \sigma &= \left\{\begin{aligned}
    &\alpha_1 \underset{\kappa \in \{\kappa^+,\kappa^-\}}{\mbox{max}}\left( \frac{\mu_{\kappa} \ell_{\kappa}^2}{h_{\kappa}}\right) \ & F \in \mathcal{F}_I,\\
    &\alpha_1 \mu_{\kappa} \ell_\kappa^2 h_{\kappa}^{-1} \ & F \in \mathcal{F}_B,\\
    \end{aligned}
    \right.
    \ \
     \xi = &&\left\{\begin{aligned}
    &\alpha_2 \underset{\kappa \in \{\kappa^+,\kappa^-\}}{\mbox{max}}\left( \frac{\lambda_{\kappa} \ell_\kappa^2}{h_{\kappa}}\right) \ & F \in \mathcal{F}_I,\\
    &\alpha_2 \lambda_{\kappa} \ell_\kappa^2 h_{\kappa}^{-1} \ & F \in \mathcal{F}_B,\\
    \end{aligned}
    \right.\\
    \zeta &= \left\{\begin{aligned}
    &\alpha_3 \underset{\kappa \in \{\kappa^+,\kappa^-\}}{\mbox{max}}\left( \frac{\ell_\kappa^2}{c_{0,\kappa} h_{\kappa}}\right) \ & F \in \mathcal{F}_I,\\
    &\alpha_3 {c_{0,\kappa}}^{-1} \ell_\kappa^2 h_{\kappa}^{-1} \ & F \in \mathcal{F}_B,\\
    \end{aligned}
    \right.
    \ \
    \varrho = && \left\{\begin{aligned}
    &\alpha_4 \underset{\kappa \in \{\kappa^+,\kappa^-\}}{\mbox{max}}\left( \frac{\theta_{\kappa} \ell_\kappa^2}{h_{\kappa}}\right) \ & F \in \mathcal{F}_I,\\
    &\alpha_4 \theta_{\kappa} \ell_\kappa^2 h_{\kappa}^{-1} \ & F \in \mathcal{F}_B,\\
    \end{aligned}
    \right.\\
    \end{aligned}   
\end{equation}
where $\alpha_1, \alpha_2, \alpha_3$ and  $\alpha_4 \in \mathbb{R}$ are positive constants to be properly defined.

\subsection{Fully-discrete scheme}
\label{sec:fully_discrete}

By fixing a basis for $\mathbf{V}_h^{\ell}, V_h^{\ell}$ and denoting by $\left[ \mathbf{U}, \mathbf{W}, \mathbf{T} \right]^T$ the vector of the expansion coefficients of the variables $(\mathbf{u}_h, \mathbf{v}_h, T_h)$, we can rewrite the semi-discrete problem \eqref{eq:semidiscr_form_tau} in the equivalent form
\begin{equation}
\label{eq:semidiscr_alg_sys}
\begin{aligned}
\begin{bmatrix}
    \rho \mathbf{M}_{uw} & \rho_f \mathbf{M}_{uw} & 0 \\
    \rho_f \mathbf{M}_{uw} & \rho_w \mathbf{M}_{uw} & 0 \\
    \tau \ \mathbf{C} & \tau \ \mathbf{C} & \tau\left(a_0 - b_0^2/c_0 \right) \mathbf{M}_T
\end{bmatrix}
\begin{bmatrix}
    \ddot{\mathbf{U}} \\
    \ddot{\mathbf{W}} \\
    \ddot{\mathbf{T}} 
\end{bmatrix} &+
\begin{bmatrix}
    0 & 0 & 0 \\
    0 & \mathbf{B} & 0 \\
    \mathbf{C} & \mathbf{C} & \mathbf{M}_T 
\end{bmatrix}
\begin{bmatrix}
    \dot{\mathbf{U}} \\
    \dot{\mathbf{W}} \\
    \dot{\mathbf{T}} 
\end{bmatrix} \\
& +
\begin{bmatrix}
    \mathbf{A}_e + \alpha^2 \mathbf{A}_p  & \alpha \mathbf{A}_p & - \mathbf{C}^T \\
    \alpha \mathbf{A}_p & \mathbf{A}_p & - \mathbf{C}^T \\
    0 & 0 & \mathbf{A}_T 
\end{bmatrix}
\begin{bmatrix}
    \mathbf{U} \\
    \mathbf{W} \\
    \mathbf{T} 
\end{bmatrix} =
\begin{bmatrix}
    \mathbf{F} \\
    \mathbf{G} \\
    \mathbf{H} 
\end{bmatrix}
\end{aligned}
\end{equation}
with initial conditions $\mathbf{U}(0) = \mathbf{U}_0, \mathbf{W}(0) = \mathbf{W}_0, T(0) = T_0, \dot{\mathbf{U}}(0) = \mathbf{U}_1, \dot{\mathbf{W}}(0) = \mathbf{W}_1, \dot{T}(0) = T_1$. The vectors $\mathbf{F}, \mathbf{G}, \mathbf{H}$ are representations of the linear functionals appearing in the right-hand side of \eqref{eq:semidiscr_form_tau}.

To integrate \eqref{eq:semidiscr_alg_sys} in time, we introduce a time-step $\Delta t = T_f/n$, with $n\in\mathbb{N}_{>0}$, discretize the interval $(0, T_f]$ as a sequence of time instants $\{ t_k \}_{0\le k\le n}$ such that $t_{k+1} - t_k = \Delta t$, and define $\mathbf{X}^k = \mathbf{X}(t^k)$, with $\mathbf{X} = \left[\mathbf{U}, \mathbf{W}, \mathbf{T} \right]^T$. Next, we rewrite \eqref{eq:semidiscr_alg_sys} in a compact form as $\mathbf{A} \ddot{\mathbf{X}} + \mathbf{B} \dot{\mathbf{X}} + \mathbf{C} \mathbf{X} = \mathbf{F}$ and derive
\begin{equation}
    \label{eq:def_L_timeint}
    \ddot{\mathbf{X}} = \mathbf{A}^{-1} \left( \mathbf{F} - \mathbf{B} \dot{\mathbf{X}} - \mathbf{C} \mathbf{X} \right) = \mathbf{A}^{-1} \mathbf{F} - \mathbf{A}^{-1} \mathbf{B} \dot{\mathbf{X}} - \mathbf{A}^{-1}\mathbf{C} \mathbf{X} = \mathcal{L}(t, \mathbf{X}, \dot{\mathbf{X}}).
\end{equation}
Last, we integrate in time \eqref{eq:def_L_timeint} with the use of Newmark-$\beta$ scheme, that exploit a Taylor expansion for $\mathbf{X}$ and $\mathbf{Y} = \dot{\mathbf{X}}$:
\begin{equation}
    \label{eq:newmark}
    \left\{
    \begin{aligned}
        & \mathbf{X}^{k+1} = \mathbf{X}^{k} + \Delta t \mathbf{Y}^{k} + \Delta t^2 \left( \beta_N \mathcal{L}^{k+1} + (\frac{1}{2} - \beta_N) \mathcal{L}^{k} \right), \\
        & \mathbf{Y}^{k+1} = \mathbf{Y}^{k} + \Delta t \left( \gamma_N \mathcal{L}^{k+1} + (1 - \gamma_N) \mathcal{L}^{k} \right),
    \end{aligned}
    \right.
\end{equation}
where $\mathcal{L}^k = \mathcal{L}(t^k, \mathbf{X}^k, \dot{\mathbf{X}}^k)$ and the Newmark parameters $\beta_N, \gamma_N$ satisfy: $0 \leq 2 \beta_N \leq 1$, $0 \leq \gamma_N \leq 1$. The typical choices for the Newmark parameters, that ensure unconditionally stability and second-order accuracy for the scheme, are $\beta_N = 1/4$ and $\gamma_N = 1/2$. These are the values used in all the numerical tests of the next two sections.

\section{Convergence tests}
\label{sec:convergence_test}
The aim of this section is to assess the performance of the proposed scheme in terms of accuracy. The results found in this section can be compared with the theoretical analysis for the semi-discrete formulation reported in Section~\ref{sec:semidiscrete_analysis}.

\begin{table}[ht]
\begin{minipage}[c]{0.4\linewidth}
    \centering
    \includegraphics[width=0.9\linewidth]{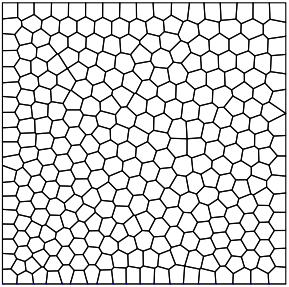}
    \captionof{figure}{Convergence test: example of a $2D$ Voronoi polygonal mesh made of 300 elements.}
    \label{fig:voronoi_mesh}
\end{minipage}\hfill
\begin{minipage}[c]{0.56\linewidth}
    \centering 
    \small
     \begin{tabular}{ l | l  c l | l}
    \textbf{Coefficient} & \textbf{Value} & & \textbf{Coefficient} & \textbf{Value} \T\B \\
    $a_0 \ [\si[per-mode = symbol]{\giga\pascal \per \kelvin\squared}]$ & 0.02 & & $k \ [\si{\dm \squared \per \giga\pascal \per \hour}]$ & 0.2  \T\B \\
    $b_0 \ [\si{\per \kelvin}]$ & 0.01 & & $\theta \ [\si{\dm \squared \giga\pascal \per \kelvin\squared \per \hour}]$ & 0.05 \T\B \\
    $c_0 \ [\si{\per \giga\pascal}]$ & 0.03 & &  $\rho_f \ [\si{\kilogram \per \meter\cubed}]$ & 0.03  \T\B \\
    $\alpha \ [-]$ & 1 & & $\rho_s \ [\si{\kilogram \per \meter\cubed}]$ & 0.03 \T\B \\
    $\beta \ [\si{\giga\pascal \per \kelvin}]$ & 0.8 & & $\phi \ [-]$ & 0.5 \T\B \\
    $\mu \ [\si{\giga\pascal}]$ & 1 & & $a \ [-]$ & 1 \T\B \\
    $\lambda \ [\si{\giga\pascal}]$ & 5 & & $\tau [\si{\second}]$ & 0.01 \T\B \\
    \end{tabular}
\caption{Convergence test: problem's parameters for the convergence analysis}
\label{tab:params_convtest}
\end{minipage}
\end{table}
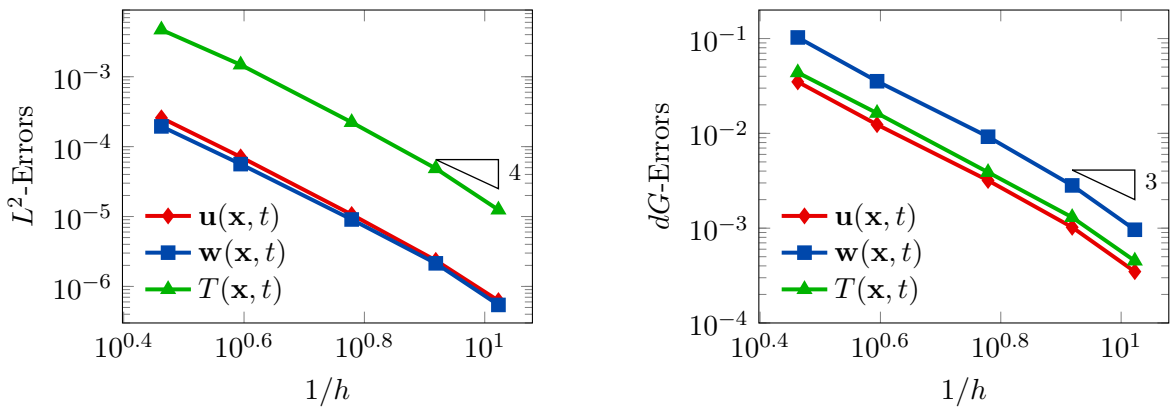
\begin{figure}[ht]
\begin{subfigure}{1\textwidth}
\centering
\begin{subfigure}[b]{0.49\textwidth}
\begin{tikzpicture}
\begin{axis}[%
width=0.65\textwidth,
height=0.5\textwidth,
at={(0\textwidth,0\textwidth)},
scale only axis,
xmode=log,
xmin=2.5,
xmax=12,
xminorticks=true,
xlabel={$1/h$},
ymode=log,
ymin=3e-07,
ymax=0.009,
yminorticks=true,
ylabel={$L^2$-Errors},
legend style={draw=none,fill=none,legend cell align=left},
legend pos=south west
]
\addplot [color=myred,solid,line width=1.5pt, mark=diamond*,mark options={color=myred}]
  table[row sep=crcr]{
    2.9017   0.00025795 \\
    3.9303   7.1154e-05 \\
    6.0098   1.0671e-05 \\
    8.2911   2.3553e-06 \\
    10.5366   6.258e-07 \\
};
\addlegendentry{$\mathbf{u}(\mathbf{x},t)$}

\addplot [color=myblue,solid,line width=1.5pt,mark=square*,mark options={color=myblue}]
  table[row sep=crcr]{
    2.9017   0.00019466 \\
    3.9303   5.6195e-05 \\
    6.0098   9.069e-06 \\
    8.2911   2.1399e-06 \\
    10.5366   5.4361e-07 \\
};
\addlegendentry{$\mathbf{w}(\mathbf{x},t)$}

\addplot [color=mygreen,solid,line width=1.5pt,mark=triangle*,mark options={color=mygreen}]
  table[row sep=crcr]{
    2.9017   0.004696 \\
    3.9303   0.0014816 \\
    6.0098   0.00022264 \\
    8.2911   4.864e-05 \\
    10.5366   1.2459e-05 \\
};
\addlegendentry{$T(\mathbf{x},t)$}

\addplot [color=black,solid,line width=0.5pt]
  table[row sep=crcr]{
 8.2911       6.5207e-05 \\
 10.5366      6.5207e-05 \\
 10.5366      2.5e-05 \\
 8.2911       6.5207e-05 \\ 
};

\node[right, align=left, text=black, font=\footnotesize]
at (axis cs:10.5366,4.2603e-05) {4}; 

\end{axis}
\end{tikzpicture}
\end{subfigure}
\begin{subfigure}[b]{0.49\textwidth}
\begin{tikzpicture}
\begin{axis}[%
width=0.65\textwidth,
height=0.5\textwidth,
at={(0\textwidth,0\textwidth)},
scale only axis,
xmode=log,
xmin=2.5,
xmax=12,
xminorticks=true,
xlabel={$1/h$},
ymode=log,
ymin=0.0001,
ymax=0.2,
yminorticks=true,
ylabel={$dG$-Errors},
legend style={draw=none,fill=none,legend cell align=left},
legend pos=south west
]
\addplot [color=myred,solid,line width=1.5pt, mark=diamond*,mark options={color=myred}]
  table[row sep=crcr]{
    2.9017   0.034954 \\
    3.9303   0.012368 \\
    6.0098   0.0031879 \\
    8.2911   0.0010222 \\
    10.5366   0.00034814 \\
};
\addlegendentry{$\mathbf{u}(\mathbf{x},t)$}

\addplot [color=myblue,solid,line width=1.5pt,mark=square*,mark options={color=myblue}]
  table[row sep=crcr]{
    2.9017   0.10273 \\
    3.9303   0.035459 \\
    6.0098   0.0092184 \\
    8.2911   0.0028243 \\
    10.5366   0.00096237 \\
};
\addlegendentry{$\mathbf{w}(\mathbf{x},t)$}

\addplot [color=mygreen,solid,line width=1.5pt,mark=triangle*,mark options={color=mygreen}]
  table[row sep=crcr]{
    2.9017   0.043787 \\
    3.9303   0.016338 \\
    6.0098   0.0038858 \\
    8.2911   0.0013011 \\
    10.5366   0.0004519 \\
};
\addlegendentry{$T(\mathbf{x},t)$}

\addplot [color=black,solid,line width=0.5pt]
  table[row sep=crcr]{
 8.2911      0.0041 \\
 10.5366     0.0041 \\
 10.5366     0.002 \\
 8.2911      0.0041 \\ 
};
\node[right, align=left, text=black, font=\footnotesize]
at (axis cs:10.5366,0.0031) {3}; 

\end{axis}
\end{tikzpicture}
\end{subfigure}
\end{subfigure}
\caption{Convergence test for $\tau = 0$ : computed errors in $L^2$-norm (left) and $dG$-norm (right) versus $1/h$ (\textit{log-log} scale). The errors are computed at the final time $T_f$. The polynomial degree of approximation is taken as $\ell = 3$.}
\label{fig:convH_tau0}
\end{figure}

The numerical implementation is carried out in MATLAB and the Voronoi meshes are generated via the \texttt{Polymesher} algorithm \cite{Talischi2012}. In all the tests the PolyDG space discretization 
is coupled with the Newmark-$\beta$ time-integration scheme (cf. Section~\ref{sec:fully_discrete}) with parameters $\gamma_N = 1/2$ and $\beta_N = 1/4$ \cite{Antonietti2021}.
In the forthcoming tests, we consider the Maxwell-Vernotte-Cattaneo relaxation time to be either null or positive. 
In Section~\ref{sec:semidiscrete_analysis} the analysis for the semi-discrete problem is carried out for the case $\tau = 0$. However, in this section, we demonstrate numerically that the results can be extended to the case $\tau \neq 0$.

We consider problem \eqref{eq:TPE_system} in the square domain $\Omega = (0,1)^2$ with manufactured analytical solutions 
\begin{equation}
    \begin{aligned}
    \mathbf{u}(x,y,t) & \ = \left( \begin{aligned}
    & x^2 \cos\left(\frac{\pi x}{2}\right) \sin(\pi x) \\
    & x^2 \cos\left(\frac{\pi x}{2}\right) \sin(\pi x)
    \end{aligned} \right) \cos(\sqrt{2} \pi t), \\
    \mathbf{w}(x,y,t) & \ = - \mathbf{u}(x,y,t), \\
    T(x,y,t) & \ = \left( x^2 \sin(\pi x) \sin(\pi y) \right) \sin(\sqrt{2} \pi t).
    \end{aligned}
\end{equation}
The initial conditions, boundary conditions, and forcing terms are inferred from the exact solutions. The model coefficients are chosen as reported in Table~\ref{tab:params_convtest} and follow from a combination of the convergence parameters for the quasi-static thermo-poroelastic problem considered in \cite{AntoniettiBonetti2022} and the two-displacements poroelasticity of \cite{Antonietti2021}.
In the first convergence test, we consider $\tau  = 0$. In this case, the third equation in \eqref{eq:semidiscr_alg_sys} reduces to a first-order differential equation for the temperature $T$. Then, to integrate in time, we consider a Newmark-$\beta$ scheme for the mass and momentum conservation equations, coupled with the Crank-Nicolson method for the energy conservation equation.
The time discretization parameters are $T_f = 0.1$, $\Delta t = \num[exponent-product=\ensuremath{\cdot}]{1e-4}$, and all the penalty coefficients $\alpha_i$, $i=1,...,4$ in \eqref{eq:stabilization_func} are set equal to $10$.

The convergence of the PolyDG-scheme is tested both with respect to the mesh size $h$ and to the polynomial approximation degree $\ell$. For the $h$-convergence a sequence of polygonal meshes as the one in Figure~\ref{fig:voronoi_mesh} is considered, while for the $\ell$-convergence we fixed a computational mesh of $100$ elements and varied the polynomial degree $\ell = 1,2,\dots,5$. 

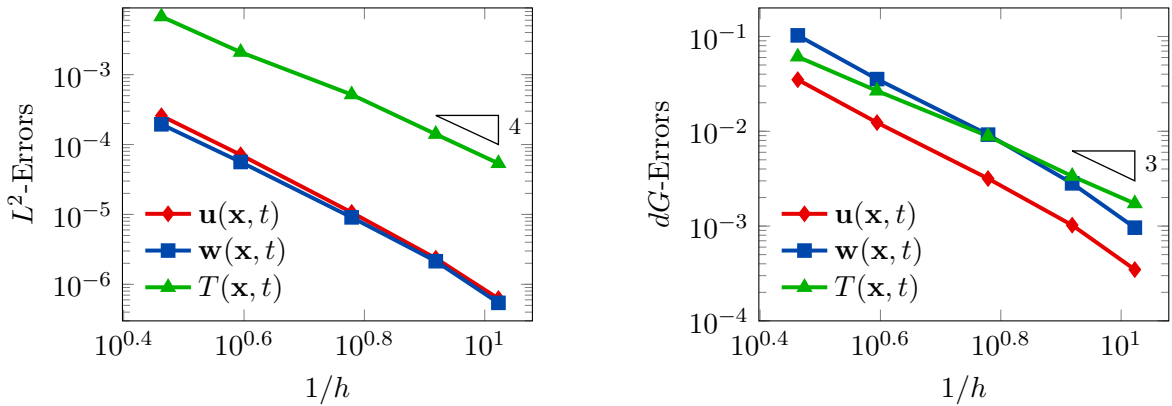
\begin{figure}[ht]
\begin{subfigure}{1\textwidth}
\centering
\begin{subfigure}[b]{0.49\textwidth}
\begin{tikzpicture}
\begin{axis}[%
width=0.65\textwidth,
height=0.5\textwidth,
at={(0\textwidth,0\textwidth)},
scale only axis,
xmode=log,
xmin=2.5,
xmax=12,
xminorticks=true,
xlabel={$1/h$},
ymode=log,
ymin=3e-07,
ymax=0.009,
yminorticks=true,
ylabel={$L^2$-Errors},
legend style={draw=none,fill=none,legend cell align=left},
legend pos=south west
]
\addplot [color=myred,solid,line width=1.5pt, mark=diamond*,mark options={color=myred}]
  table[row sep=crcr]{
    2.9017   0.00025789 \\
    3.9303   7.1161e-05 \\
    6.0098   1.0666e-05 \\
    8.2911   2.3539e-06 \\
    10.5366   6.2509e-07 \\
};
\addlegendentry{$\mathbf{u}(\mathbf{x},t)$}

\addplot [color=myblue,solid,line width=1.5pt,mark=square*,mark options={color=myblue}]
  table[row sep=crcr]{
    2.9017   0.00019458 \\
    3.9303   5.6217e-05 \\
    6.0098   9.0737e-06 \\
    8.2911   2.1407e-06 \\
    10.5366   5.4331e-07 \\
};
\addlegendentry{$\mathbf{w}(\mathbf{x},t)$}

\addplot [color=mygreen,solid,line width=1.5pt,mark=triangle*,mark options={color=mygreen}]
  table[row sep=crcr]{
    2.9017   0.0067536 \\
    3.9303   0.0020942 \\
    6.0098   0.0005178 \\
    8.2911   0.00014052 \\
    10.5366   5.3701e-05 \\
};
\addlegendentry{$T(\mathbf{x},t)$}

\addplot [color=black,solid,line width=0.5pt]
  table[row sep=crcr]{
 8.2911       2.6083e-04 \\
 10.5366      2.6083e-04 \\
 10.5366      1e-04 \\
 8.2911       2.6083e-04 \\ 
};

\node[right, align=left, text=black, font=\footnotesize]
at (axis cs:10.5366,1.8041e-04) {4}; 

\end{axis}
\end{tikzpicture}
\end{subfigure}
\begin{subfigure}[b]{0.49\textwidth}
\begin{tikzpicture}
\begin{axis}[%
width=0.65\textwidth,
height=0.5\textwidth,
at={(0\textwidth,0\textwidth)},
scale only axis,
xmode=log,
xmin=2.5,
xmax=12,
xminorticks=true,
xlabel={$1/h$},
ymode=log,
ymin=0.0001,
ymax=0.2,
yminorticks=true,
ylabel={$dG$-Errors},
legend style={draw=none,fill=none,legend cell align=left},
legend pos=south west
]
\addplot [color=myred,solid,line width=1.5pt, mark=diamond*,mark options={color=myred}]
  table[row sep=crcr]{
    2.9017   0.035003 \\
    3.9303   0.012347 \\
    6.0098   0.0031719 \\
    8.2911   0.0010225 \\
    10.5366   0.0003475 \\
};
\addlegendentry{$\mathbf{u}(\mathbf{x},t)$}

\addplot [color=myblue,solid,line width=1.5pt,mark=square*,mark options={color=myblue}]
  table[row sep=crcr]{
    2.9017   0.10271 \\
    3.9303   0.035519 \\
    6.0098   0.0092151 \\
    8.2911   0.0028222 \\
    10.5366   0.00096199 \\
};
\addlegendentry{$\mathbf{w}(\mathbf{x},t)$}

\addplot [color=mygreen,solid,line width=1.5pt,mark=triangle*,mark options={color=mygreen}]
  table[row sep=crcr]{
    2.9017   0.061029 \\
    3.9303   0.026673 \\
    6.0098   0.0088726 \\
    8.2911   0.0033583 \\
    10.5366   0.0017277 \\
};
\addlegendentry{$T(\mathbf{x},t)$}

\addplot [color=black,solid,line width=0.5pt]
  table[row sep=crcr]{
 8.2911      0.0062 \\
 10.5366     0.0062 \\
 10.5366     0.003 \\
 8.2911      0.0062 \\ 
};
\node[right, align=left, text=black, font=\footnotesize]
at (axis cs:10.5366,0.0046) {3}; 

\end{axis}
\end{tikzpicture}
\end{subfigure}
\end{subfigure}
\caption{Convergence test for $\tau = 0.01$ : computed errors in $L^2$-norm (left) and $dG$-norm (right) versus $1/h$ (\textit{log-log} scale). The errors are computed at the final time $T_f$. The polynomial degree of approximation is taken as $\ell = 3$.}
\label{fig:convH}
\end{figure}
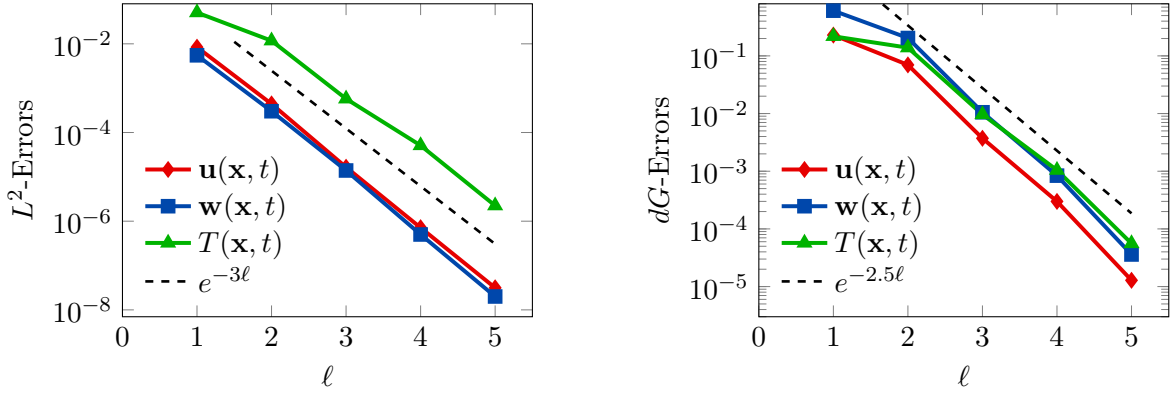
\begin{figure}[ht]
\begin{subfigure}{1\textwidth}
\centering
\begin{subfigure}[b]{0.49\textwidth}
\begin{tikzpicture}
\begin{axis}[%
width=0.65\textwidth,
height=0.5\textwidth,
at={(0\textwidth,0\textwidth)},
scale only axis,
xmin=0,
xmax=5.5,
xminorticks=true,
xlabel={$\ell$},
ymode=log,
ymin=7e-09,
ymax=0.08,
yminorticks=true,
ylabel={$L^2$-Errors},
legend style={draw=none,fill=none,legend cell align=left},
legend pos=south west
]
\addplot [color=myred,solid,line width=1.5pt, mark=diamond*,mark options={color=myred}]
  table[row sep=crcr]{
    1   0.0083972 \\
    2   0.0004411 \\
    3   1.6478e-05 \\
    4   7.134e-07 \\
    5   3.1443e-08 \\
};
\addlegendentry{$\mathbf{u}(\mathbf{x},t)$}

\addplot [color=myblue,solid,line width=1.5pt,mark=square*,mark options={color=myblue}]
  table[row sep=crcr]{
    1   0.0055074 \\
    2   0.00030197 \\
    3   1.3912e-05 \\
    4   5.0098e-07 \\
    5   2.0029e-08 \\
};
\addlegendentry{$\mathbf{w}(\mathbf{x},t)$}

\addplot [color=mygreen,solid,line width=1.5pt,mark=triangle*,mark options={color=mygreen}]
  table[row sep=crcr]{
    1   0.050701 \\
    2   0.011744 \\
    3   0.00057611 \\
    4   5.0986e-05 \\
    5   2.2409e-06 \\
};
\addlegendentry{$T(\mathbf{x},t)$}

\addplot [color=black,dashed,line width=1.0pt]
  table[row sep=crcr]{
    1.5 11.1090e-003
    2   2.4788e-03 \\
    3   123.4098e-06 \\
    4   6.1442e-06 \\
    5   305.9023e-09 \\
};
\addlegendentry{$e^{-3 \ell}$}

\end{axis}
\end{tikzpicture}
\end{subfigure}
\begin{subfigure}[b]{0.49\textwidth}
\begin{tikzpicture}
\begin{axis}[%
width=0.65\textwidth,
height=0.5\textwidth,
at={(0\textwidth,0\textwidth)},
scale only axis,
xmin=0,
xmax=5.5,
xminorticks=true,
xlabel={$\ell$},
ymode=log,
ymin=3e-6,
ymax=0.8,
yminorticks=true,
ylabel={$dG$-Errors},
legend style={draw=none,fill=none,legend cell align=left},
legend pos=south west
]
\addplot [color=myred,solid,line width=1.5pt, mark=diamond*,mark options={color=myred}]
  table[row sep=crcr]{
    1   0.22939 \\
    2   0.069741 \\
    3   0.0037224 \\
    4   0.00029955 \\
    5   1.2845e-05 \\
};
\addlegendentry{$\mathbf{u}(\mathbf{x},t)$}

\addplot [color=myblue,solid,line width=1.5pt,mark=square*,mark options={color=myblue}]
  table[row sep=crcr]{
    1   0.60782 \\
    2   0.20409 \\
    3   0.010525 \\
    4   0.00083873 \\
    5   3.6e-05 \\
};
\addlegendentry{$\mathbf{w}(\mathbf{x},t)$}

\addplot [color=mygreen,solid,line width=1.5pt,mark=triangle*,mark options={color=mygreen}]
  table[row sep=crcr]{
    1   0.21848 \\
    2   0.13877 \\
    3   0.0096452 \\
    4   0.0010444 \\
    5   5.5884e-05 \\
};
\addlegendentry{$T(\mathbf{x},t)$}

\addplot [color=black,dashed,line width=1.0pt]
  table[row sep=crcr]{
     1.5    1.1759 \\
     2      336.8973e-3 \\
     3      27.6542e-3 \\
     4      2.2700e-3 \\
     5      186.3327e-6 \\
};
\addlegendentry{$e^{-2.5 \ell}$}

\end{axis}
\end{tikzpicture}
\end{subfigure}
\end{subfigure}
\caption{Convergence test for $\tau = 0.01$ : computed errors in $L^2$-norm (left) and $dG$-norm (right) versus $\ell$ (\textit{semi-log} scale). The errors are computed at the final time $T_f$. The computational mesh is made of $100$ polygons.}
\label{fig:convEll}
\end{figure}

\begin{figure}[ht]
\begin{subfigure}{1\textwidth}
\centering
\begin{subfigure}[b]{0.49\textwidth}
\begin{tikzpicture}
\begin{axis}[%
width=0.65\textwidth,
height=0.5\textwidth,
at={(0\textwidth,0\textwidth)},
scale only axis,
xmode=log,
xmin=2.5,
xmax=12,
xminorticks=true,
xlabel={$1/h$},
ymode=log,
ymin=3e-4,
ymax=0.05,
yminorticks=true,
ylabel={$L^2$-Errors},
legend style={draw=none,fill=none,legend cell align=left},
legend pos=south west
]

\addplot [color=myyellow,solid,line width=1.5pt,mark=triangle*,mark options={color=myyellow}]
  table[row sep=crcr]{
    2.9017   0.043611 \\
    3.9303   0.016623 \\
    6.0098   0.0047426 \\
    8.2911   0.0015573 \\
    10.5366   0.00051298 \\
};
\addlegendentry{$\tau = 0$}

\addplot [color=mypurple,solid,line width=1.5pt,mark=square*,mark options={color=mypurple}]
  table[row sep=crcr]{
    2.9017   0.038517 \\
    3.9303   0.017446 \\
    6.0098   0.0047434 \\
    8.2911   0.0016325 \\
    10.5366   0.00054019 \\
};
\addlegendentry{$\tau \neq 0$}

\addplot [color=black,solid,line width=0.5pt]
  table[row sep=crcr]{
 8.2911       0.0021 \\
 10.5366      0.0021 \\
 10.5366      0.0008 \\
 8.2911       0.0021 \\ 
};

\node[right, align=left, text=black, font=\footnotesize]
at (axis cs:10.5366,0.0014) {4}; 

\end{axis}
\end{tikzpicture}
\end{subfigure}
\begin{subfigure}[b]{0.49\textwidth}
\begin{tikzpicture}
\begin{axis}[%
width=0.65\textwidth,
height=0.5\textwidth,
at={(0\textwidth,0\textwidth)},
scale only axis,
xmode=log,
xmin=2.5,
xmax=12,
xminorticks=true,
xlabel={$1/h$},
ymode=log,
ymin=0.2,
ymax=5.5,
yminorticks=true,
ylabel={$dG$-Errors},
legend style={draw=none,fill=none,legend cell align=left},
legend pos=south west
]

\addplot [color=myyellow,solid,line width=1.5pt,mark=triangle*,mark options={color=myyellow}]
  table[row sep=crcr]{
    2.9017   4.2873 \\
    3.9303   2.6689 \\
    6.0098   1.204 \\
    8.2911   0.56923 \\
    10.5366   0.27159 \\
};
\addlegendentry{$\tau = 0$}

\addplot [color=mypurple,solid,line width=1.5pt,mark=square*,mark options={color=mypurple}]
  table[row sep=crcr]{
    2.9017   4.8006 \\
    3.9303   2.4533 \\
    6.0098   1.177 \\
    8.2911   0.53374 \\
    10.5366   0.24903 \\
};
\addlegendentry{$\tau \neq 0$}

\addplot [color=black,solid,line width=0.5pt]
  table[row sep=crcr]{
 8.2911      0.7183 \\
 10.5366     0.7183 \\
 10.5366     0.35 \\
 8.2911      0.7183 \\ 
};
\node[right, align=left, text=black, font=\footnotesize]
at (axis cs:10.5366,0.5342) {3}; 

\end{axis}
\end{tikzpicture}
\end{subfigure}
\end{subfigure}
\caption{Convergence test: computed errors for the reconstructed pressure vs $h$, considering both the cases $\tau = 0$ and $\tau = 0.01$. The polynomial degree of approximation is takes as $\ell = 3$}
\label{fig:convPressure}
\end{figure}
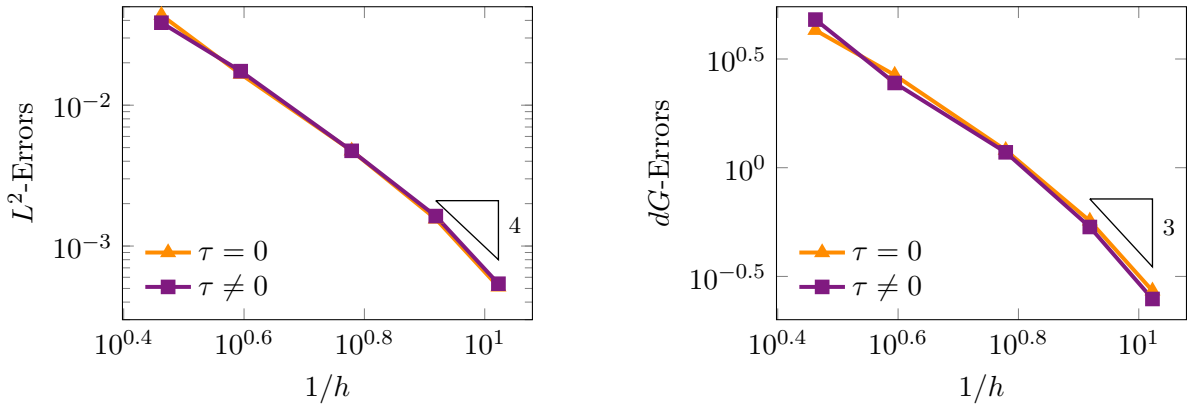

In Figure~\ref{fig:convH_tau0} we report the $L^2$ and $dG$-errors for the three variables with respect to the mesh size (log-log scale). 
In agreement with Theorem~\ref{thm:error_est}, we observe that, as we are using $\ell = 3$, the errors show a convergence rate proportional to $h^3$. Moreover, for what concerns the $L^2$-errors, we observe that we reach $h^{\ell + 1}$ convergence. Notice that this behavior is not covered by our theoretical analysis. In Figure~\ref{fig:convH} we report the same quantities, but for the case $\tau > 0$ ($\tau = 0.01$ according to Table~\ref{tab:params_convtest}). We observe that, even by considering the fully-hyperbolic problem, we recover the theoretical results, cf. Theorem~\ref{thm:error_est}.
We have also computed the $L^2$ and $dG$-errors for the pressure field $p_h$, which is recovered in the post-processing procedure by the use of \eqref{eq:const_law_press}. 
The initial condition for the pressure is taken such that $p_0(\mathbf{x}) = -c_0^{-1}(\alpha \div{\mathbf{u}_0} + \div{\mathbf{w}_0})$, being $\mathbf{u}_0$ and $\mathbf{w}_0$ the initial conditions for the displacement and the filtration displacement, respectively. 
For computing the $dG$-error, we have considered the following $dG$-norm for the pressure \cite{AntoniettiBonetti2022}: $\|p_h\|_{dG,prs}^2 = \|\sqrt{k} \nabla_h p_h\|^2 + \sum_{F\in\mathcal{F}_h} \| \sqrt{\gamma} \jump{p_h} \ \|_{F}^2$, where $\gamma$ is defined in a way similar to \eqref{eq:stabilization_func}. As for the three unknowns of the problem we observe a decay of the error proportional to $h^3$ ($h^4$ in $L^2$-norm), cf. Figure~\ref{fig:convPressure}. We remark that, due to the formulation we are using, we do not have an explicit bound on the $dG$-norm for the pressure (cf. \eqref{eq:const_law_press} and hypotheses of Theorem~\ref{thm:error_est}); however, we observe an optimal rate of decay also in this case. We think that the effect of not having control over the $dG$-norm is seen in the magnitude of the errors, which is higher compared with the orders of magnitude of the errors shown in Figure~\ref{fig:convH_tau0} and Figure~\ref{fig:convH}. Last, the $\ell$-convergence test shows that, in agreement with the theoretical estimates, we reach an exponential decrease of the error, the results of this test are presented in Figure~\ref{fig:convEll}.

\section{Physical tests}
\label{sec:numerical_results}
The aim of this section is to evaluate the proposed scheme with respect to a wide set of physically-sound test cases. We compare the results obtained by our method with the ones presented in \cite{Carcione2019} simulating the wave propagation in a homogeneous and heterogeneous media. A comparison with the results obtained via the poroelastic model is presented too. The set up for numerical implementation is the same as the one described in the preamble of Section~\ref{sec:convergence_test}.

\subsection{Test case 1: homogeneous media}
In this section we consider a wave propagation problem in a homogeneous thermo-poroelastic medium inspired by \cite{Carcione2019}. The aim of this simulation is to prove that our scheme can reproduce known results present in the literature. We consider a domain $\Omega = (0, 1500)^2$ \si{\meter \squared} and in Table~\ref{tab:TPEcoeff_test1} we report the thermo-poroelastic properties of the medium.
\begin{table}[ht]
    \centering 
    \begin{tabular}{ l | l  c c l | l}
    \textbf{Coefficient} & \textbf{Value} & & & \textbf{Coefficient} & \textbf{Value} \T\B \\
    $a_0$ [\si[per-mode = symbol]{\pascal \per \kelvin \squared}] & 4.1695 & & & $k$ [\si[per-mode = symbol]
    {\metre\squared \per \pascal \per \second}] & \num[exponent-product=\ensuremath{\cdot}]{1e-9} \\
    $b_0$ [\si{\per \kelvin}] & \num[exponent-product=\ensuremath{\cdot}]{1.4361e-5} & & & $\theta$ [\si[per-mode = symbol]{\metre\squared \pascal \per \kelvin\squared \per \second}] & \num[exponent-product=\ensuremath{\cdot}]{1.5e+4}  \\
    $c_0$ [\si{\per\pascal}] & \num[exponent-product=\ensuremath{\cdot}]{1.4361e-10} & & & $\rho_f$ [\si[per-mode = symbol]
    {\kilogram \per \meter\cubed}] & 1000 \\
    $\alpha$ [-] & 0.9514 & & & $\rho_s$ [\si[per-mode = symbol]
    {\kilogram \per \meter\cubed}] & 2650 \\
    $\beta$ [\si[per-mode = symbol]{\pascal \per \kelvin}] & \num[exponent-product=\ensuremath{\cdot}]{2.4857e+4} & & & $\phi$ [-] & 0.3 \\ 
    $\mu$ [\si{\pascal}] & \num[exponent-product=\ensuremath{\cdot}]{1.885e+9} & & & $a$ [-] & 2 \\
    $\lambda$ [\si{\pascal}] & \num[exponent-product=\ensuremath{\cdot}]{4.433e+8} & & & $\tau$ [\si{\second}] & \num[exponent-product=\ensuremath{\cdot}]{1.5e-2} \\
    \end{tabular}
    \\[10pt]
    \caption{Test case 1: TPE medium properties}
    \label{tab:TPEcoeff_test1}
\end{table}
We model the shear source in terms of a moment tensor $\mathbf{M}$ as $\mathbf{f} = - \mathbf{M} \, \div{\delta(\mathbf{x} - \mathbf{x}_s)} \, h(t)$ \cite{Morency2008}, where $\mathbf{x}_s$ is the point-source location, $\delta(\mathbf{x} - \mathbf{x}_s)$ is the Kronecker delta located in $\mathbf{x}_s$, and $h(t)$ is the time-history. This form of $\mathbf{f}$ is often used in the context of earthquakes. In our case, the time evolution is given by \cite{Carcione2019} $h(t) = A_0 \cos\left[ 2 \pi (t - t_0) f_0 \right] \exp\left[ -2 (t - t_0)^2 f_0^2 \right]$, where $A_0 = 10$ \si{\meter} is the amplitude, $f_0 = 5$ \si{\hertz} is the peak-frequency, and $t_0 = 3/(2 f_0) = 0.3$ \si{\second} is the time-shift. To discretize our domain we choose a polygonal mesh with mesh size $h \sim 60$ \si{\meter} (\# Elements = 2500) and polynomial degree $\ell = 4$. As a time stepping scheme we employ the Newmark-$\beta$ scheme, with $\Delta t = \num[exponent-product=\ensuremath{\cdot}]{1e-2}$ and $T_f = 1$ \si{\second}. Finally, we complete our problem with homogeneous Dirichlet boundary conditions and with null initial conditions.
In the following, we denote by $\mathbf{v}_h$ the solid velocity (i.e. $\dot{ \mathbf{u}}_h$), by $v_{h,y}$ its vertical component, and by $q_{h,y}$ the vertical component of the filtration velocity (i.e. $\left( \dot{ \mathbf{w}}_h \right)_y$ ).

\begin{figure}[ht]
\begin{subfigure}[b]{.33\textwidth}
    \centering
    \includegraphics[width=1\textwidth]{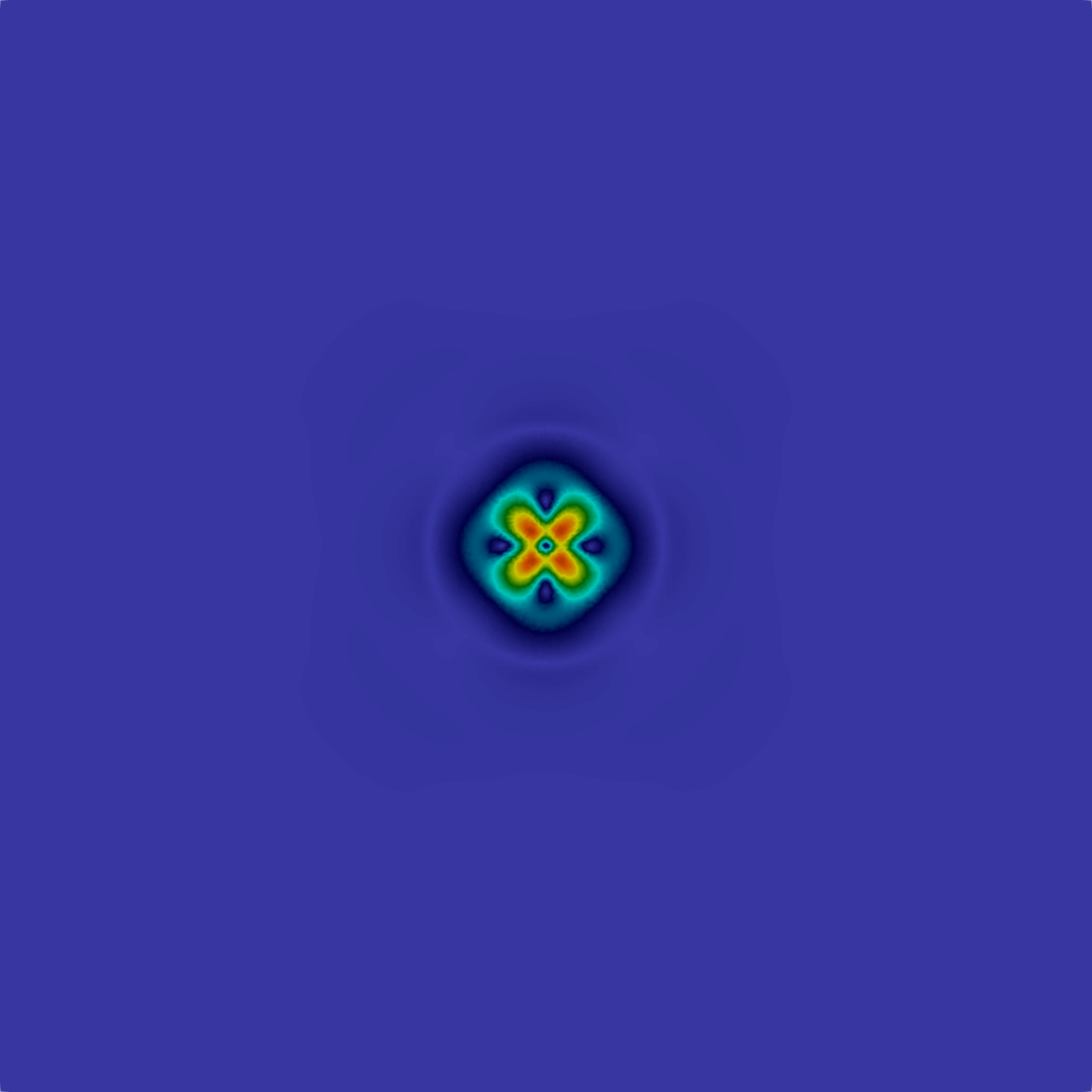}
    \label{fig:v_test1_02}
\end{subfigure}
\begin{subfigure}[b]{.33\textwidth}
    \centering
    \includegraphics[width=1\textwidth]{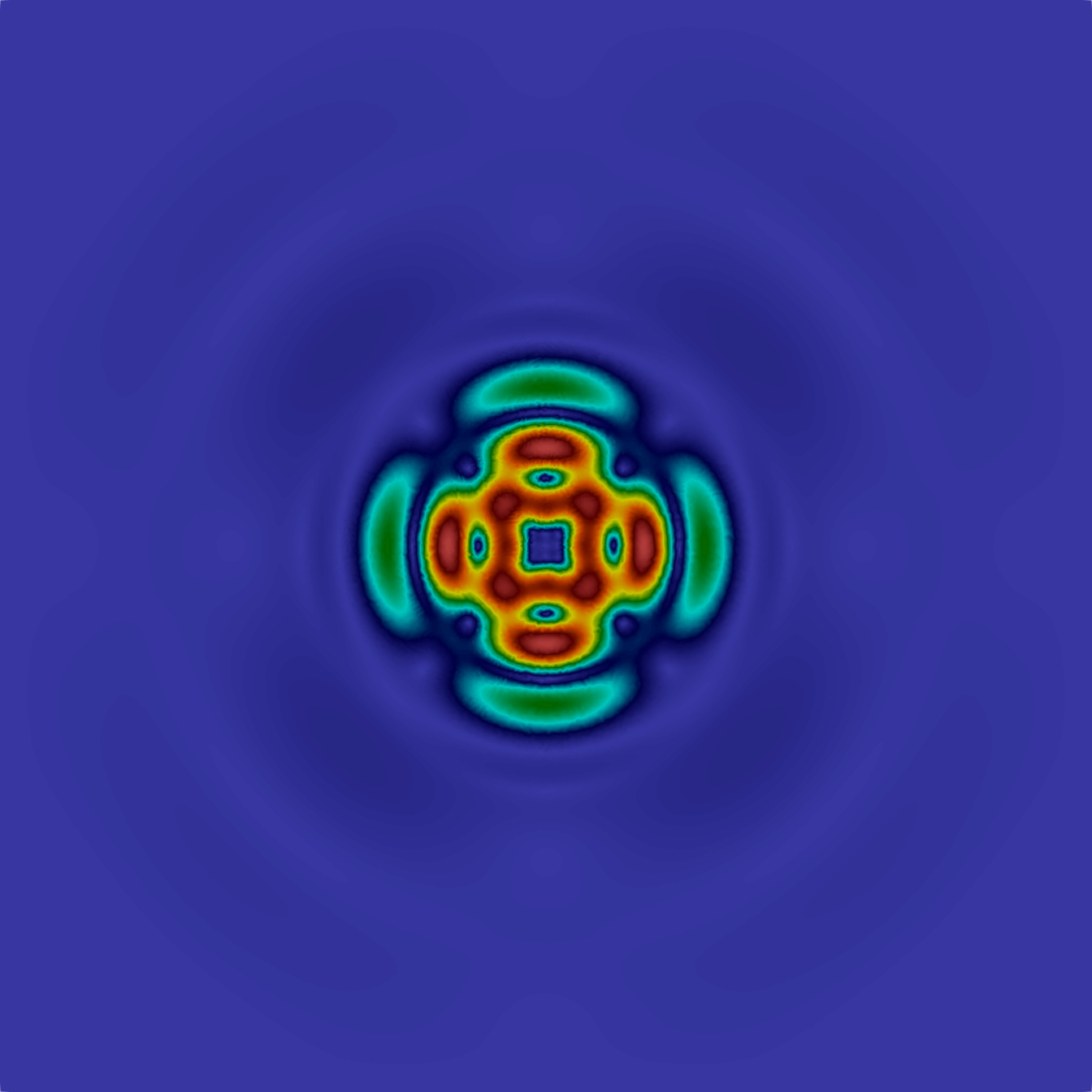}
    \label{fig:v_test1_04}
\end{subfigure}
\begin{subfigure}[b]{.33\textwidth}
    \centering
    \includegraphics[width=1\textwidth]{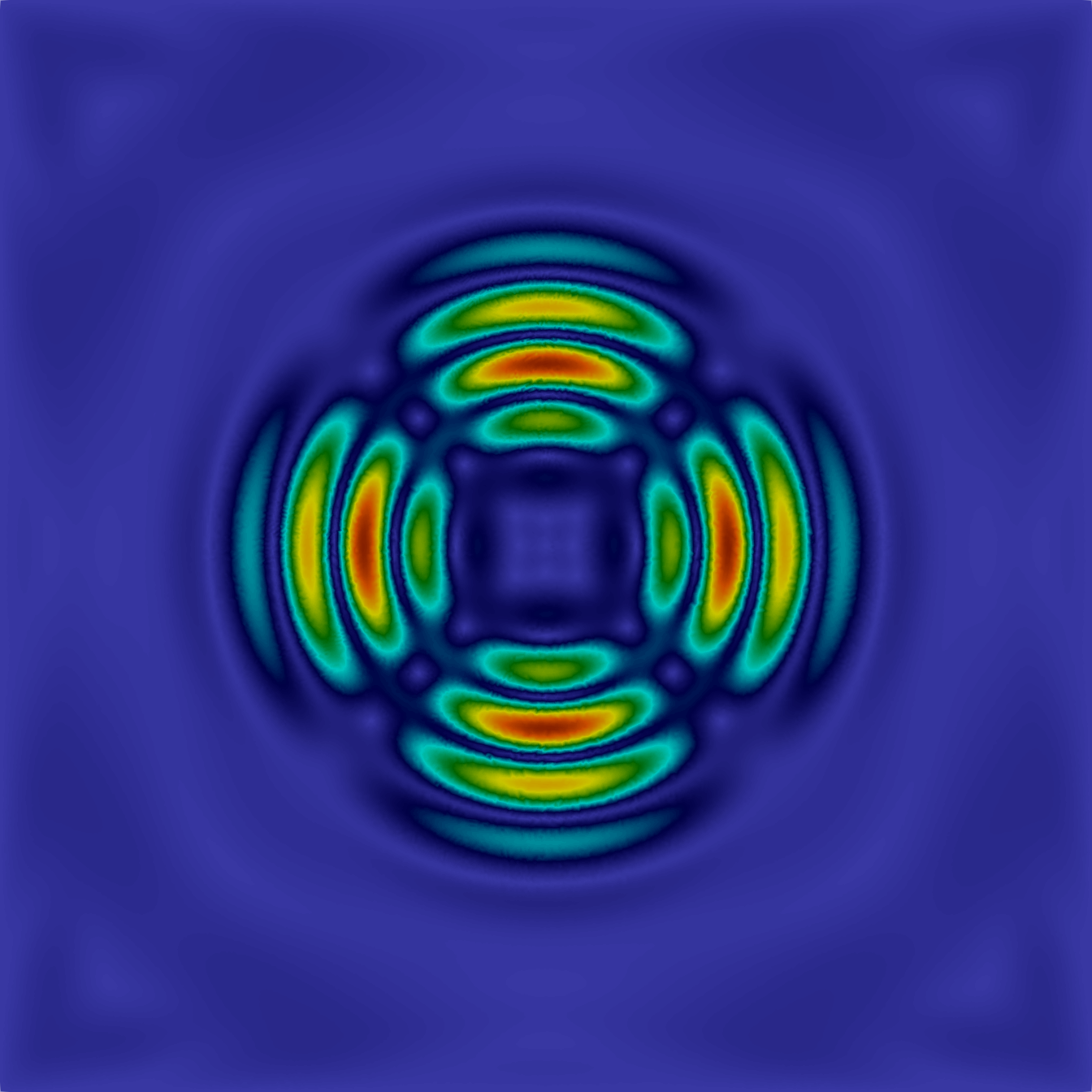}
    \label{fig:v_test1_06}
\end{subfigure}

\vspace{-0.4cm}
\centering
\includegraphics[width=0.4\textwidth]{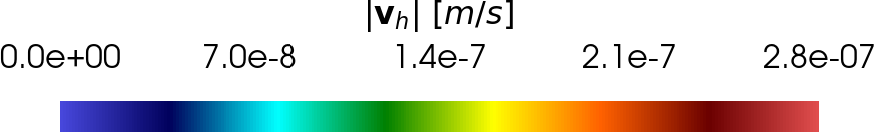}

\caption{Test case 1: computed velocity field $|\mathbf{v}_h|$ at the time instants $t=0.2 \si{\second}$ (left), $t=0.4 \si{\second}$ (center), $t=0.6 \si{\second}$ (right).}
\label{fig:v_test1}
\end{figure}

\begin{figure}[ht]
\begin{subfigure}[b]{.33\textwidth}
    \centering
    \includegraphics[width=1\textwidth]{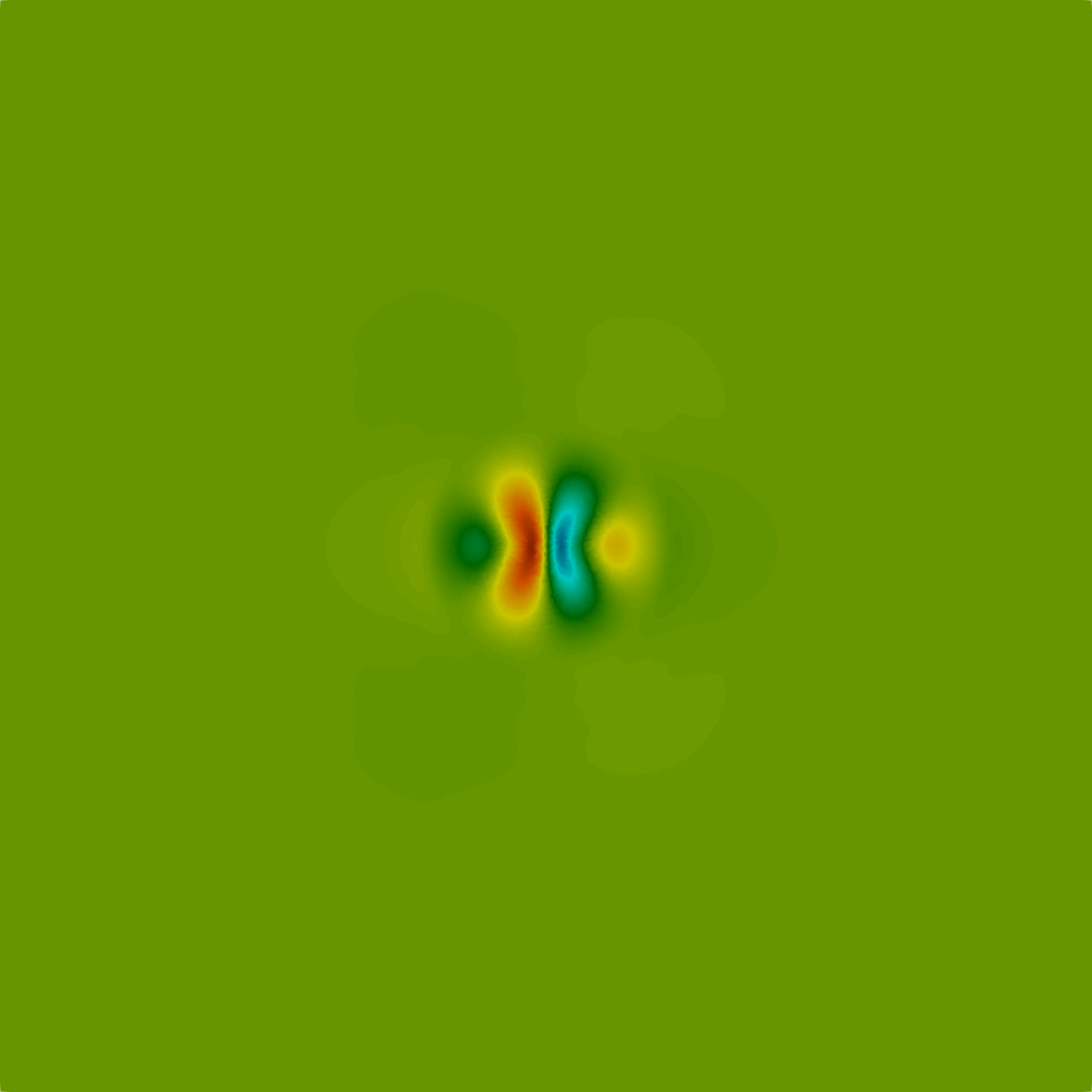}
    \label{fig:v2_test1_02}
\end{subfigure}
\begin{subfigure}[b]{.33\textwidth}
    \centering
    \includegraphics[width=1\textwidth]{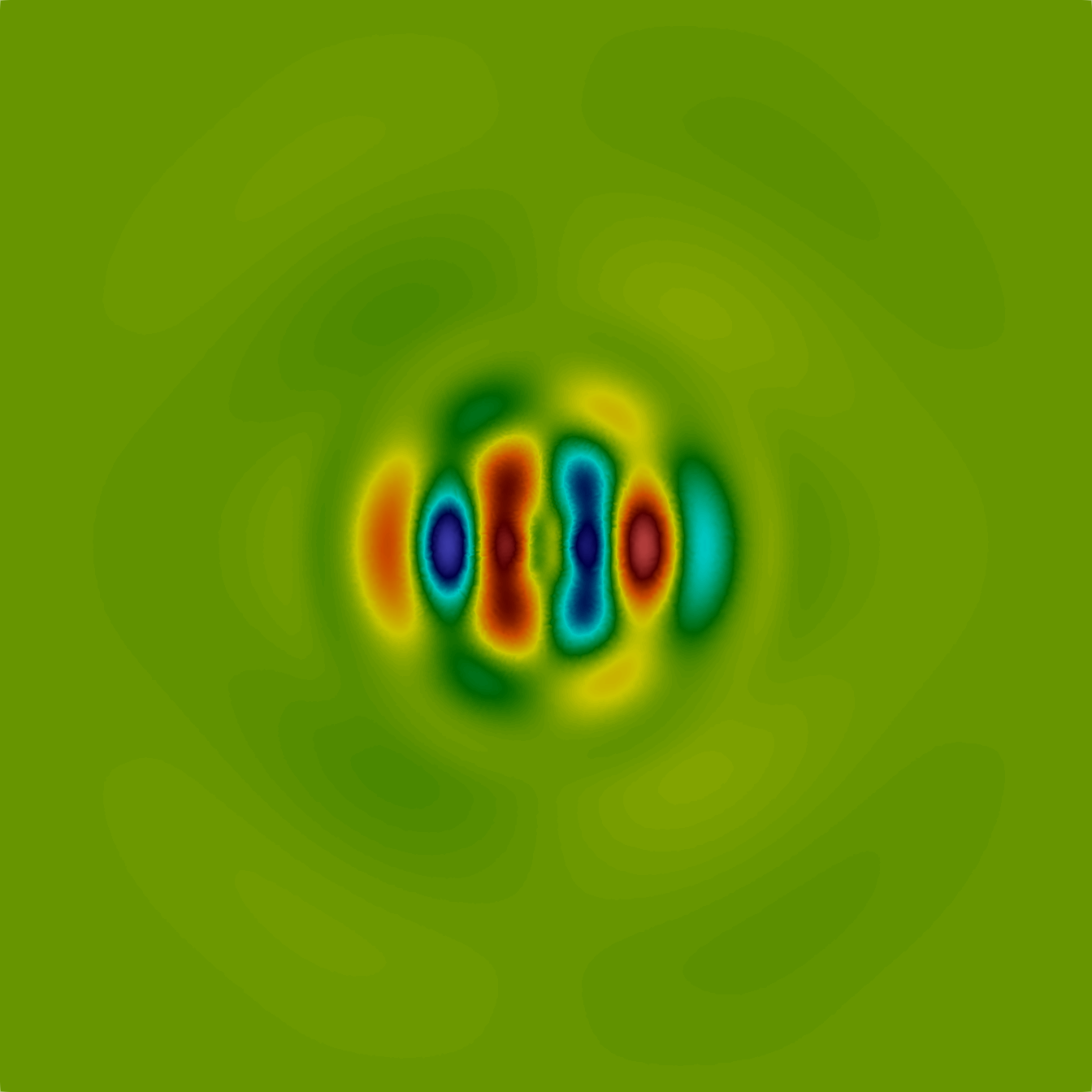}
    \label{fig:v2_test1_04}
\end{subfigure}
\begin{subfigure}[b]{.33\textwidth}
    \centering
    \includegraphics[width=1\textwidth]{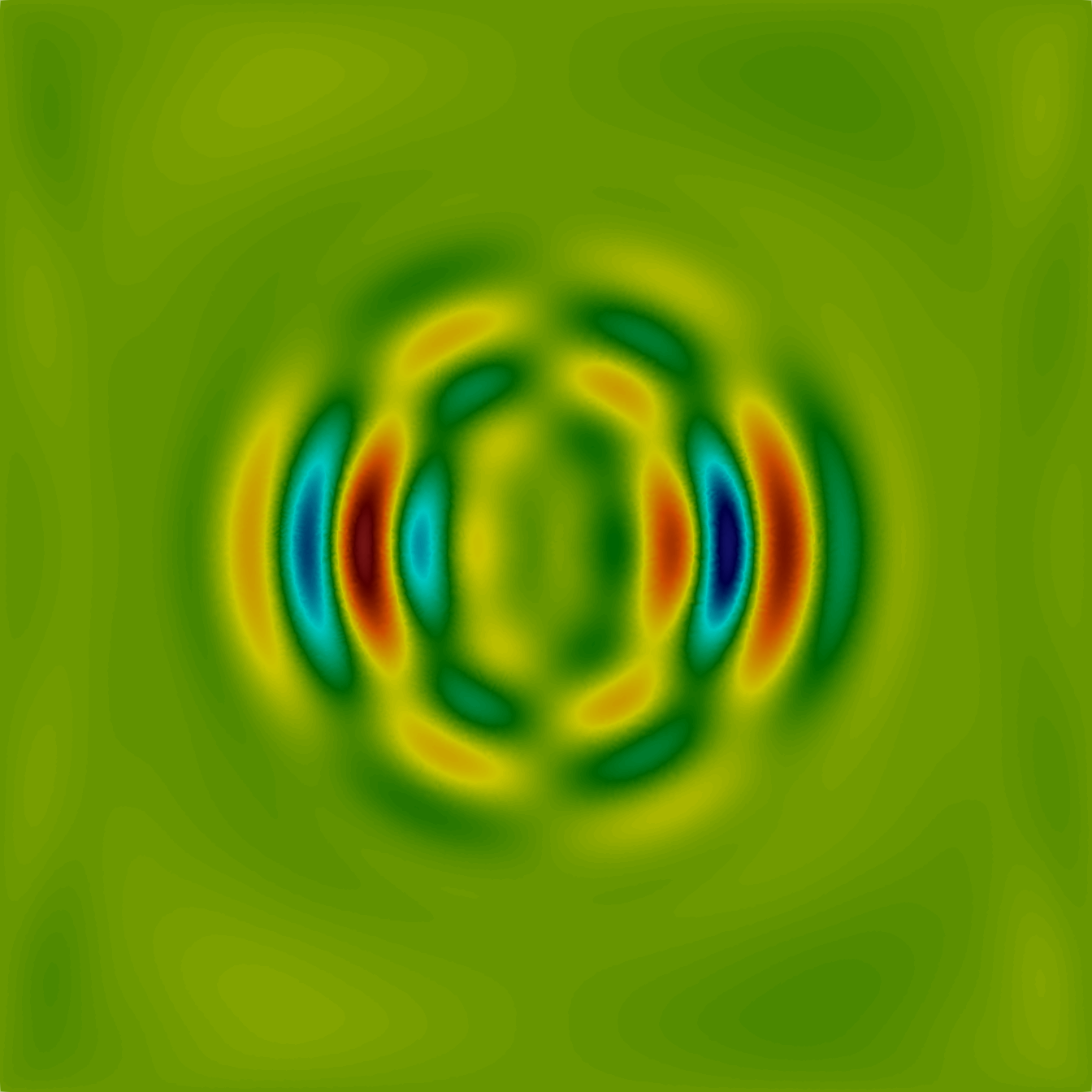}
    \label{fig:v2_test1_06}
\end{subfigure}

\vspace{-0.4cm}
\centering
\includegraphics[width=0.4\textwidth]{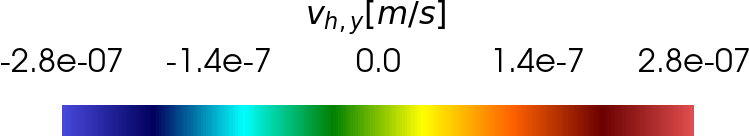}

\caption{Test case 1: computed vertical component of the velocity field $v_{h,y}$ at the time instants $t=0.2s$ (left), $t=0.4s$ (center), $t=0.6s$ (right)}
\label{fig:v2_test1}
\end{figure}

\begin{figure}[ht]
\begin{subfigure}[b]{.33\textwidth}
    \centering
    \includegraphics[width=1\textwidth]{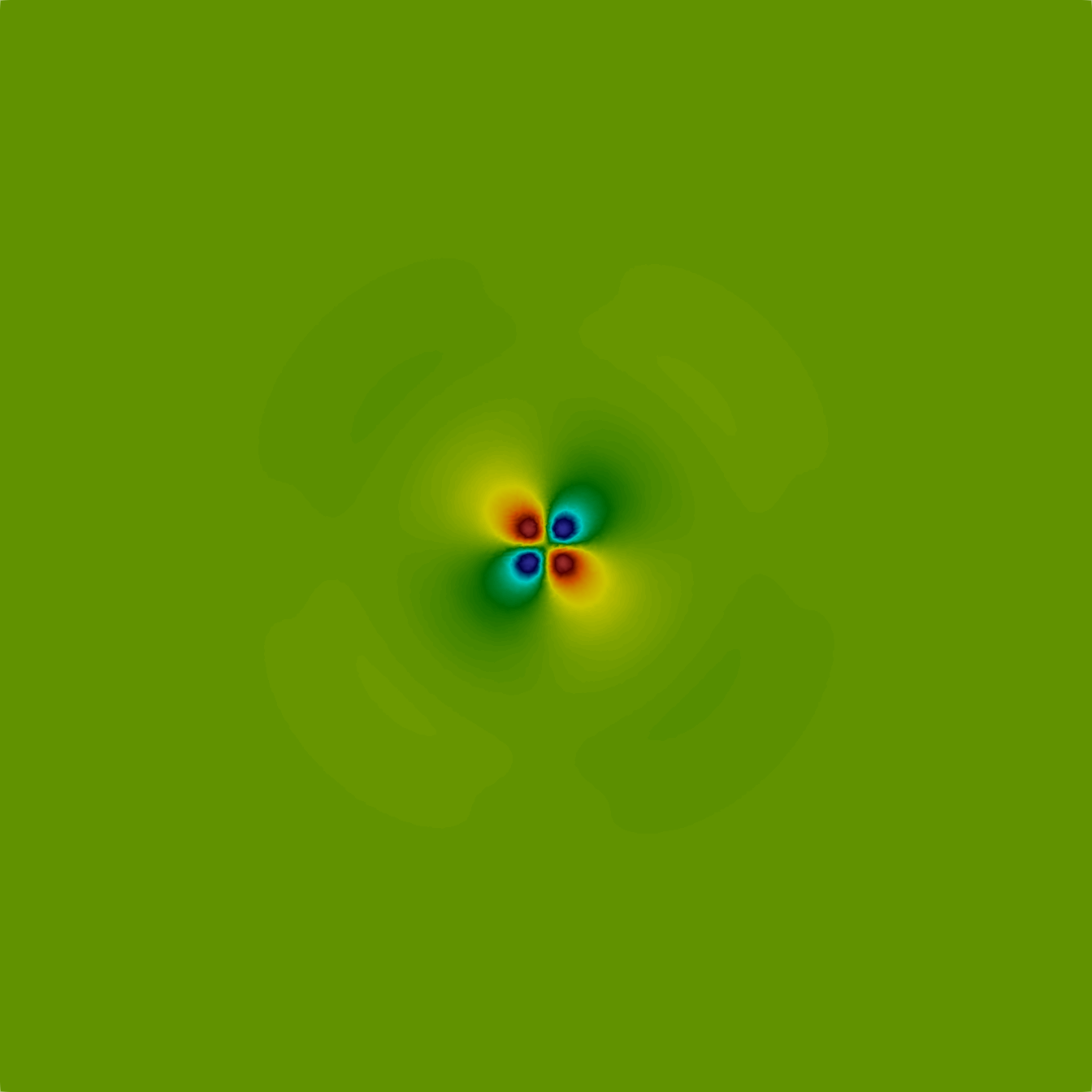}
    \label{fig:T_test1_02}
\end{subfigure}
\begin{subfigure}[b]{.33\textwidth}
    \centering
    \includegraphics[width=1\textwidth]{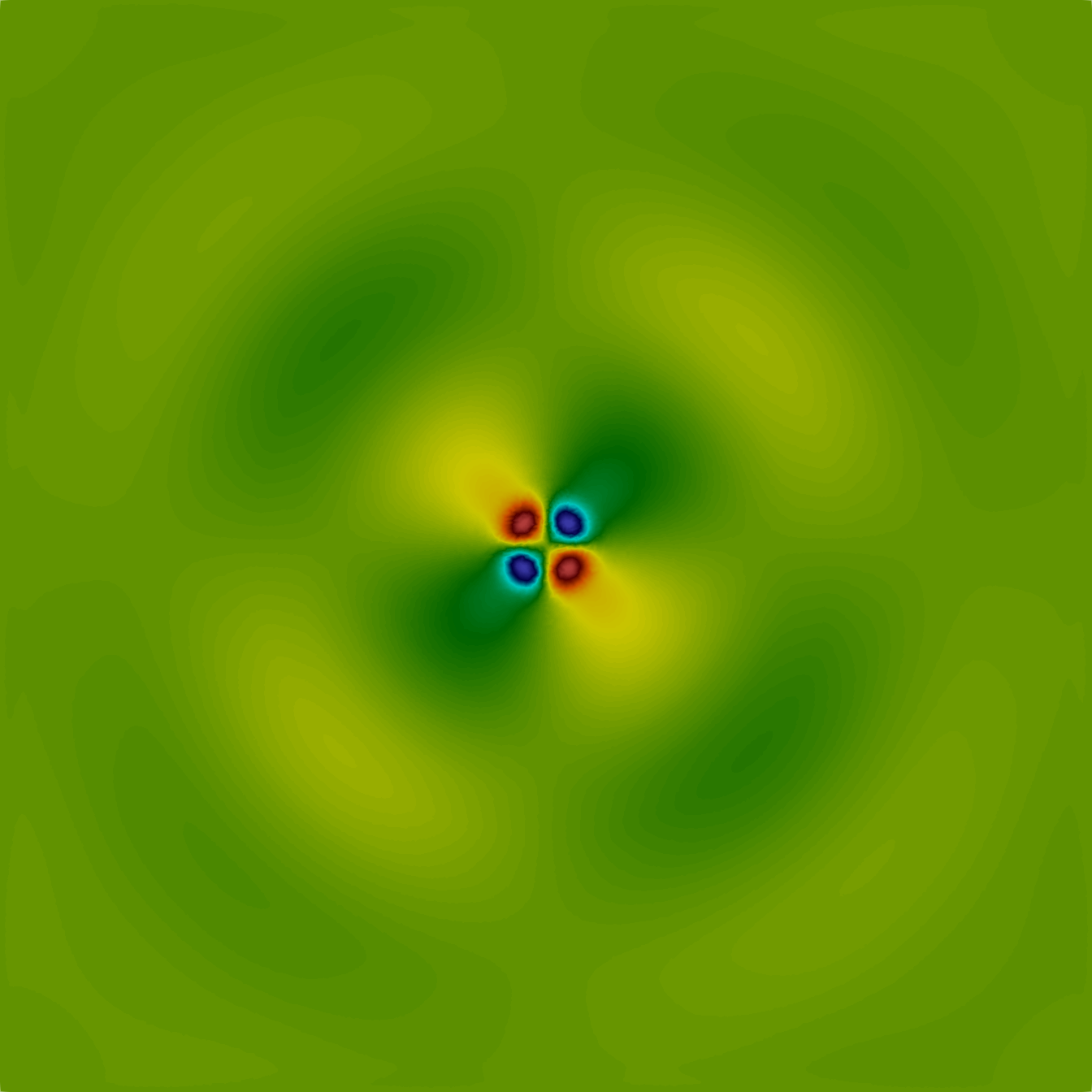}
    \label{fig:T_test1_04}
\end{subfigure}
\begin{subfigure}[b]{.33\textwidth}
    \centering
    \includegraphics[width=1\textwidth]{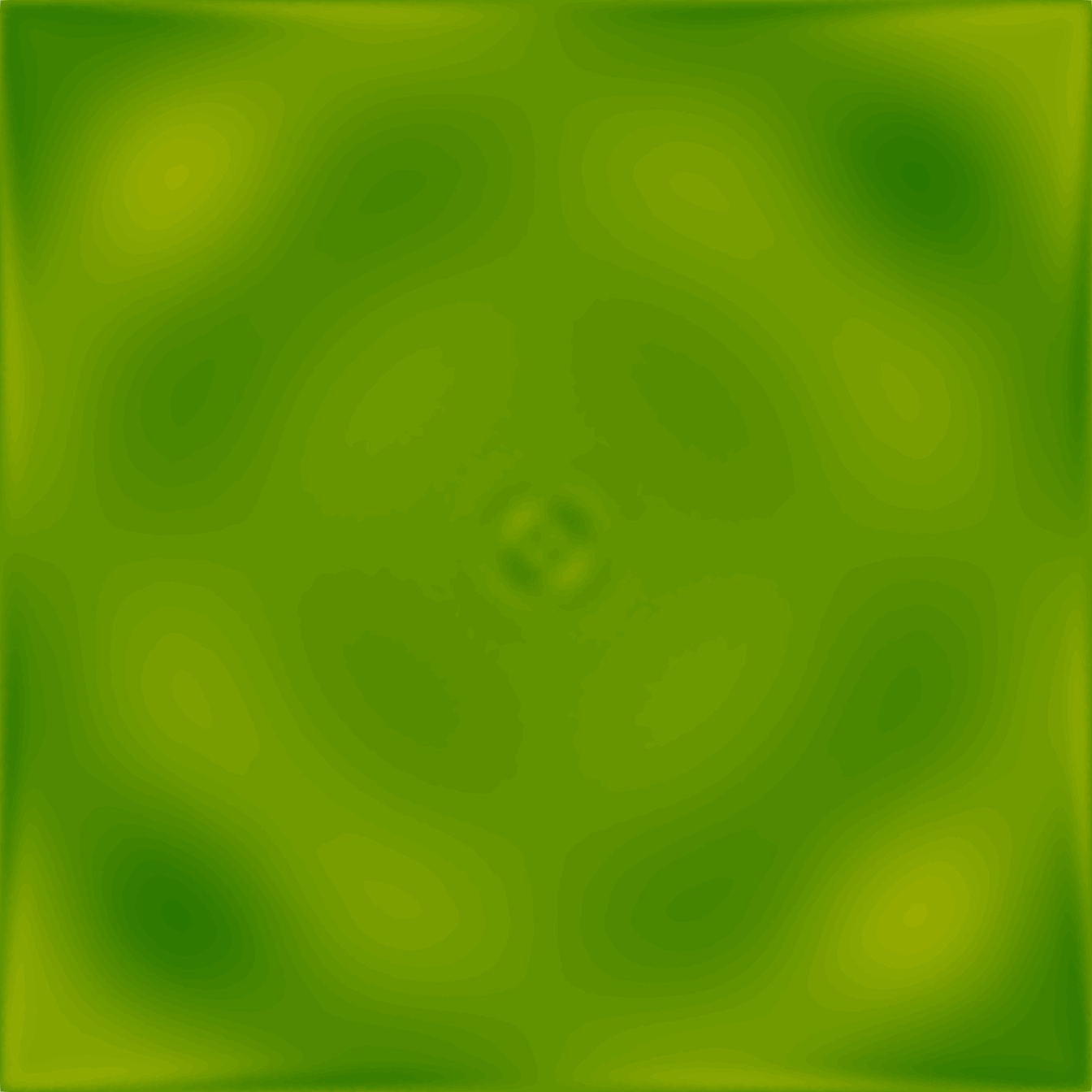}
    \label{fig:T_test1_06}
\end{subfigure}

\vspace{-0.4cm}
\centering
\includegraphics[width=0.4\textwidth]{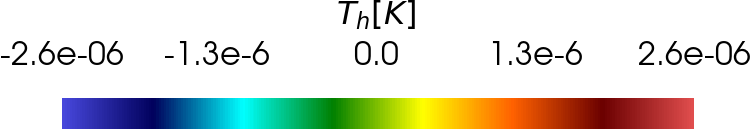}

\caption{Test case 1: computed temperature field $T_h$ at the time instants $t=0.2s$ (left), $t=0.4s$ (center), $t=0.6s$ (right)}
\label{fig:T_test1}
\end{figure}

We report in Figure~\ref{fig:v_test1}, Figure~\ref{fig:v2_test1}, and Figure~\ref{fig:T_test1}, the computed quantities $|\mathbf{v}_h|$, $v_{h,y}$, and $q_{h,y}$ at selected time instants, respectively.

From the results of Figure~\ref{fig:v_test1} we notice a 
symmetric wavefront that detaches from the center of the domain at $t = 0.4$ \si{\second}; this is due to the homogeneity of the thermo-poroelastic material in which it propagates. We can see that the axes of symmetry of our wavefront are the diagonals of the square domain. This behavior is correct and is due to the form of the forcing term we are imposing. 
Looking at the three snapshots we can also observe the presence of the fast $P$-wave captured by our scheme, even if its amplitude is far lower than the ones of the slow $P$-wave and the $S$-wave. We can also notice that in the last frame, the fast $P$-wave reaches the border of the domain, causing some reflection effects in the corners. 

From the results in Figure~\ref{fig:v2_test1} and Figure~\ref{fig:T_test1} we can qualitatively compare our results with the ones presented in \cite{Carcione2019} with very similar parameters. We highlight that the main difference between our test and the one proposed in \cite{Carcione2019} lies in the lower frequency content of the source term, i.e. lower $f_0$. This generates a wavefield with a larger wavelength and makes the identification of the $P$ and $S$ waves more difficult.  In Figure~\ref{fig:v2_test1}, we observe more plainly the appearance of the shear waves and the anti-symmetric pattern of the wave fronts with respect to the $y$-axis. In Figure~\ref{fig:T_test1} we can see the presence of the diffusive thermal $T$-wave. In conclusion, taking into account the value of the central peak frequency $f_0$ that we considered, we can see a good agreement between our results and the ones presented in \cite{Carcione2019}.

\subsection{Test case 2: comparison with the poroelastic model}
\begin{figure}[ht]
\begin{subfigure}[b]{.33\textwidth}
    \centering
    \includegraphics[width=1\textwidth]{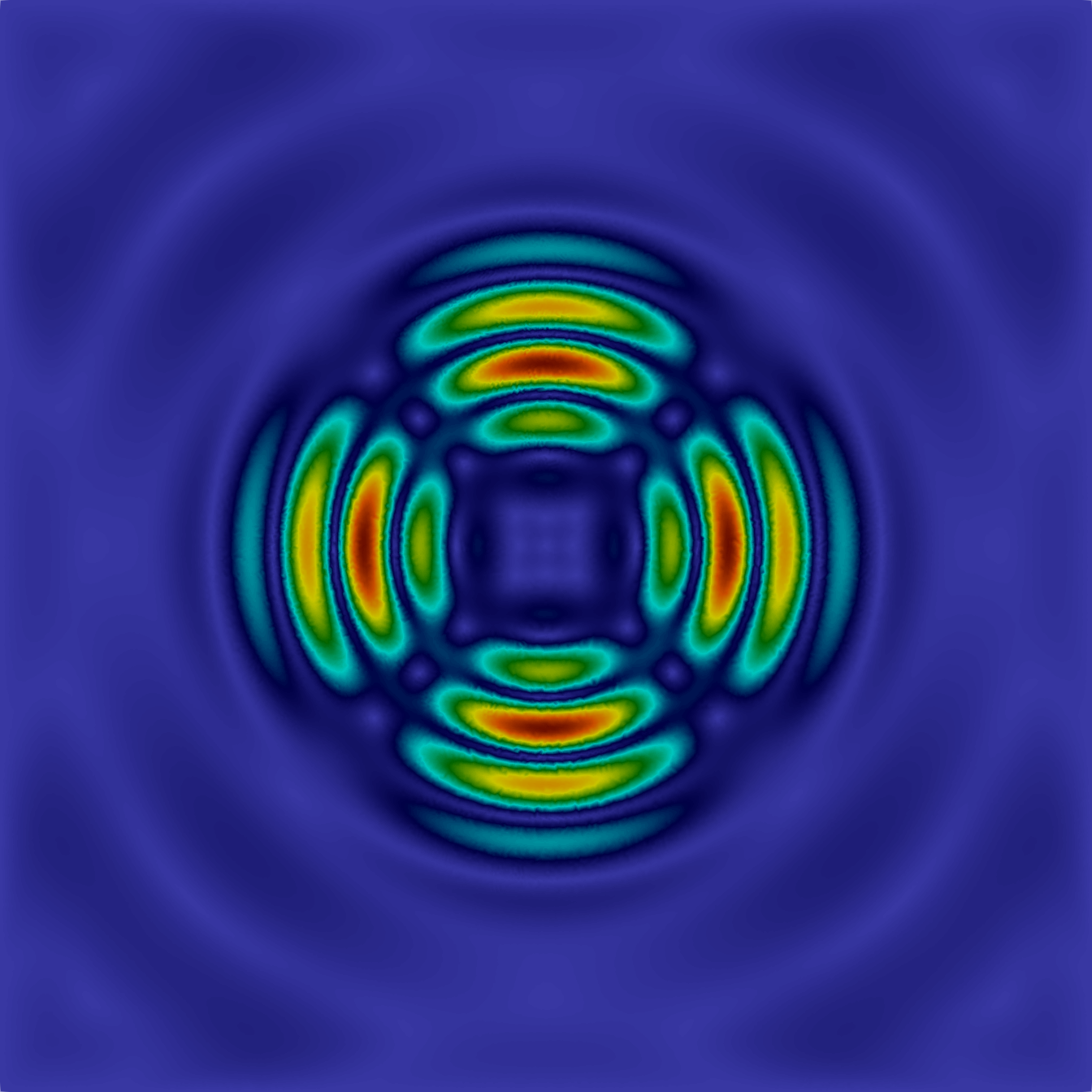}
    \label{fig:v_PE_test1_06}
\end{subfigure}
\begin{subfigure}[b]{.33\textwidth}
    \centering
    \includegraphics[width=1\textwidth]{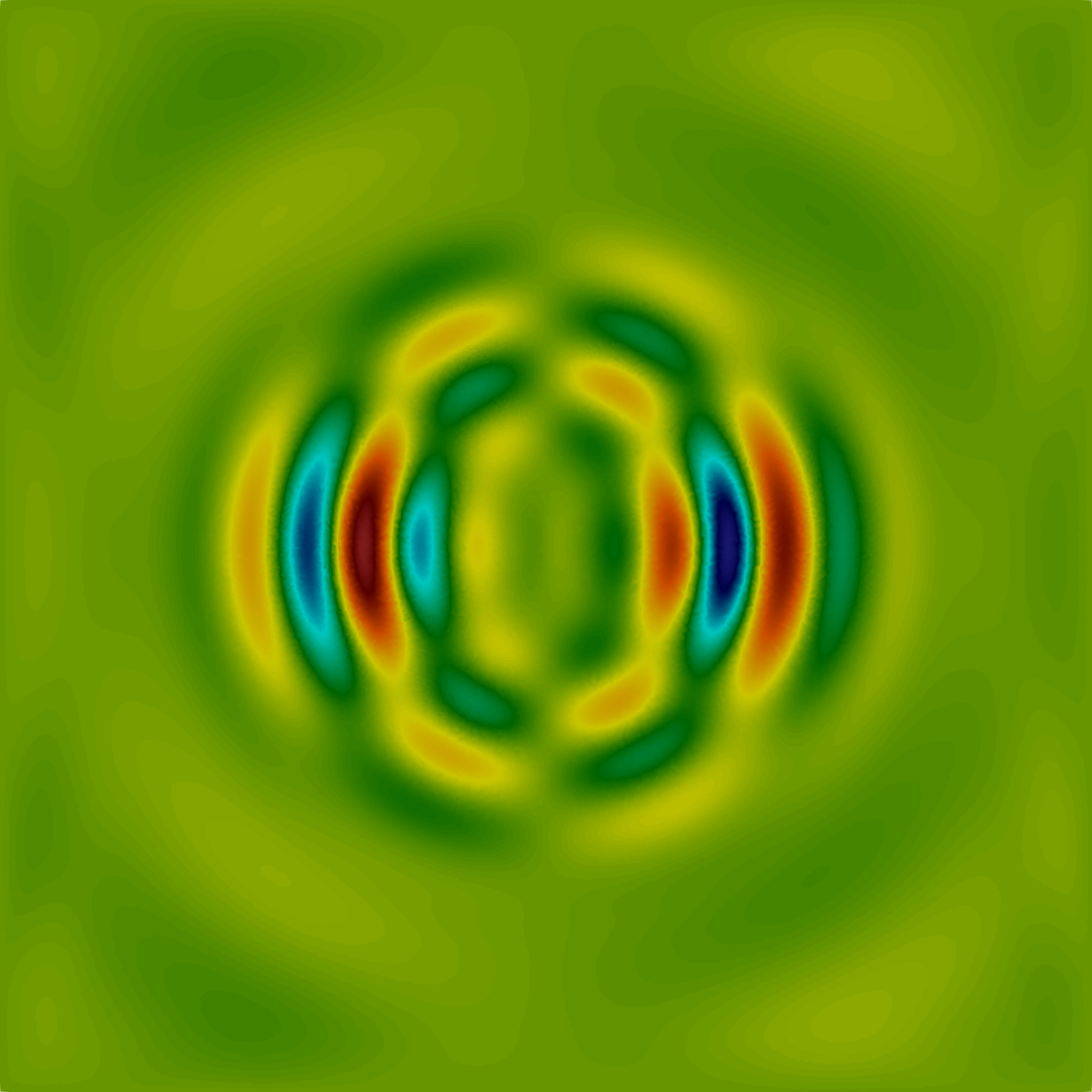}
    \label{fig:v2_PE_test1_06}
\end{subfigure}
\begin{subfigure}[b]{.33\textwidth}
    \centering
    \includegraphics[width=1\textwidth]{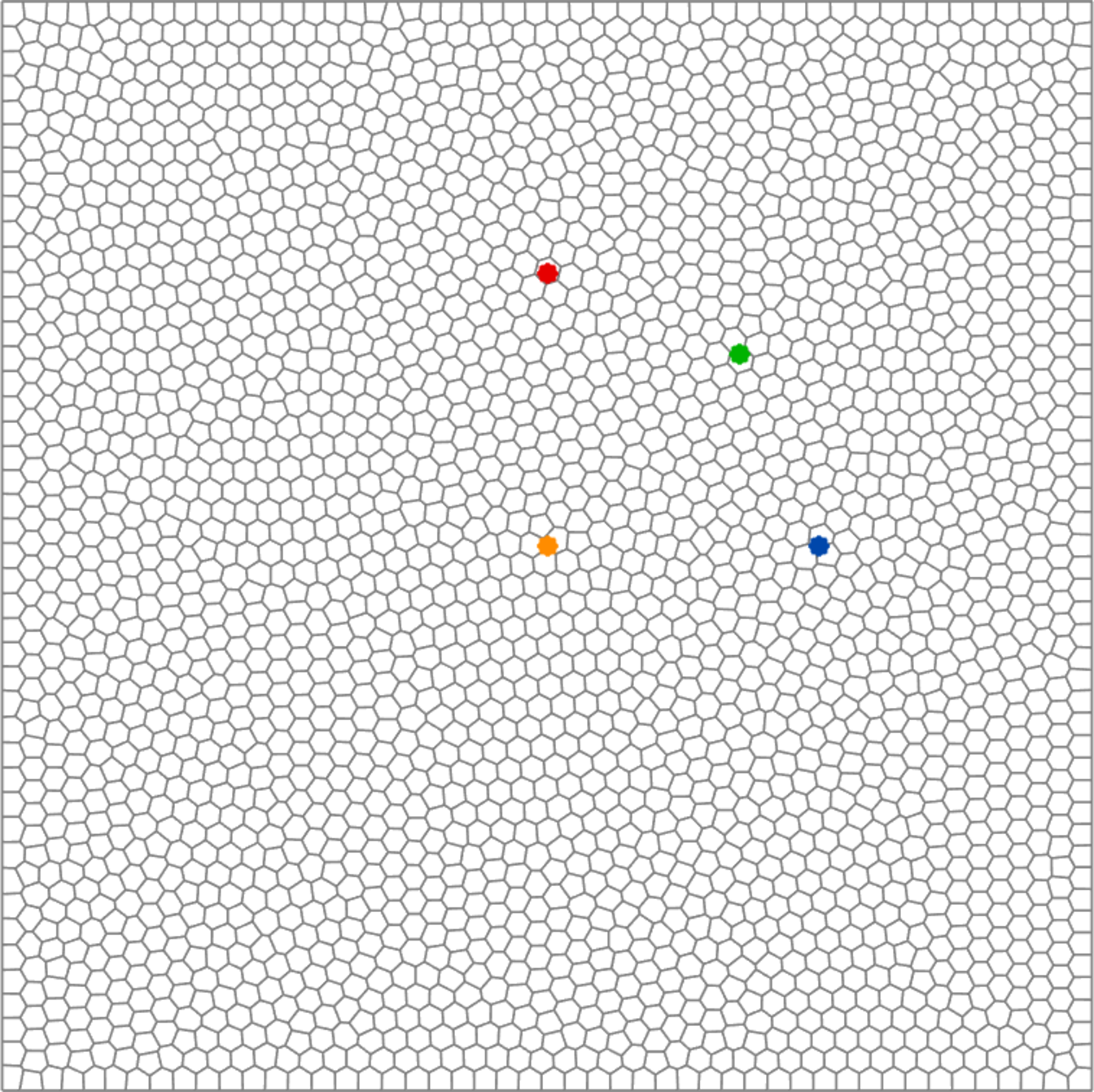}
    \label{fig:mesh_comparison}
\end{subfigure}

\begin{subfigure}[b]{.33\textwidth}
    \centering
    \includegraphics[width=0.8\textwidth]{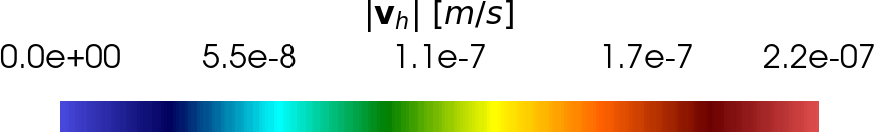}
    \label{fig:v_PE_test1_colormap}
\end{subfigure}
\begin{subfigure}[b]{.33\textwidth}
    \centering
    \includegraphics[width=0.8\textwidth]{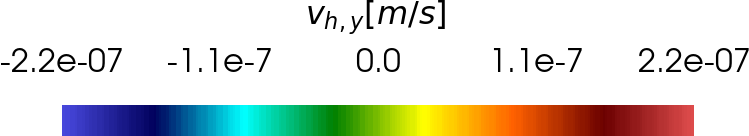}
    \label{fig:v2_PE_test1_colormap}
\end{subfigure}

\caption{Test case 2: computed velocity field $|\mathbf{v}_h|$ for the poroelastic case (left) and its vertical component $v_{h,y}$ (center) at time instant $t=0.6$ \si{\second}. Voronoi polygonal mesh (\#Elements = 2500) generated via \texttt{PolyMesher} algorithm \cite{Talischi2012} used in Test case 1, Test case 2 (right).}
\label{fig:Comparison_PE_test1}
\end{figure}

The aim of this subsection is to compare the results obtained via the thermo-poroelastic model, with the ones obtained through the poroelastic model presented in \cite{Antonietti2021}. In the poroelastic setting we consider the  parameters taken from Table~\ref{tab:TPEcoeff_test1} and the same mesh as the one displayed in Figure~\ref{fig:Comparison_PE_test1} (right). The yellow dot represents the point in which the forcing term $\mathbf{f}$ is located, while the red, green, and blue dots represent the points $\mathbf{x}_1 = (750 \si{\metre}, 1125 \si{\meter})$, $\mathbf{x}_2 = (1015 \si{\metre}, 1015 \si{\meter})$, and $\mathbf{x}_3 = (1125 \si{\metre}, 750 \si{\meter})$, respectively, where the solution is recorded. 
In Figure~\ref{fig:Comparison_PE_test1} we report the snapshots of the magnitude of the velocity field (left) and of its vertical component (center) computed for the poroelastic problem at the time instant $t = 0.6$ \si{\second}. As one can see, they are qualitatively similar to those reported in Figures~\ref{fig:v_test1}-\ref{fig:v2_test1}.

\begin{figure}[ht]
\begin{subfigure}[b]{.33\textwidth}
    \centering
    \includegraphics[width=0.9\textwidth]{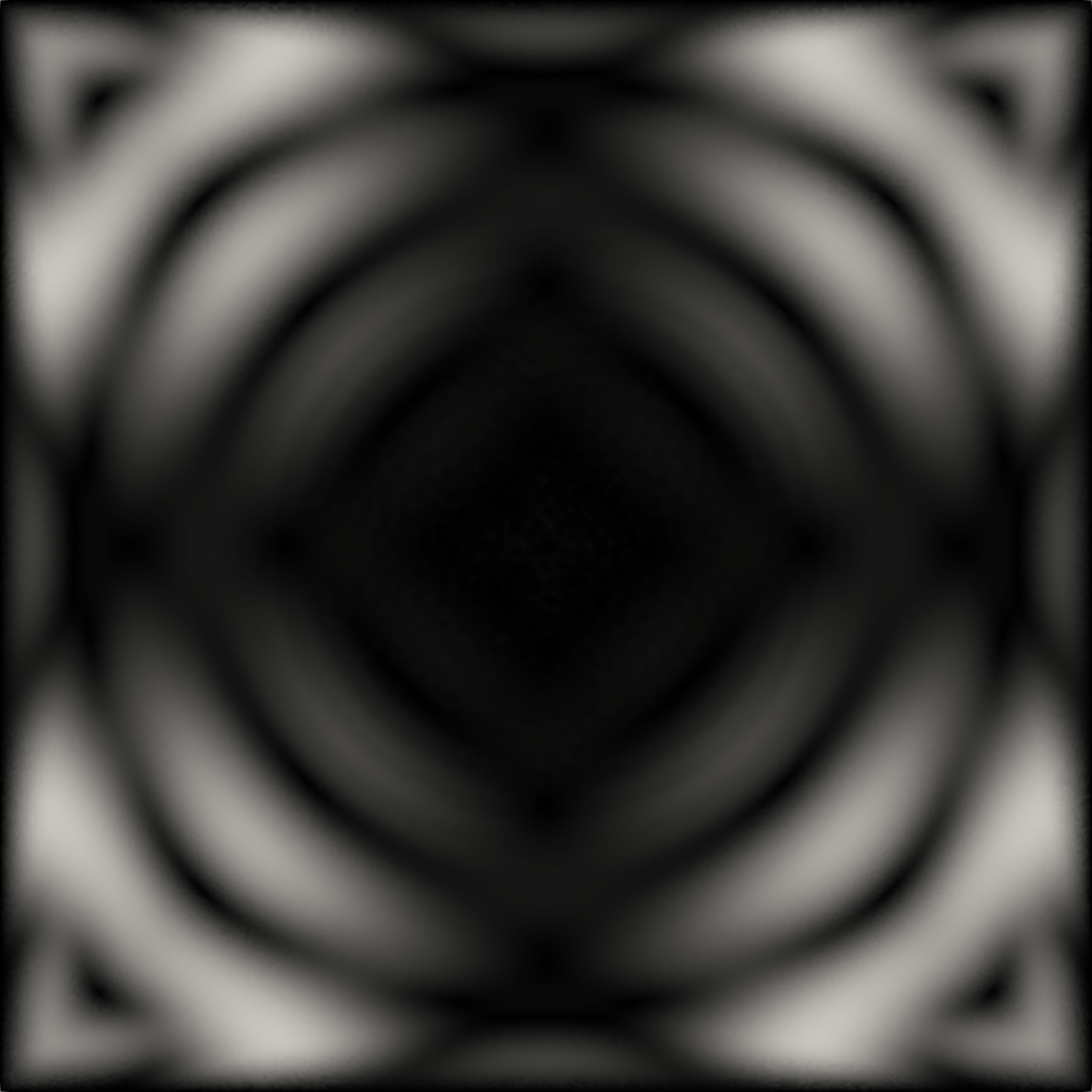}
    \label{fig:mag_err_v_test1}
\end{subfigure}
\begin{subfigure}[b]{.33\textwidth}
    \centering
    \includegraphics[width=0.9\textwidth]{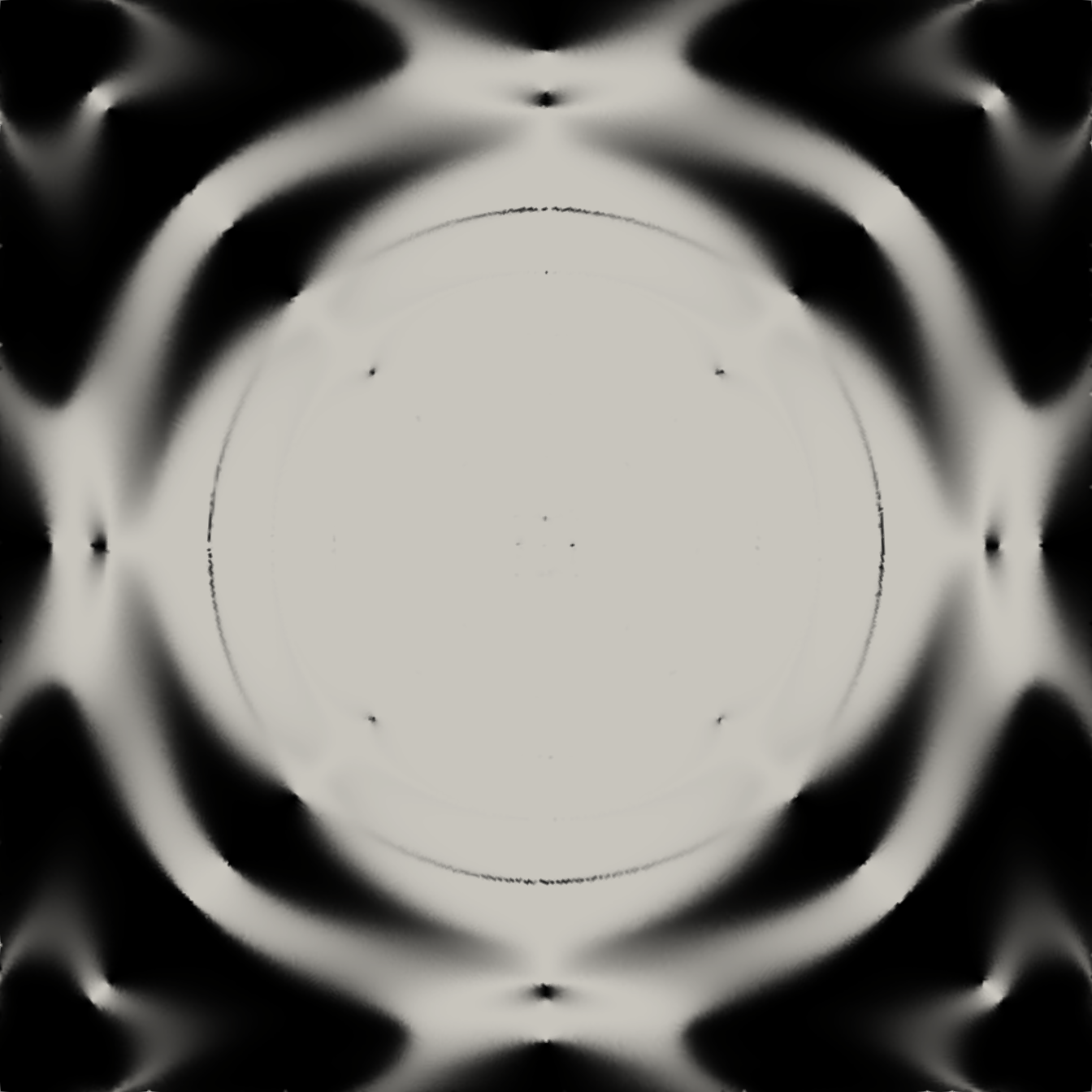}
    \label{fig:err_theta_v_test1}
\end{subfigure}
\begin{subfigure}[b]{.33\textwidth}
    \centering
    \includegraphics[width=0.9\textwidth]{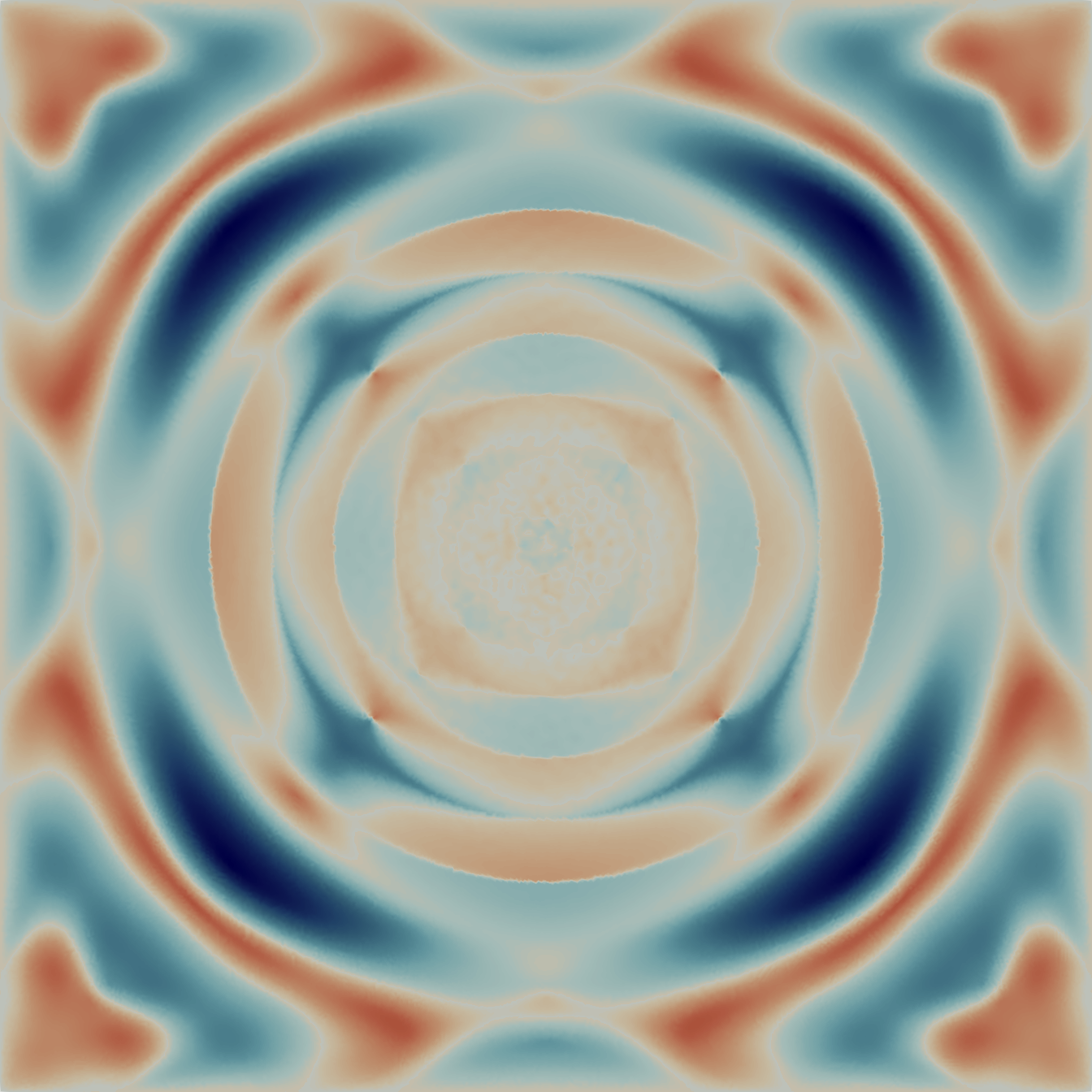}
    \label{fig:err_mag_v_test1}
\end{subfigure}

\begin{subfigure}[b]{.33\textwidth}
    \centering
    \includegraphics[width=0.9\textwidth]{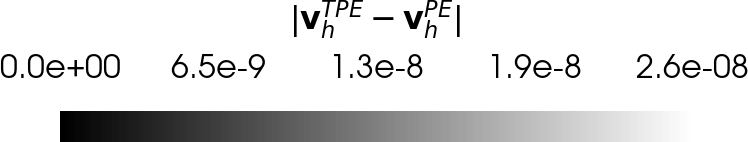}
    \caption{}
    \label{fig:mag_err_v_test1_colorbar}
\end{subfigure}
\begin{subfigure}[b]{.33\textwidth}
    \centering
    \includegraphics[width=0.9\textwidth]{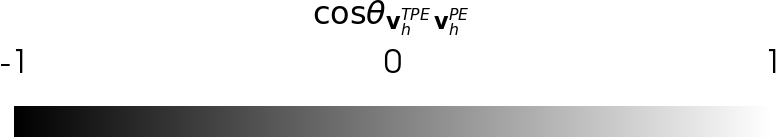}
    \caption{}
    \label{fig:err_theta_v_test1_colorbar}
\end{subfigure}
\begin{subfigure}[b]{.33\textwidth}
    \centering
    \includegraphics[width=0.9\textwidth]{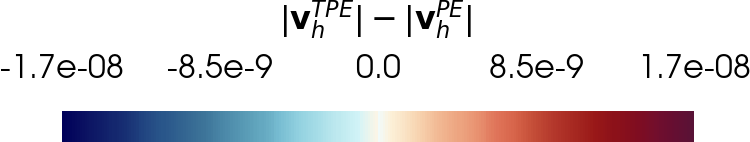}
    \caption{}
    \label{fig:err_mag_v_test1_colorbar}
\end{subfigure}

\caption{Test case 2: comparison of the velocity field $|\mathbf{v}_h|$ between the thermo-poroelastic (TPE) and poroelastic (PE) model in terms of magnitude of the difference (left), cosine of the subtended angle (center), and difference of magnitudes (right) at the time instant $t=0.6s$.}
\label{fig:diff_v_PE_test1}
\end{figure}

\begin{figure}[ht]
\begin{subfigure}[b]{.33\textwidth}
    \centering
    \includegraphics[width=0.9\textwidth]{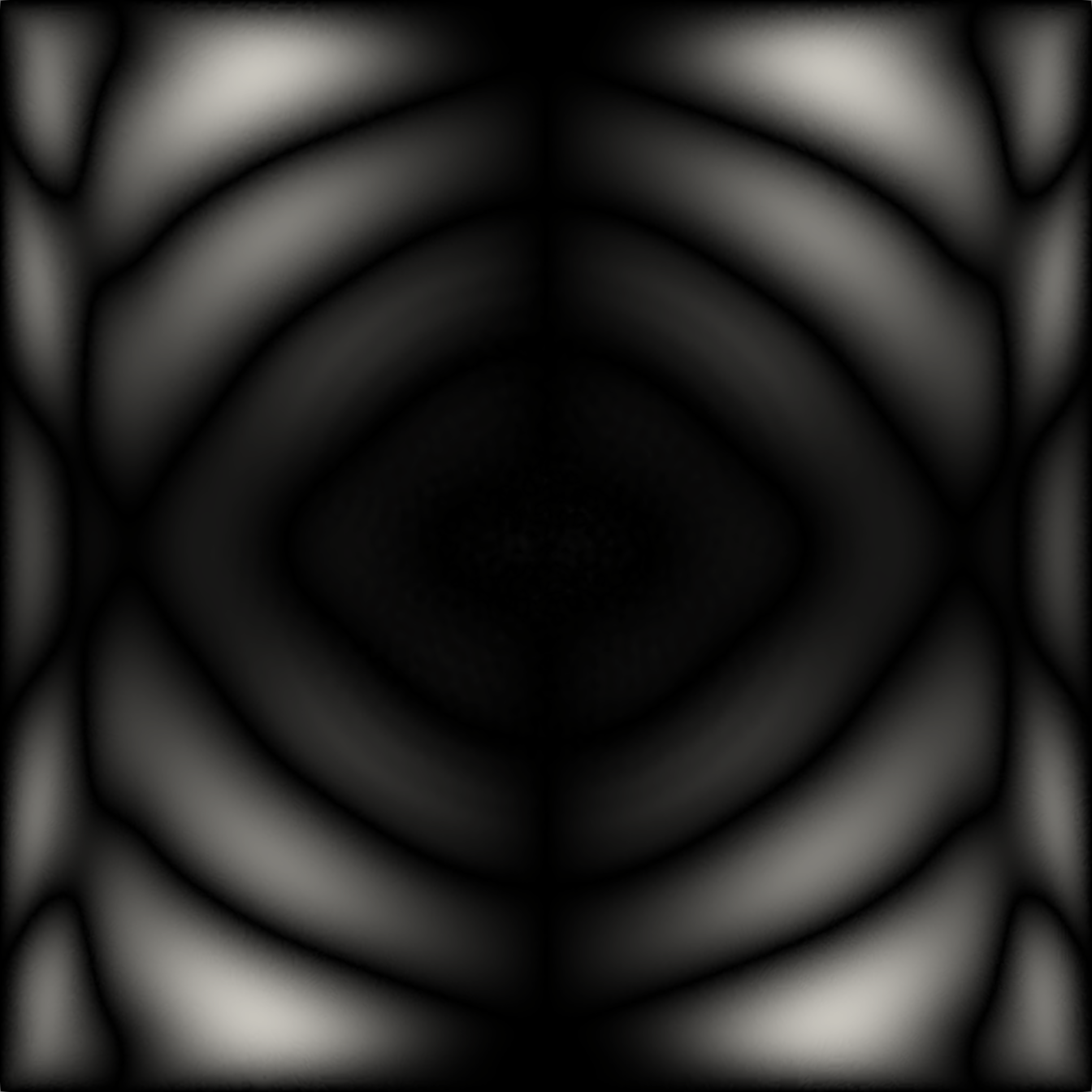}
    \label{fig:mag_err_v2_test1}
\end{subfigure}
\begin{subfigure}[b]{.33\textwidth}
    \centering
    \includegraphics[width=0.9\textwidth]{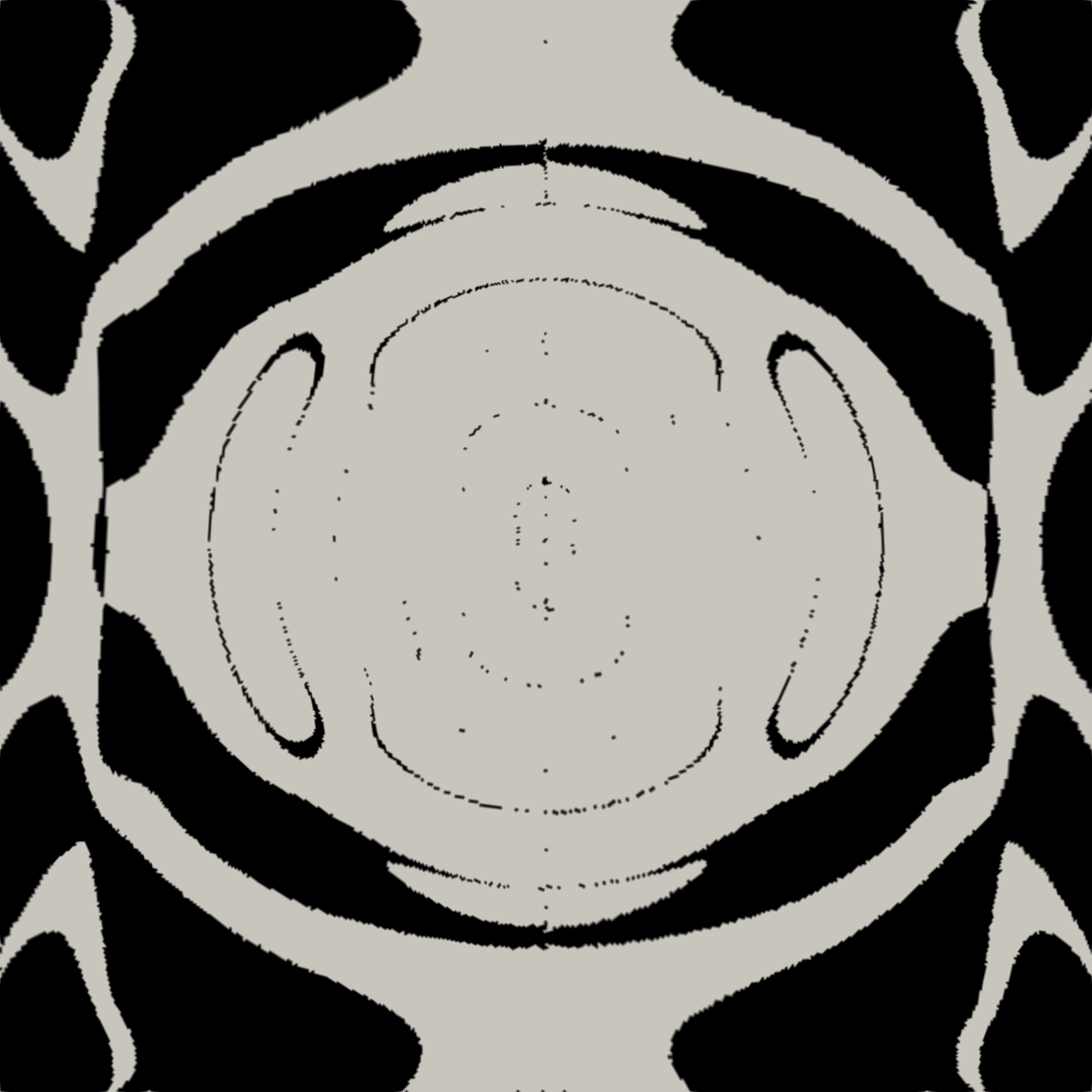}
    \label{fig:err_theta_v2_test1}
\end{subfigure}
\begin{subfigure}[b]{.33\textwidth}
    \centering
    \includegraphics[width=0.9\textwidth]{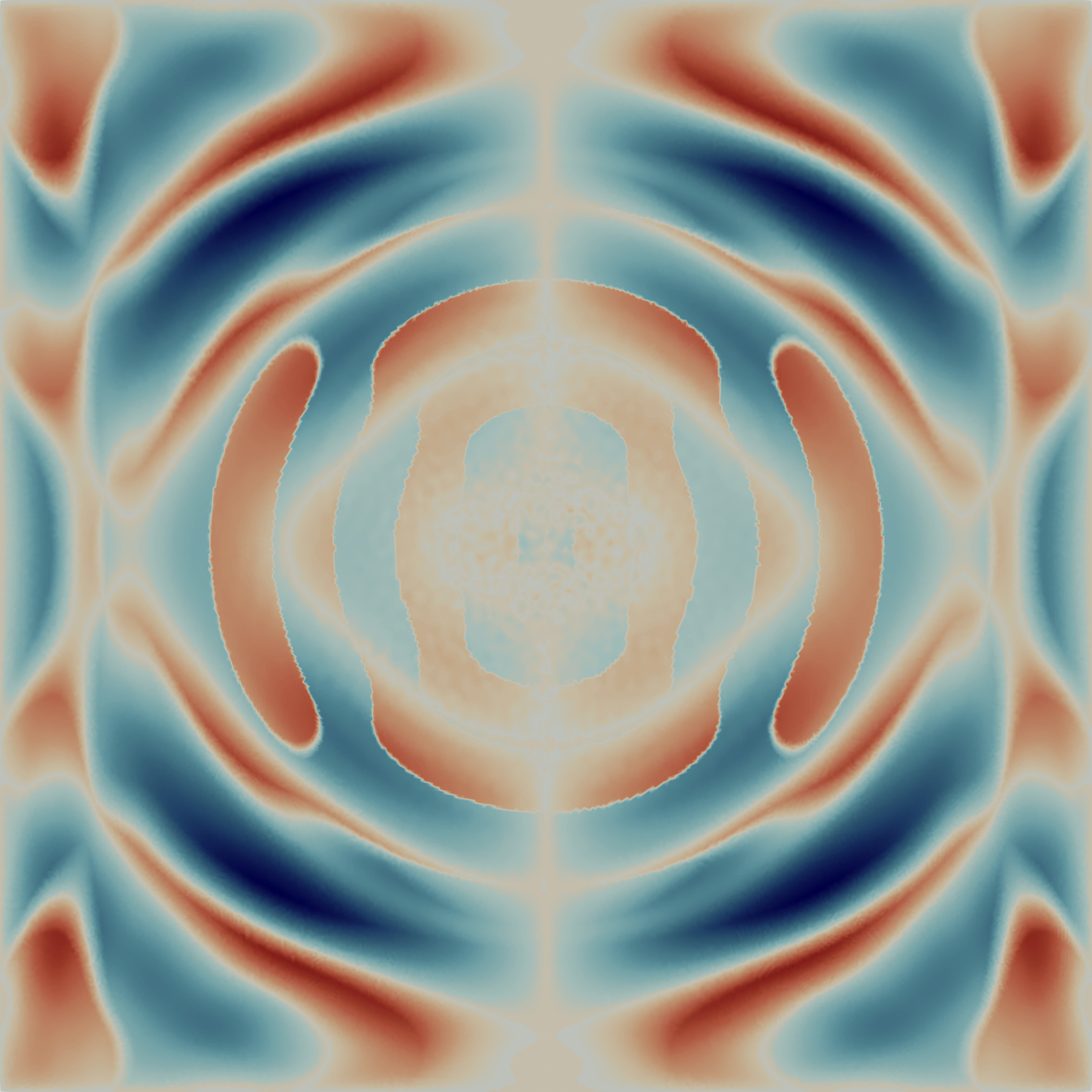}
    \label{fig:err_mag_v2_test1}
\end{subfigure}

\begin{subfigure}[b]{.33\textwidth}
    \centering
    \includegraphics[width=0.9\textwidth]{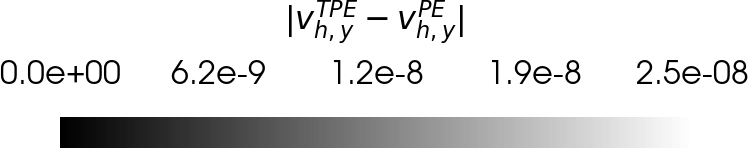}
    \caption{}
    \label{fig:mag_err_v2_test1_colorbar}
\end{subfigure}
\begin{subfigure}[b]{.33\textwidth}
    \centering
    \includegraphics[width=0.9\textwidth]{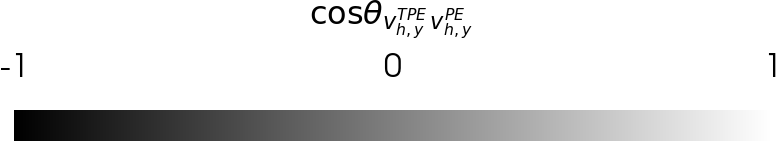}
    \caption{}
    \label{fig:err_theta_v2_test1_colorbar}
\end{subfigure}
\begin{subfigure}[b]{.33\textwidth}
    \centering
    \includegraphics[width=0.9\textwidth]{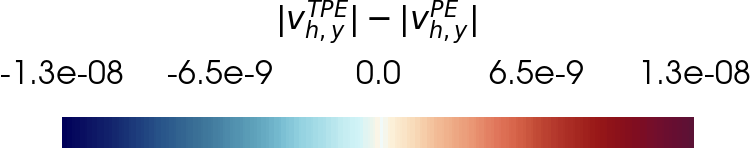}
    \caption{}
    \label{fig:err_mag_v2_test1_colorbar}
\end{subfigure}

\caption{Test case 2: comparison of the vertical component of the velocity field $v_{h,y}$ between the thermo-poroelastic (TPE) and poroelastic (PE) model in terms of magnitude of the difference (left), cosine of the subtended angle (center), and difference of magnitudes (right) at the time instant $t=0.6s$.}
\label{fig:diff_v2_PE_test1}
\end{figure}

In Figures~\ref{fig:diff_v_PE_test1} - \ref{fig:comparison_PE_q2_time} we compare qualitatively the solutions obtained with the thermo-poroelastic and the poroelastic model, respectively. Namely, in Figure~\ref{fig:diff_v_PE_test1} and Figure~\ref{fig:diff_v2_PE_test1} we plot the magnitude of the difference, the cosine of the subtended angle, and the difference of the magnitudes for the solid velocity and its vertical component, respectively. Looking at the six figures we can see that the patterns of differences regarding the entire velocity and its vertical component are similar. From Figure~\ref{fig:diff_v_PE_test1}, Figure~\ref{fig:diff_v2_PE_test1} we can see that the biggest differences can be observed in the areas where shear waves are present, while in the longitudinal and transverse areas of the domain, the behavior of the solutions is similar. Particularly interesting is the fact that, in the corners, the velocity fields computed by the thermo-poroelastic and the poroelastic models have opposite directions (cf. Figure~\ref{fig:err_theta_v_test1_colorbar} and Figure~\ref{fig:err_theta_v2_test1_colorbar}); this is due to the phase shift between the two wave fields. By looking at Figure~\ref{fig:err_mag_v_test1_colorbar} and Figure~\ref{fig:err_mag_v2_test1_colorbar} we can see that in the corners the magnitude of the velocity field generated by the poroelastic model is greater than the one of the thermo-poroelastic. Even in transversal and longitudinal directions, we can observe some slight differences. For instance, in Figure~\ref{fig:err_mag_v2_test1_colorbar} we see that in the propagating direction of the transversal waves, at the wave front with larger amplitude, that magnitude of the thermo-poroelastic wave is greater than the poroelastic one.
In general, we see that the difference is one order of magnitude lower with respect to the amplitude of the two waves.

\input{IMG/comparison_v.tikz}
\input{IMG/comparison_v2.tikz}
\input{IMG/comparison_q2.tikz}

In Figure~\ref{fig:comparison_PE_v_time}, \ref{fig:comparison_PE_v2_time}, \ref{fig:comparison_PE_q2_time} we report the time evolution of $\mathbf{v}_h, v_{h,y}$, and $q_{h,y}$  during time. First of all, we observe that for both of the cases, the amplitude of the filtration velocity has several orders of magnitude of difference with respect to the solid one; this is due to the choice of the forcing terms $\mathbf{f}, \mathbf{g}$ (cf. \cite{Morency2008}).

We see that for the monitored points  $\mathbf{x}_1$ and $\mathbf{x}_3$ (the ones on the vertical and horizontal direction, respectively) the results are comparable. Instead, for what concerns the point $\mathbf{x}_2$, we can observe considerable differences, not only a reduction of the amplitude of the wave, but also a phase shift. These remarks are in agreement with Figures~\ref{fig:diff_v_PE_test1}-\ref{fig:diff_v2_PE_test1} and motivate the major discrepancies on the corners of the domain.

To sum up, from the comparison of the results of the thermo-poroelastic model and the poroelastic one, it seems that adding the temperature to our model does not have a great effect on compressional waves, while it has an important impact on the computation of shear waves.

\subsection{Test case 3: heterogeneous media}
As a third test case, we consider wave propagation in a heterogeneous media. We split $\Omega = (-750, 750) \times (0, 1500) $ \si{\meter \squared} into two vertical layers. The left part of the domain is characterized by the same thermo-poroelastic properties of \textbf{Test case 1}, while in the right part, we consider the following, cf. Table~\ref{tab:TPEcoeff_test2} (the parameters that are not listed there are taken as in Table~\ref{tab:TPEcoeff_test1}). 
\begin{table}[H]
    \centering 
    \begin{tabular}{ l | l  c c l | l}
    \textbf{Coefficient} & \textbf{Value} & & & \textbf{Coefficient} & \textbf{Value} \T\B \\
    $a_0$ [\si[per-mode = symbol]{\pascal \per \kelvin \squared}] & 4.1017 & & & $\beta$ [\si[per-mode = symbol]{\pascal \per \kelvin}] & \num[exponent-product=\ensuremath{\cdot}]{4.8571e+4} \\
    $b_0$ [\si{\per \kelvin}] & \num[exponent-product=\ensuremath{\cdot}]{1.3684e-5} & & &  $\mu$ [\si{\pascal}] & \num[exponent-product=\ensuremath{\cdot}]{9e+9} \\
    $c_0$ [\si{\per\pascal}] & \num[exponent-product=\ensuremath{\cdot}]{1.3684e-10} & & &  $\lambda$ [\si{\pascal}] & \num[exponent-product=\ensuremath{\cdot}]{4e+9} \\
    $\alpha$ [-] & 0.7143 & & \\
    \end{tabular}
    \\[10pt]
    \caption{Test case 3: thermo-poroelastic properties of the medium (right layer)}
    \label{tab:TPEcoeff_test2}
\end{table}
The forcing terms, time-integration scheme, polynomial degree of approximation, and discretizations in space and time are the same of \textbf{Test case 1}.

\begin{figure}[ht]
\begin{subfigure}[b]{.33\textwidth}
    \centering
    \includegraphics[width=1\textwidth]{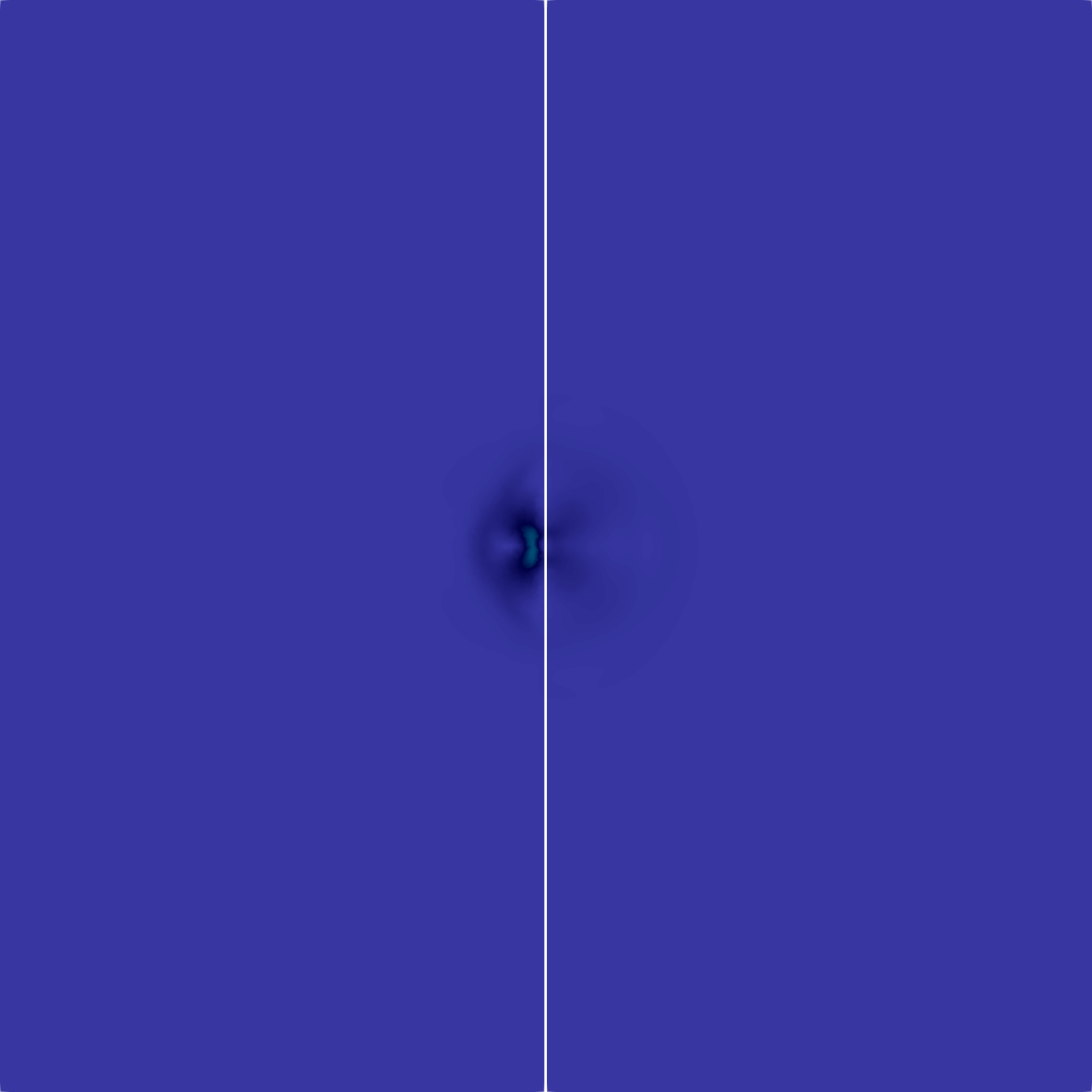}
    \label{fig:v_test2_01}
\end{subfigure}
\begin{subfigure}[b]{.33\textwidth}
    \centering
    \includegraphics[width=1\textwidth]{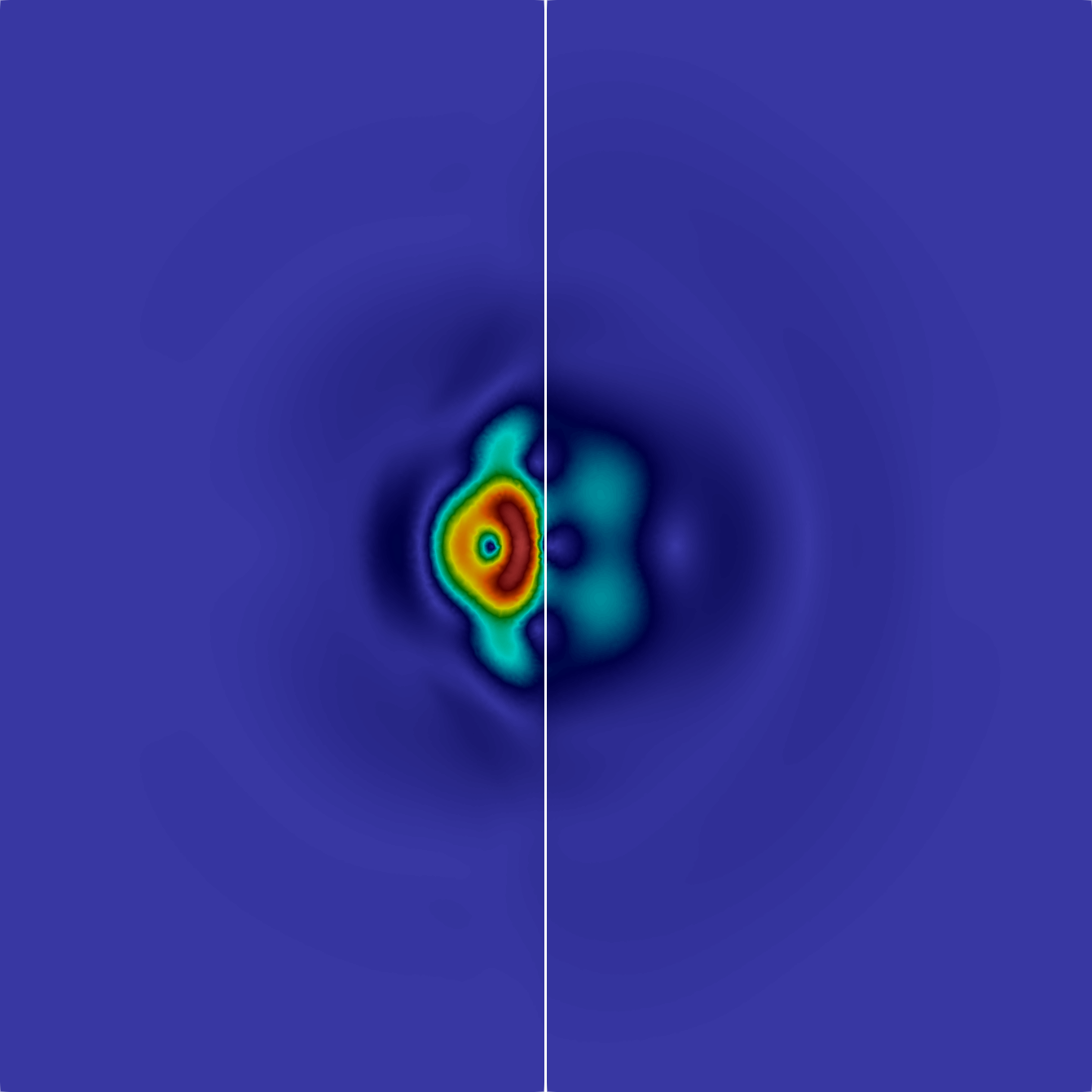}
    \label{fig:v_test2_03}
\end{subfigure}
\begin{subfigure}[b]{.33\textwidth}
    \centering
    \includegraphics[width=1\textwidth]{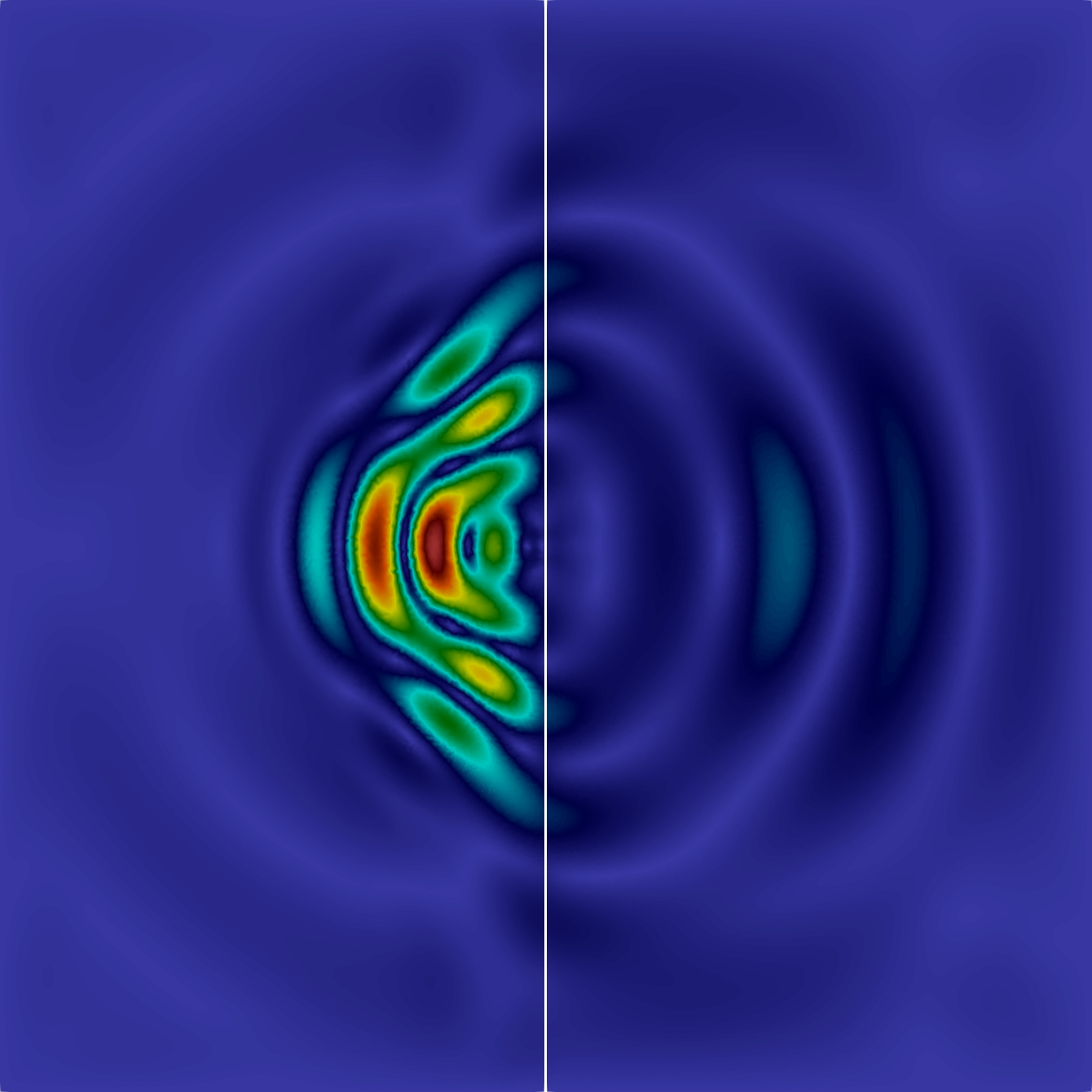}
    \label{fig:v_test2_05}
\end{subfigure}

\vspace{-0.4cm}
\centering
\includegraphics[width=0.4\textwidth]{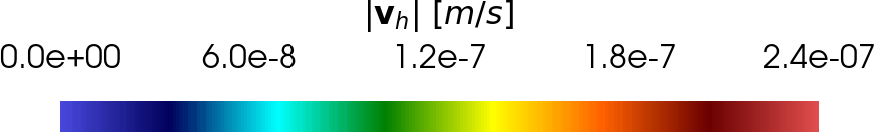}

\caption{Test case 3: computed velocity field $|\mathbf{v}_h|$ at the time instants $t=0.1s$ (left), $t=0.3s$ (center), $t=0.5s$ (right)}
\label{fig:v_test2}
\end{figure}

\begin{figure}[ht]
\begin{subfigure}[b]{.33\textwidth}
    \centering
    \includegraphics[width=1\textwidth]{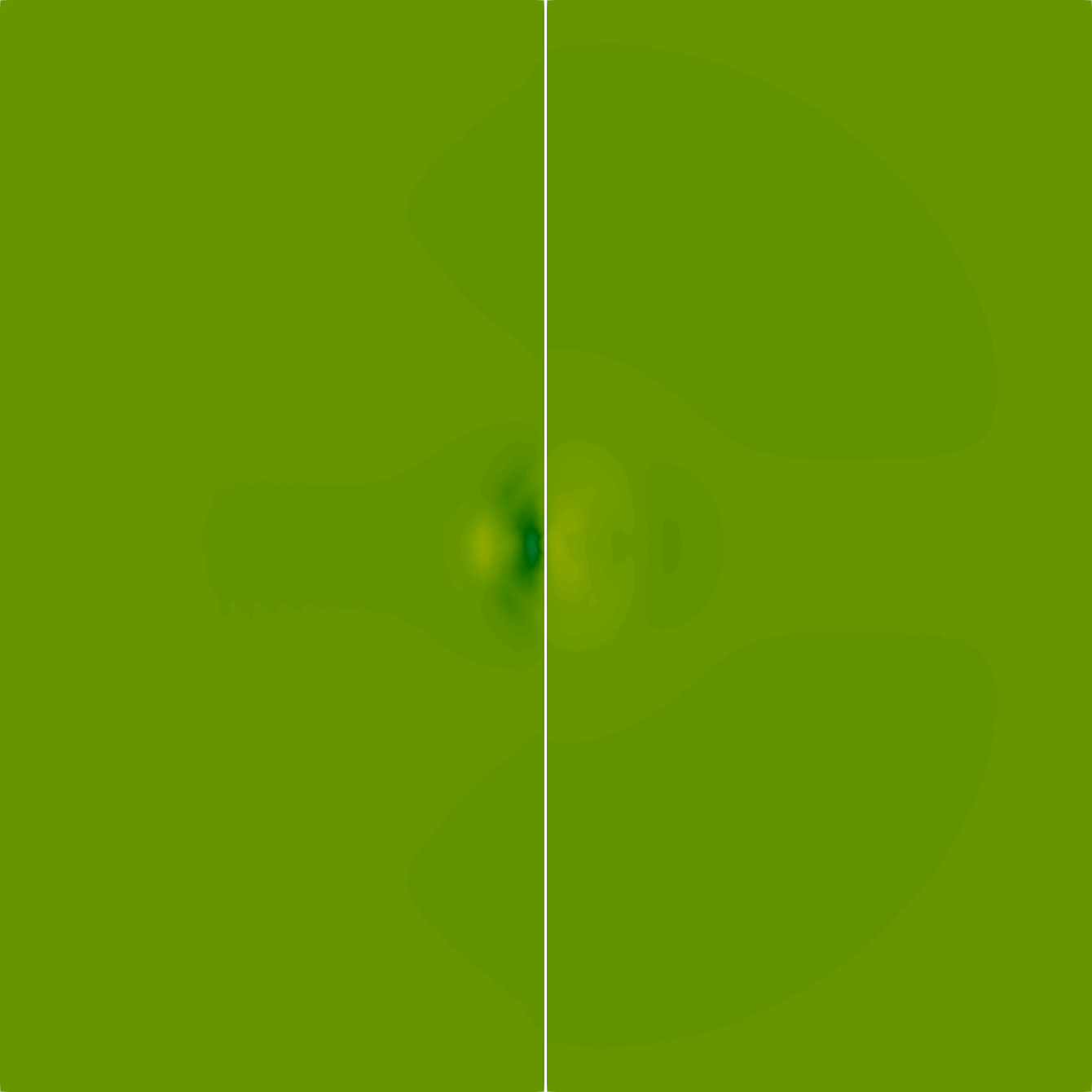}
    \label{fig:v2_test2_01}
\end{subfigure}
\begin{subfigure}[b]{.33\textwidth}
    \centering
    \includegraphics[width=1\textwidth]{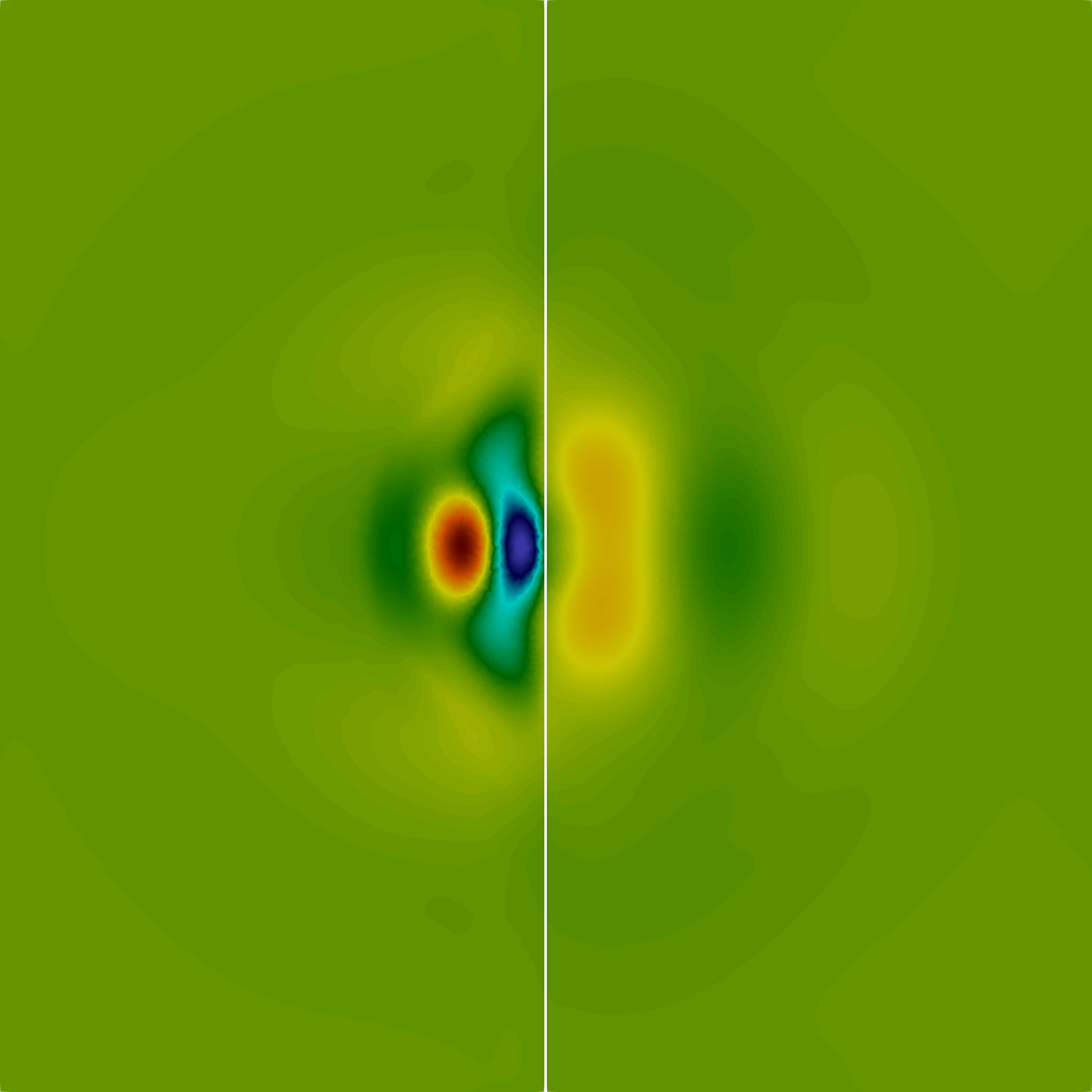}
    \label{fig:v2_test2_03}
\end{subfigure}
\begin{subfigure}[b]{.33\textwidth}
    \centering
    \includegraphics[width=1\textwidth]{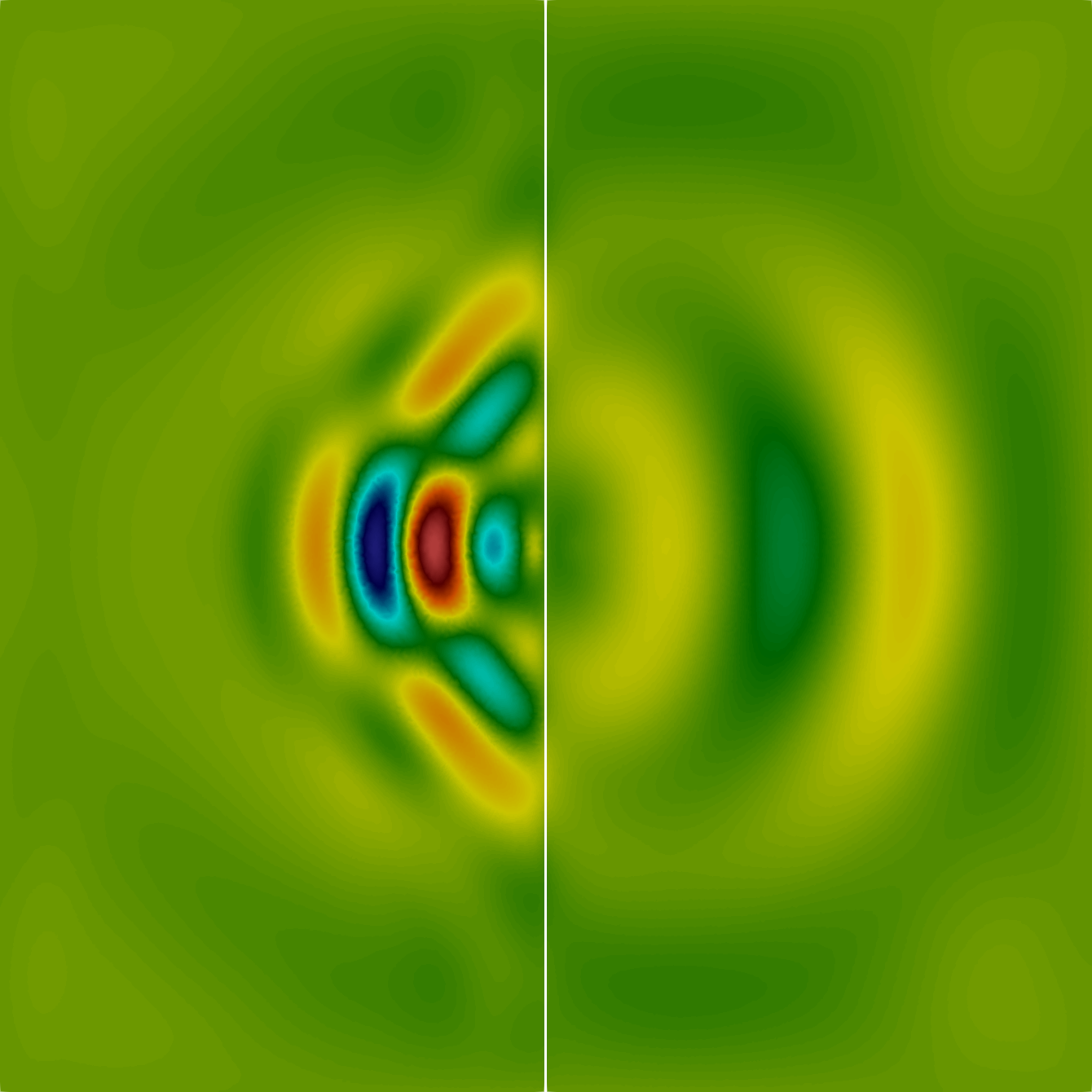}
    \label{fig:v2_test2_05}
\end{subfigure}

\vspace{-0.4cm}
\centering
\includegraphics[width=0.4\textwidth]{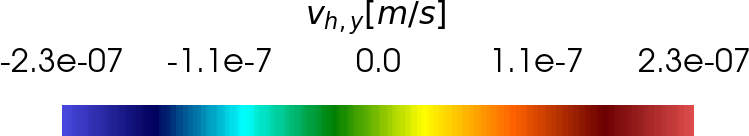}

\caption{Test case 3: computed vertical component of the velocity field $v_{h,y}$ at the time instants $t=0.1s$ (left), $t=0.3s$ (center), $t=0.5s$ (right)}
\label{fig:v2_test2}
\end{figure}

\begin{figure}[ht]
\begin{subfigure}[b]{.33\textwidth}
    \centering
    \includegraphics[width=1\textwidth]{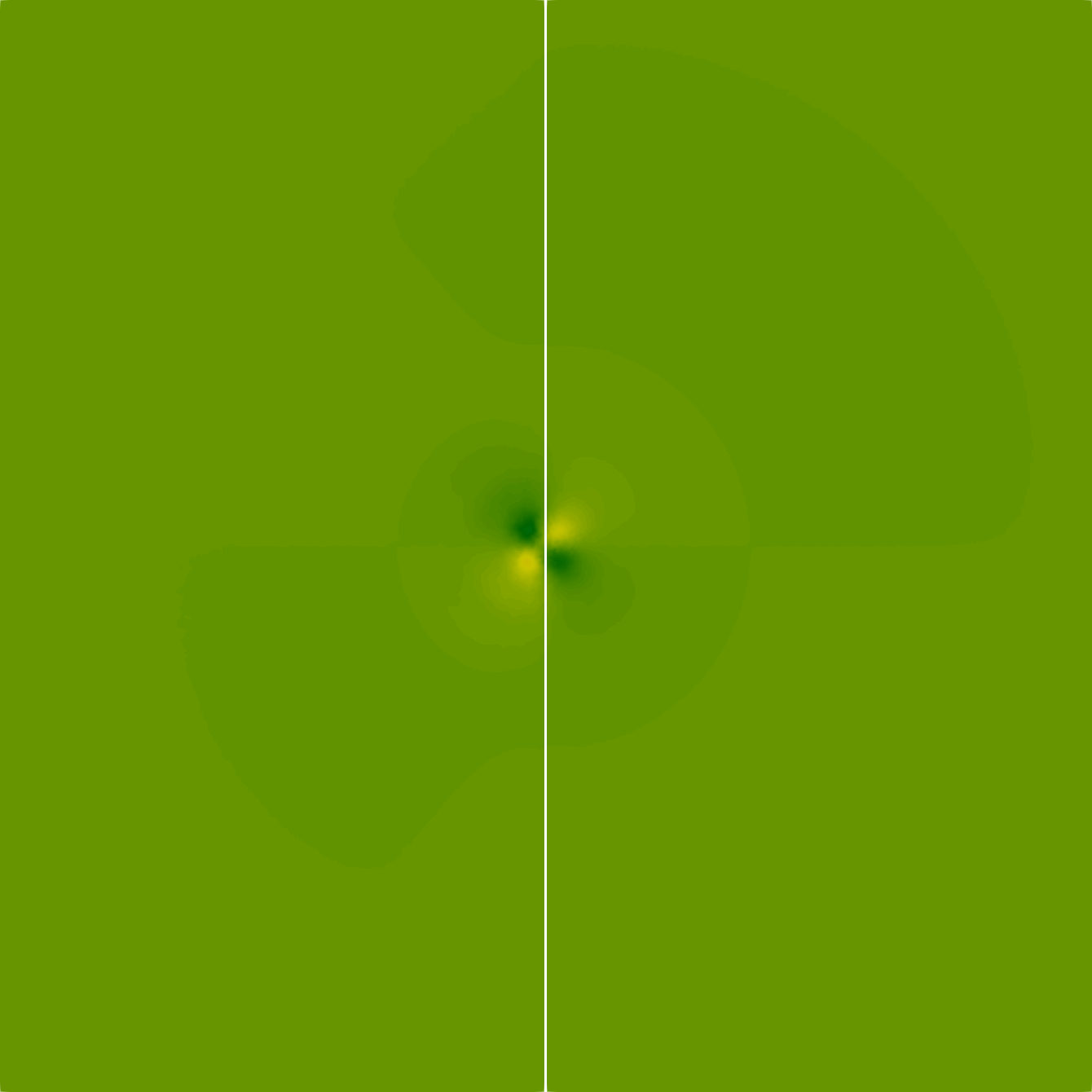}
    \label{fig:T_test2_01}
\end{subfigure}
\begin{subfigure}[b]{.33\textwidth}
    \centering
    \includegraphics[width=1\textwidth]{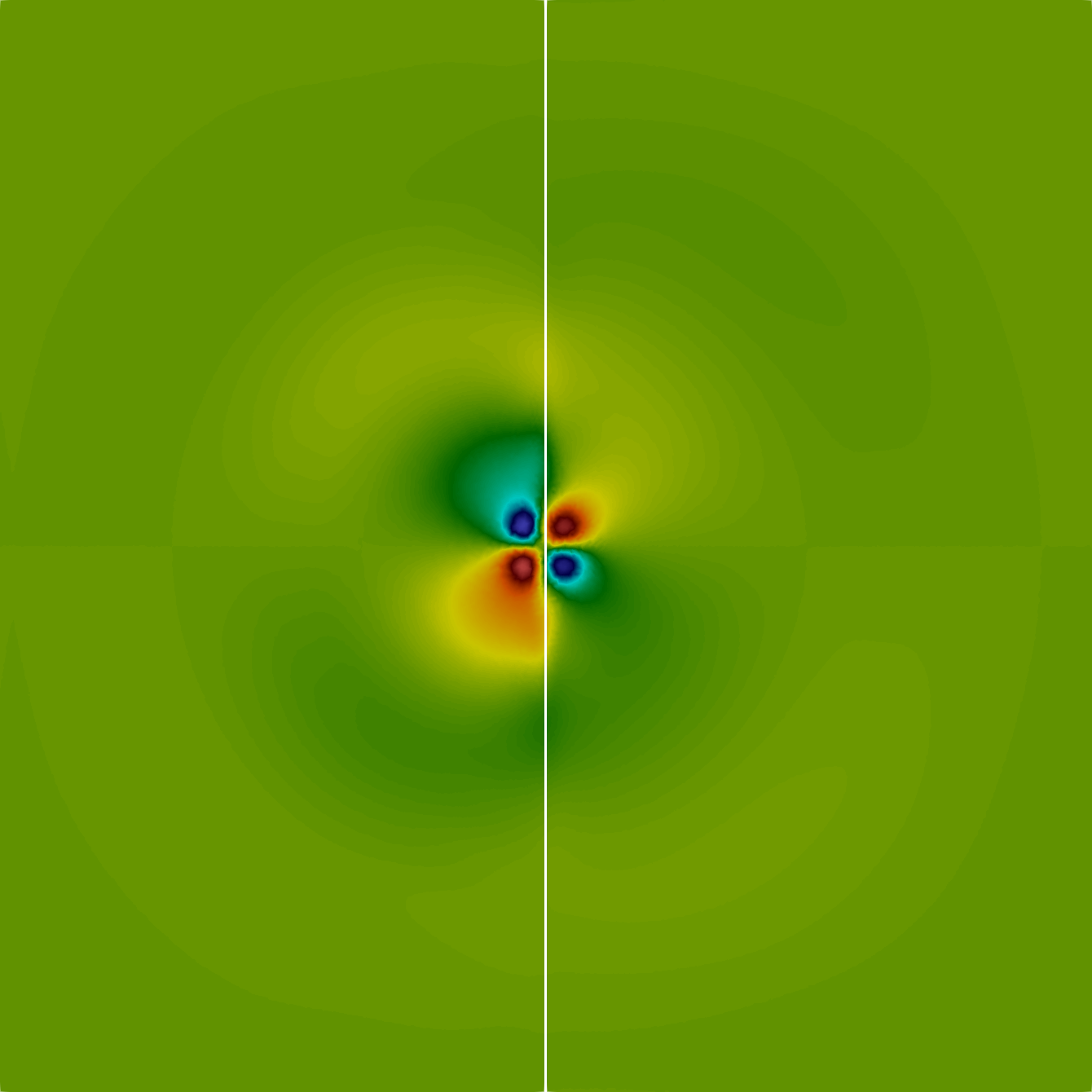}
    \label{fig:T_test2_03}
\end{subfigure}
\begin{subfigure}[b]{.33\textwidth}
    \centering
    \includegraphics[width=1\textwidth]{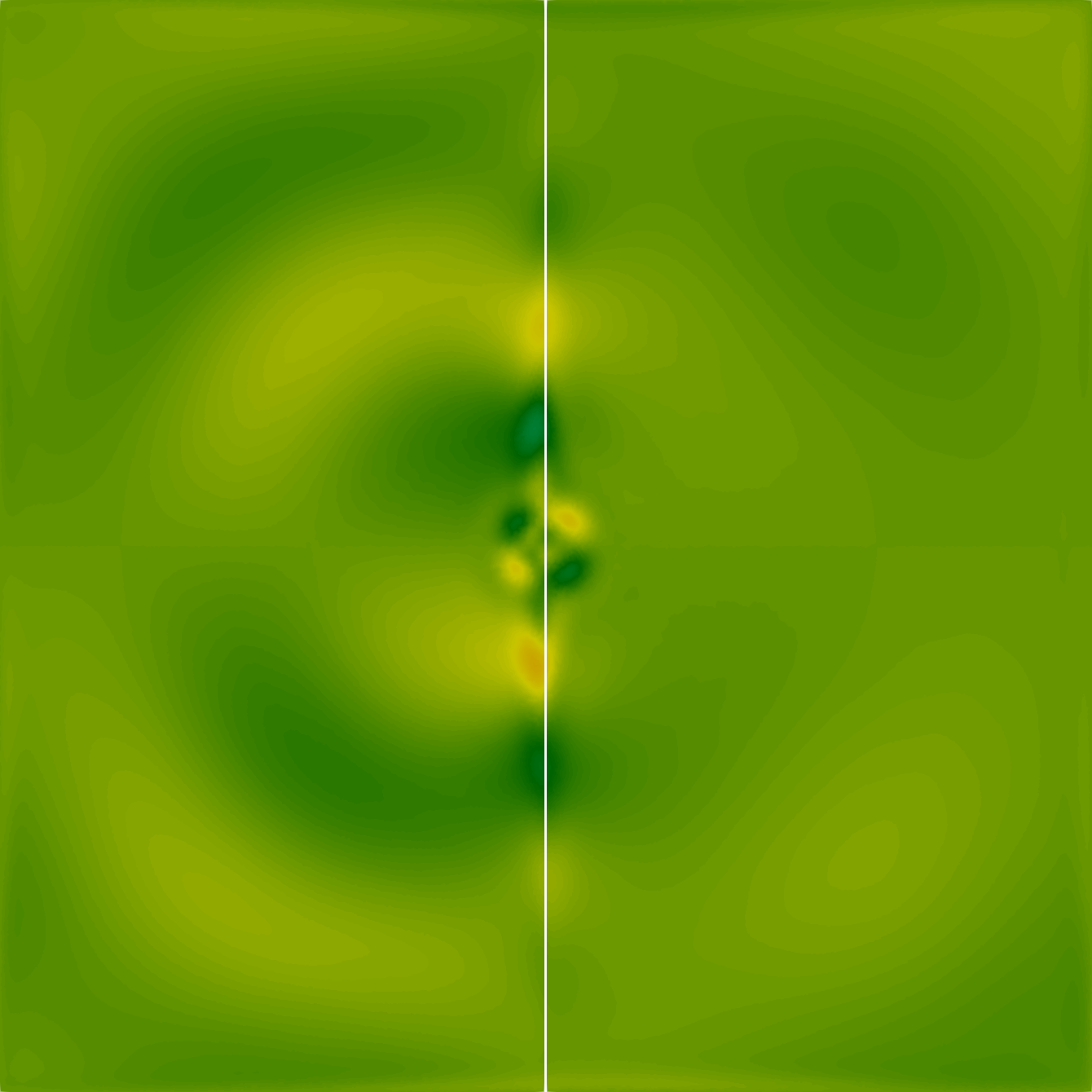}
    \label{fig:T_test2_05}
\end{subfigure}

\vspace{-0.4cm}
\centering
\includegraphics[width=0.4\textwidth]{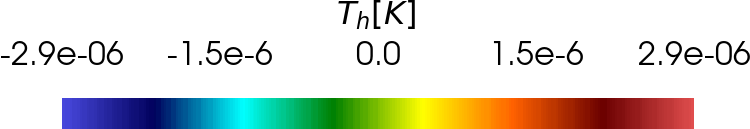}

\caption{Test case 3: computed temperature field $T_h$ at the time instants $t=0.1s$ (left), $t=0.3s$ (center), $t=0.5s$ (right)}
\label{fig:T_test2}
\end{figure}

The focus of this test case is to investigate how the heterogeneity of the media can affect wave propagation. In terms of the velocity field, the main difference with respect to the homogeneous case is the presence of the head waves, which are particularly evident by looking at the vertical component (cf. last frame of Figure~\ref{fig:v2_test2}). The field that is mainly affected by the change of medium is the temperature one. By looking at Figure~\ref{fig:T_test2} we can observe that the behavior of $T_h$ is quite different with respect to Figure~\ref{fig:T_test1}; moreover, even if the interface is very simple (i.e. a straight interface) we can observe the presence of strong boundary effects, notably stronger than the ones observed in Figure~\ref{fig:v_test2}, Figure~\ref{fig:v2_test2}. Also in this case, we can observe a qualitatively good agreement with respect to the results presented in \cite{Carcione2019}.

\section{Analysis of the semi-discrete problem}
\label{sec:semidiscrete_analysis}
The aim of this section is to provide the key ingredients and the main results regarding the analysis of the semi-discrete formulation \eqref{eq:semidiscr_form} in terms of stability and a-priori error estimates. The proofs of the stability and error estimates are reported in \ref{sec:staberr_proofs}.

In this section, we will focus on the case $\tau = 0$. We remark that, from a physical point of view, this assumption may lead to nonphysical results. However, the analysis carried out in this particular situation still can provide us crucial information about the performance of the method, even in the case $\tau \neq 0$, as numerically assessed in Section~\ref{sec:convergence_test}. As shown in \cite{Santos2021}, the analysis can be generalized to the case $\tau > 0$. The proof of the existence and uniqueness provided in the aforementioned work relies on the following steps: $(i)$ construct a sequence of approximate solutions of Problem~\ref{eq:TPE_system} using the Galerkin method; $(ii)$ derive a-priori bounds for the approximate solutions, in terms of boundary conditions, initial data, and forcing terms; $(iii)$ show the existence of the limit of a subsequence of approximate solutions in the weak$-^*$ topology via compactness arguments; and $(iv)$ show that the limit found satisfy the initial and boundary conditions for the problem. We think that, by combining these steps with arguments concerning the PolyDG-discretization, the analysis can be extended to he fully-hyperbolic model. The theoretical analysis of the full problem will be focus of future works.

In the case $\tau = 0$, problem \eqref{eq:semidiscr_form_tau} reduces to: \textit{for any $t \in (0, T_f]$, find $(\mathbf{u}_h, \mathbf{w}_h, T_h)(t) \in \mathbf{V}^{\ell}_h \times \mathbf{V}^{\ell}_h \times V^{\ell}_h $ such that $\forall \ (\mathbf{v}_h, \mathbf{z}_h, S_h) \in \mathbf{V}^{\ell}_h \times \mathbf{V}^{\ell}_h \times V^{\ell}_h $ it holds:}
\begin{equation}
\label{eq:semidiscr_form}
\begin{aligned}
& \mathcal{M}_{uw}((\ddot{\mathbf{u}}_h, \ddot{\mathbf{w}}_h), \left( \mathbf{v}_h, \mathbf{z}_h \right)) + \mathcal{B}( \dot{\mathbf{w}}_h, \mathbf{z}_h) +  \mathcal{M}_T(\dot{T}_h,S_h) + \mathcal{C}_h(\left(\dot{\mathbf{u}}_h, \dot{\mathbf{w}}_h \right), S_h) + \mathcal{A}_{uw,h}((\mathbf{u}_h, \mathbf{w}_h), \left( \mathbf{v}_h, \mathbf{z}_h \right)) \\
& + \mathcal{A}_{T,h}(T_h,S_h) - \mathcal{C}_h(\left( \mathbf{v}_h, \mathbf{z}_h \right), T_h) = ((\mathbf{f}, \mathbf{g}, H), (\mathbf{v}_h, \mathbf{z}_h, S_h)).
\end{aligned}
\end{equation}
We remark that the initial condition on $\dot{T}(0)$ is not needed for problem \eqref{eq:semidiscr_form}.

\subsection{Stability analysis}
\label{sec:stability_analysis}
In order to establish the stability of the proposed method, we mainly refer to the techniques that are used in \cite{Antonietti2021, AntoniettiBonetti2022}. 
We introduce the shorthand notation $||\cdot||_{\mathcal{F}}=\left(\sum_{F\in\mathcal{F}_h}||\cdot||_F^2\right)^{\frac12}$ and define the following $dG$-norms and semi-norms:
\begin{equation}
    \begin{aligned}
    \|\mathbf{v}\|_{dG,e}^2 & = \| \sqrt{2\mu} \ \boldsymbol{\epsilon}_h(\mathbf{v})||^2 + \| \sqrt{\sigma} \jump{\mathbf{v}} \ ||_{\mathcal{F}}^2 \ && \forall \ \mathbf{v} \in \mathbf{V}_h^{\ell}, \\[2pt]
    |\mathbf{z}|_{dG,p}^2 & = \| c_0^{-1} \, \divh{\mathbf{z}} \|^2 + \| \sqrt{\zeta} \jump{\mathbf{z}} \ ||_{\mathcal{F}}^2 \ && \forall \ \mathbf{z} \in \mathbf{V}_h^{\ell}, \\[2pt]
    \|S\|_{dG,T}^2 & = \|\sqrt{\theta}\nabla_h S\|^2 + \| \sqrt{\varrho} \jump{S} \ \|_{\mathcal{F}}^2 \ && \forall \ S \in V_h^{\ell}.
    \end{aligned}
\end{equation}
\begin{remark}
$| \cdot |_{dG,p} : \mathbf{V}_h^{\ell} \rightarrow \mathbb{R}^+$ is a semi-norm. By proceeding as in \cite{Antonietti2021} it is possible to prove that for any $\mathbf{v}$ and $ \mathbf{w}$ in $\mathbf{V}_h^{\ell} \times \mathbf{V}_h^{\ell}$, $\|\left( \mathbf{v}, \mathbf{w} \right)\|_{dG,*}^2 = \|\mathbf{v}\|_{dG,e}^2 + |\alpha \mathbf{v} + \mathbf{w}|_{dG,p}^2 + \mathcal{B}(\mathbf{w}, \mathbf{w})$ is a norm on $\mathbf{V}_h^{\ell} \times \mathbf{V}_h^{\ell}$.
\end{remark}
Then, we introduce two auxiliary norms for our problem for all $(\mathbf{v}, \mathbf{z}) \in C^1((0,T_f], \mathbf{V}_h^{\ell} \times \mathbf{V}_h^{\ell})$ and $S \in C^0((0,T_f], V_h^{\ell})$:
\begin{equation}
    \label{eq:energy_norm}
    \begin{aligned}
    & \|\left(\mathbf{v}, \mathbf{z}, S \right)(t)\|_{\mathcal{E}}^2 = \mathcal{M}_{uw}(\left(\dot{\mathbf{v}}, \dot{\mathbf{z}} \right), \left(\dot{\mathbf{v}}, \dot{\mathbf{z}} \right))(t) + \mathcal{M}_T(S,S)(t) +
    \|\left( \mathbf{v}, \mathbf{w} \right)(t)\|_{dG,*}^2 \\
    & \|\left(\mathbf{v}, \mathbf{z}, S \right)(t)\|_{\mathcal{E},*}^2 = \|\left(\mathbf{v}, \mathbf{z}, S \right)(t)\|_{\mathcal{E}}^2 + \int_0^t \| S(s) \|_{dG,T}^2 ds 
    \end{aligned}
\end{equation}
Now, we look at the boundedness and coercivity properties of the bilinear forms appearing in \eqref{eq:semidiscr_form}. Since we rely on standard steps in the discontinuous Galerkin framework, we refer to \cite[Section 3]{Antonietti2018}, \cite{Antonietti2021} for detailed proof.
\begin{lemma}
\label{lemma:cont_coer_stab}
Let Assumption~\ref{ass:c0} and Assumption~\ref{ass:mesh_Th1} be satisfied. Let $\alpha_1, \alpha_2, \alpha_3$ and $ \alpha_4$ in \eqref{eq:stabilization_func} be sufficiently large. Then, we have:
\begin{equation}
\begin{aligned}
    & \mathcal{A}_{e,h}(\mathbf{u}, \mathbf{v}) \lesssim \|\mathbf{u}\|_{dG,e} \|\mathbf{v}\|_{dG,e}, && \mathcal{A}_{e,h}(\mathbf{v}, \mathbf{v}) \gtrsim \|\mathbf{v}\|_{dG,e}^2 && \forall \ \mathbf{u}, \mathbf{v} \in \mathbf{V}_h^{\ell}, \\
    & \mathcal{A}_{T,h}(T, S) \lesssim \|T\|_{dG,T} \|S\|_{dG,T}, && \mathcal{A}_{T,h}(S,S) \gtrsim \|S\|_{dG,T}^2 && \forall \ T,S \in V_h^{\ell}, \\
    & \mathcal{M}_{uw}((\mathbf{u}, \mathbf{w}), \left( \mathbf{v}, \mathbf{z} \right)) \lesssim 
    \|(\mathbf{u},\mathbf{w})\| \|(\mathbf{v},\mathbf{z})\|, && \mathcal{M}_{uw}((\mathbf{v}, \mathbf{z}), \left( \mathbf{v}, \mathbf{z} \right)) \gtrsim 
    \|(\mathbf{v},\mathbf{z})\|^2 && \forall \ \mathbf{u}, \mathbf{w}, \mathbf{v}, \mathbf{z} \in \mathbf{V}_h^{\ell}, \\
    & \mathcal{M}_T(T,S) \lesssim 
    \|T\| \|S\|, && \mathcal{M}_T(S,S) \gtrsim 
    \|S\|^2 && \forall \ T,S \in V_h^{\ell}, \\
    \end{aligned}
\end{equation}
\vspace{-2mm}
\begin{equation}
    \begin{aligned}
    & \begin{aligned}
    \mathcal{A}_{uw,h}((\mathbf{u}, \mathbf{w}), ( \mathbf{v}, \mathbf{z})) + \mathcal{B}(\mathbf{w}, \mathbf{z}) \lesssim \|\left( \mathbf{u}, \mathbf{w} \right)\|_{dG,*} \ \|\left( \mathbf{v}, \mathbf{z} \right)\|_{dG,*}  
    \end{aligned} \\
    & \mathcal{A}_{uw,h}((\mathbf{v}, \mathbf{z}), ( \mathbf{v}, \mathbf{z})) + \mathcal{B}(\mathbf{z}, \mathbf{z}) \gtrsim \|\left( \mathbf{v}, \mathbf{z} \right)\|_{dG,*}^2
    \end{aligned} \hspace{2.37cm} \forall \ \mathbf{u}, \mathbf{w}, \mathbf{v}, \mathbf{z} \in \mathbf{V}_h^{\ell}.
\end{equation}
\end{lemma}
Note that, from Lemma~\ref{lemma:cont_coer_stab} it follows the well-posedness $\forall t \in (0, T_f]$ of Problem~\ref{eq:semidiscr_form}. We can state now the main Theorem of this section:

\begin{theorem}
\label{thm:stability_est}
Let Assumption~\ref{ass:c0} and Assumption~\ref{ass:mesh_Th1} hold, suppose that the parameters $\alpha_1$, $\alpha_2$, $\alpha_3$, and $\alpha_4$ appearing in \eqref{eq:stabilization_func} are large enough and let $\mathbf{X}_h = (\mathbf{u}_h, \mathbf{w}_h, T_h)(t)\in \mathbf{V}^{\ell}_h \times \mathbf{V}^{\ell}_h \times V^{\ell}_h$ be the solution of \eqref{eq:semidiscr_form} for any $t \in (0,T_f]$. Then, it holds
\begin{equation}
\label{eq:stability_est_thm}
\begin{aligned}
    \sup_{t\in(0,T_f]}\|\mathbf{X}_h(t)\|_{\mathcal{E},*}
    \lesssim & \ \|\mathbf{X}_h(0)\|_{\mathcal{E}} + \int_0^{T_f} \|  \left( \mathbf{f}, \mathbf{g}, H \right)(s) \| ds,
\end{aligned}
\end{equation}
where the hidden constant depends on the material properties and the final time $T_f$, but it does not depend on the mesh size $h$ and the polynomial degree $\ell$.
\end{theorem}

\textit{Proof.} The proof of Theorem~\ref{thm:stability_est} can be found in \ref{sec:stability_analysis_proof}.

\subsection{Error analysis}
\label{sec:error_analysis}
We start by defining, for a given element-wise constant $l:\mathcal{T}_h\to\mathbb{N}_{>0}$, the broken Sobolev spaces with variable regularity, which are needed to establish the error bounds in the $hp$-framework. We set
\begin{equation}
   \begin{aligned}
   H^{l}(\mathcal{T}_h) & \ = \left\{ v_h \in L^2(\Omega) : v_h |_{\kappa} \in H^{l_{\kappa}}(\kappa) \ \ \forall \kappa \in \mathcal{T}_h \right\}, \quad \mathbf{H}^{l}(\mathcal{T}_h) = \left[ H^{l}(\mathcal{T}_h) \right]^d.
   \end{aligned}
\end{equation}
where as usual $l_{\kappa}=l_{|\kappa}$. Next, we introduce the \textit{stronger} $dG$-norms
\begin{equation}
    \label{eq:DG_triple_norms}
    \begin{aligned}
    & |||\mathbf{v}|||^2_{dG,e} = ||\mathbf{v}||^2_{dG,e} +  ||\sigma^{-\frac12}\avg{2\mu\boldsymbol{\epsilon}_h(\mathbf{v})}||_{\mathcal{F}}^2 + \|\xi^{-\frac12} \avg{\lambda\divh{\mathbf{v}}}\|_{\mathcal{F}}^2
    \ && \forall \ \mathbf{v} \in \mathbf{H}^2(\mathcal{T}_h), \\
    & |||\mathbf{z}|||^2_{dG,p} = |\mathbf{z}|^2_{dG,p} +  ||\zeta^{-\frac12}\avg{{c_0}^{-1} \divh{\mathbf{z}}}||_{\mathcal{F}}^2 \ && \forall \ \mathbf{z} \in \mathbf{H}^2(\mathcal{T}_h), \\
    &|||S|||^2_{dG,T} = ||S||^2_{dG,T} + ||\varrho^{-\frac12} \avg{ \theta \, \nabla_h S}||_{\mathcal{F}}^2 \ &&  \forall \ S \in H^2({\mathcal{T}_h}).\\ 
    \end{aligned}
\end{equation}
and the following Lemma stating the boundedness of the bilinear forms \eqref{eq:bil_form_discr} in the dG-norms \eqref{eq:DG_triple_norms}:
\begin{lemma}
\label{lemma:bound_bilinearform_error}
Let Assumption~\ref{ass:c0} and Assumption~\ref{ass:mesh_Th1} be satisfied.  Let $\alpha_1, \alpha_2, \alpha_3, \alpha_4$ in \eqref{eq:stabilization_func} be sufficiently large. Then, 
\begin{equation}
    \begin{aligned}
    & \mathcal{A}_{T,h}(T,S) \lesssim \ |||T|||_{dG,T} ||S||_{dG,T} \quad && \forall \ T \in H^2(\mathcal{T}_h), \ \forall \ S \in V_h^{\ell},\\
    & \mathcal{A}_{e,h}(\mathbf{u},\mathbf{v}) \lesssim \ |||\mathbf{u}|||_{dG,e} ||\mathbf{v}||_{dG,e} \quad && \forall \ \mathbf{u} \in \mathbf{H}^2(\mathcal{T}_h),\ \forall \ \mathbf{v} \in \mathbf{V}_h^{\ell},\\
    & \mathcal{A}_{p,h}(\mathbf{w},\mathbf{z}) \lesssim \ |||\mathbf{w}|||_{dG,p} |\mathbf{z}|_{dG,p} \quad && \forall \ \mathbf{w} \in \mathbf{H}^2(\mathcal{T}_h),\ \forall \ \mathbf{z} \in \mathbf{V}_h^{\ell},\\
    & \mathcal{C}_{h}((\mathbf{u},\mathbf{w}),S) \lesssim \ \big( |||\mathbf{u}|||_{dG,e} + |||\mathbf{w}|||_{dG,p} \big) \| S \|_{dG,T} \quad && \forall \ \mathbf{u}, \mathbf{w} \in \mathbf{H}^2(\mathcal{T}_h),\ \forall \ S \in V_h^{\ell},\\
    & \mathcal{C}_{h}((\mathbf{v},\mathbf{z}),T) \lesssim \ \big( \| \mathbf{v} \|_{dG,e} + | \mathbf{z} |_{dG,p} \big) |||T|||_{dG,T} \quad && \forall \ T \in H^2(\mathcal{T}_h),\ \forall \ \mathbf{v}, \mathbf{z} \in \mathbf{V}_h^{\ell},\\
    & \begin{aligned}
    \mathcal{A}_{uw,h}((\mathbf{u},\mathbf{w}),(\mathbf{v},\mathbf{z})) \lesssim & \ |||\mathbf{u}|||_{dG,e}^2 +  |||\mathbf{u}|||_{dG,p}^2 + |||\mathbf{w}|||_{dG,p}^2 \\
    & + \|\mathbf{v}\|_{dG,e}^2 + |\mathbf{v}|_{dG,p}^2 +  |\mathbf{z}|_{dG,p}^2 \\  
    \end{aligned} \quad && \forall \ \mathbf{u}, \mathbf{w} \in \mathbf{H}^2(\mathcal{T}_h),\ \forall \ \mathbf{v}, \mathbf{z} \in \mathbf{V}_h^{\ell}.\\
    \end{aligned}
\end{equation}
\end{lemma}
For the sake of readability, we also define
\begin{equation}
\begin{aligned}
& |||(\mathbf{v}, \mathbf{z})(t)|||_{dG,*}^2 = |||\mathbf{v}(t)|||_{dG,e}^2 + |||\mathbf{z}(t)|||_{dG,p}^2, \\
& |||\left(\mathbf{v}, \mathbf{z}, S \right)(t)|||_{dG}^2 = |||\mathbf{v}(t)|||_{dG,e}^2 + |||\mathbf{v}(t)|||_{dG,p}^2 + |||\mathbf{z}(t)|||_{dG,p}^2 + |||S(t)|||_{dG,T}^2, \\ 
& |||\left(\mathbf{v}, \mathbf{z}, S \right)(t)|||_{\mathcal{E}}^2 = \mathcal{M}_e(\left(\dot{\mathbf{v}}, \dot{\mathbf{z}} \right), \left(\dot{\mathbf{v}}, \dot{\mathbf{z}} \right))(t) + \mathcal{M}_T(S,S)(t) + \mathcal{B}(\mathbf{z}, \mathbf{z})(t) + |||\left(\mathbf{v}, \mathbf{z}, S \right)(t)|||_{dG}^2.
\end{aligned}
\end{equation}
Then, we need to introduce the interpolants $\left( \mathbf{u}_I, \mathbf{
w}_I, T_I \right)$ of the solutions to \eqref{eq:DYN_weak_form}. In order to properly treat the interpolation errors, we introduce the Stein extension operator. For a polytopic mesh $\mathcal{T}_h$ satisfying Assumption~\ref{ass:A2}, the Stein operator $\mathcal{E}: H^m(\kappa) \rightarrow H^m(\mathbb{R}^d)$ is defined for any $\kappa \in \mathcal{T}_h$ and $m \in \mathbb{N}_{>0}$ such that
\begin{equation}
    \label{eq:stein_operator}
    \mathcal{E}v \rvert_{\kappa} = v, \quad ||\mathcal{E}v||_{H^m(\mathbb{R}^d)} \lesssim ||v||_{H^m(\kappa)} \qquad \forall v \in H^m(\kappa).  
\end{equation}
Analogously, we can define the vector-valued version that acts component-wise and is denoted in the same way. In what follows, for any $\kappa\in\mathcal{T}_h$, we will denote by $\mathcal{K}_\kappa$ the simplex belonging to $\mathcal{T}_h^*$ such that $\kappa\subset\mathcal{K}_\kappa$. 
Then, we state the following approximation estimate: 
\begin{lemma}
    \label{lemma:interp}
    Let Assumption~\ref{ass:mesh_Th1} be fulfilled. For any $\mathbf{v} \in \mathbf{H}^l(\mathcal{T}_h)$, $\mathbf{z} \in \mathbf{H}^m(\mathcal{T}_h)$, and $S \in H^n(\mathcal{T}_h)$ with $l_{\kappa}, m_{\kappa}, n_{\kappa} \geq 2$ for all $\kappa\in\mathcal{T}_h$, there exist $\mathbf{v}_I \in \mathbf{V}_h^{\ell}$, $\mathbf{z}_I \in \mathbf{V}_h^{\ell}$, and $S_I \in V_h^{\ell}$ such that
    \begin{equation}
        \begin{aligned}
        & |||\mathbf{v} - \mathbf{v}_I|||_{dG,e}^2 \lesssim \sum_{\kappa \in \mathcal{T}_h} \frac{h_{\kappa}^{2(q_{\kappa} - 1)}}{\ell_{\kappa}^{2l_{\kappa}-3}} \ ||\mathcal{E}\mathbf{v}||_{\mathbf{H}^{l_{\kappa}}(\mathcal{K}_\kappa)}^2 \\
        & |||\mathbf{z} - \mathbf{z}_I|||_{dG,p}^2 \lesssim \sum_{\kappa \in \mathcal{T}_h} \frac{h_{\kappa}^{2(r_{\kappa} - 1)}}{\ell_{\kappa}^{2m_{\kappa}-3}} \ ||\mathcal{E}\mathbf{z}||_{\mathbf{H}^{m_{\kappa}}(\mathcal{K}_\kappa)}^2 \\
        & |||S - S_I|||_{dG,T}^2 \lesssim \sum_{\kappa \in \mathcal{T}_h} \frac{h_{\kappa}^{2(s_{\kappa} - 1)}}{\ell_{\kappa}^{2n_{\kappa}-3}} \ ||\mathcal{E}S||_{H^{n_{\kappa}}(\mathcal{K}_\kappa)}^2 \\
        \end{aligned}
    \end{equation}
    where $q_{\kappa} = \min\{\ell_{\kappa}+1,l_{\kappa}\}$, $r_{\kappa} = \min\{\ell_{\kappa}+1,m_{\kappa}\}$, and $s_{\kappa} = \min\{\ell_{\kappa}+1,n_{\kappa}\}$. 
\end{lemma} 
For the detailed proof of Lemma~\ref{lemma:interp} refer to \cite[Lemma 3.6]{Antonietti2018}, \cite[Theorem 36]{Cangiani2017}, and \cite[Corollary 5.1]{antonietti2020_stokesDG}. Last, we need to introduce the notion of error that is used in the analysis. We consider the discretization errors $E = (\mathbf{e}^u, \mathbf{e}^w, e^T)$, where $\textbf{e}^u(t) = \mathbf{u}(t) - \mathbf{u}_h(t), \quad \textbf{e}^w(t) = \mathbf{w}(t) - \mathbf{w}_h(t), \quad e^T(t) = T(t) - T_h(t).$
Denoting by $\mathbf{X}_h(t) = (\mathbf{u}_h, \mathbf{w}_h, T_h)(t) \in \mathbf{V}_h^{\ell} \times \mathbf{V}_h^{\ell} \times V_h^{\ell}$ and $\mathbf{X}(t)=(\mathbf{u}, \mathbf{w}, T)(t) \in  \mathbf{V} \times \mathbf{W} \times V$ for all $t\in (0,T_f]$ the solutions to \eqref{eq:semidiscr_form} and \eqref{eq:DYN_weak_form} (with $\tau = 0$), respectively, we split the errors as $E(t) = E_I(t)-E_h(t)$, with 
\begin{align}
    E_I(t) & = \mathbf{X}(t)-\mathbf{X}_I(t) = (\textbf{e}_I^u(t),\textbf{e}_I^w(t),e_I^T(t))=(\mathbf{u}(t) - \mathbf{u}_I(t), \mathbf{w}(t) - \mathbf{w}_I(t),T(t) - T_I(t))  \\
    E_h(t) & = \mathbf{X}_I(t)-\mathbf{X}_h(t) = (\textbf{e}_h^u(t),\textbf{e}_h^w(t),e_h^T(t))=(\mathbf{u}_I(t) - \mathbf{u}_h(t), \mathbf{w}_I(t) - \mathbf{w}_h(t),T_I(t) - T_h(t)).
\end{align}

Exploiting all of the previous ingredients, we can now state the main result of this section.
\begin{theorem}
\label{thm:error_est}
Let Assumption~\ref{ass:c0} and Assumption~\ref{ass:mesh_Th1} be valid and assume that the parameters $\alpha_1$, $\alpha_2$, $\alpha_3$, and $\alpha_4$ appearing in \eqref{eq:stabilization_func} are large enough. Let the exact solutions of problem \eqref{eq:TPE_system} be such that
\begin{equation}
\begin{aligned}
(\mathbf{u}, \mathbf{w}) \in C^2((0,T_f]; \mathbf{H}^l(\mathcal{T}_h) \times \mathbf{H}^m(\mathcal{T}_h)) \, \cap \, C^1((0,T_f]; \mathbf{V} \times \mathbf{W}), \quad T \in C^1((0,T_f]; H^n(\mathcal{T}_h)\cap V)
\end{aligned}
\end{equation}
with $l, m, n \geq 2$, and let $(\mathbf{u}_h, \mathbf{w}_h) \in C^2((0,T_f]; \mathbf{V}_{h}^{\ell} \times \mathbf{V}_{h}^{\ell})$, $T_h \in C^1((0,T_f]; V_{h}^{\ell})$ be the solutions of the semi-discrete problem \eqref{eq:semidiscr_form}. Then, for all $t \in (0,T_f]$ the discretization error $E_h = (\mathbf{e}_h^u, \mathbf{e}_h^w, e_h^T)$ satisfies
\begin{equation}
    \label{eq:error_est}
    \begin{aligned}
    \|E_h(t)\|_{\mathcal{E}}^2 + \int_0^t \|e_h^T(s)\|_{dG,T}^2 \lesssim
    \ & \sum_{\kappa \in \mathcal{T}_h} \frac{h_{\kappa}^{2q_{\kappa}-2}}{\ell_{\kappa}^{2l_{\kappa}-3}} \bigg( \| \mathcal{E} \mathbf{u} \|_{\mathbf{H}^{l_\kappa}(\mathcal{K}_\kappa)}^2 + \int_0^t \| \mathcal{E} \dot{\mathbf{u}} \|_{\mathbf{H}^{l_\kappa}(\mathcal{K}_\kappa)}^2 + \int_0^t \| \mathcal{E} \ddot{\mathbf{u}} \|_{\mathbf{H}^{l_\kappa}(\mathcal{K}_\kappa)}^2 \bigg) \\
    + & \sum_{\kappa \in \mathcal{T}_h} \frac{h_{\kappa}^{2r_{\kappa}-2}}{\ell_{\kappa}^{2m_\kappa-3}} \bigg( \| \mathcal{E} \mathbf{w} \|_{\mathbf{H}^{m_\kappa}(\mathcal{K}_\kappa)}^2 + \int_0^t \| \mathcal{E} \dot{\mathbf{w}} \|_{\mathbf{H}^{m_\kappa}(\mathcal{K}_\kappa)}^2 + \int_0^t \| \mathcal{E} \ddot{\mathbf{w}} \|_{\mathbf{H}^{m_\kappa}(\mathcal{K}_\kappa)}^2 \bigg) \\
    + & \sum_{\kappa \in \mathcal{T}_h} \frac{h_{\kappa}^{2s_{\kappa}-2}}{\ell_{\kappa}^{2n_{\kappa}-3}} \left( \| \mathcal{E} T \|_{\mathbf{H}^{n_\kappa}(\mathcal{K}_\kappa)}^2 + \| \mathcal{E} \dot{T} \|_{\mathbf{H}^{n_\kappa}(\mathcal{K}_\kappa)}^2 + \int_0^t \| \mathcal{E} \dot{T} \|_{\mathbf{H}^{n_\kappa}(\mathcal{K}_\kappa)}^2  \right),
    \end{aligned}
\end{equation}
where $q_\kappa, r_\kappa$, and $s_\kappa$ are defined as in Lemma \ref{lemma:interp}. 
The hidden constant depends on the time $t$ and on the material properties, but do not depend on the discretization parameters.  
\end{theorem}

\textit{Proof.} The proof of Theorem~\ref{thm:error_est} can be found in \ref{sec:error_analysis_proof}.

\section{Conclusions and further developments}
\label{sec:conclusions}
In this work, we have proposed a new PolyDG discretization method for the fully-dynamic thermo-poroelastic problem. The stability and error analysis for the semi-discrete problem have been performed, establishing a-priori $hp$-error bounds. A wide set of numerical simulations is presented. First, we demonstrated the convergence error bounds of our scheme with respect to both the mesh size and the polynomial degree of approximation. Second, we assessed the capabilities of the proposed formulation addressing literature test cases. Last, we test our approach in physically-sound test cases by observing the wave-propagation phenomenon in thermo-poroelastic media. A comparison with the poroelastic model is presented too, showing the crucial role of temperature in the behavior of the shear waves.

Further developments of this work are possible. First of all, from the point of view of the theoretical analysis, it would be very interesting to include the third-order terms in the energy equation, while from the numerical point of view, the use of effective splitting schemes could allow coping with the high computational cost required for the resolution of the problem in its monolithic formulation. This will be the subject of future research.

\appendix

\section{Proofs of main theorems for the semi-discrete analysis}
\label{sec:staberr_proofs}
In this appendix, we report the proofs of the main results constituting the analysis of the PolyDG-semidiscrete formulation carried out in Section~\ref{sec:semidiscrete_analysis}.

\subsection{Stability estimate (proof of Theorem~\ref{thm:stability_est})}
\label{sec:stability_analysis_proof}

Take $\left( \dot{\mathbf{u}}_h, \dot{\mathbf{v}}_h, T_h \right)$ as test functions in \eqref{eq:semidiscr_form}. We use the skew-simmetry property of the bilinear form $\mathcal{C}_h$ and the symmetry property of the bilinear forms $\mathcal{M}_{uw}, \mathcal{M}_T, \mathcal{A}_{uw,h}$ to get
\begin{equation}
\label{eq:stab_est_1}
\begin{aligned}
& \frac{1}{2}\frac{d}{dt} \bigg[ \mathcal{M}_{uw}((\dot{\mathbf{u}}_h, \dot{\mathbf{w}}_h), \left( \dot{\mathbf{u}}_h, \dot{\mathbf{w}}_h \right)(t) +  \mathcal{M}_T(T_h,T_h)(t) + \mathcal{A}_{uw,h}((\mathbf{u}_h, \mathbf{w}_h), (\mathbf{u}_h, \mathbf{w}_h))(t) \bigg] + \mathcal{B}( \dot{\mathbf{w}}_h, \dot{\mathbf{w}}_h)(t) \\
& + \mathcal{A}_{T,h}(T_h,T_h)(t) = ((\mathbf{f}, \mathbf{g}, H), (\dot{\mathbf{u}}_h, \dot{\mathbf{w}}_h, T_h))(t).
\end{aligned}
\end{equation}
Next, we integrate in time from $0$ to $t\le T_f$ and we obtain
\begin{equation}
\label{eq:stab_est_2}
\begin{aligned}
& \mathcal{M}_{uw}((\dot{\mathbf{u}}_h, \dot{\mathbf{w}}_h), \left( \dot{\mathbf{u}}_h, \dot{\mathbf{w}}_h \right))(t) +  \mathcal{M}_{T}(T_h,T_h)(t) + \mathcal{A}_{uw,h}((\mathbf{u}_h, \mathbf{w}_h), (\mathbf{u}_h, \mathbf{w}_h))(t) + 2 \int_0^t \mathcal{B}( \dot{\mathbf{w}}_h, \dot{\mathbf{w}}_h)(s)ds \\
& + 2 \int_0^t \mathcal{A}_{T,h}(T_h,T_h)(s) ds = 2 \int_0^t ((\mathbf{f}, \mathbf{g}, H), (\dot{\mathbf{u}}_h, \dot{\mathbf{w}}_h, T_h))(s)ds + \mathcal{M}_{uw}((\dot{\mathbf{u}}_h, \dot{\mathbf{w}}_h), \left( \dot{\mathbf{u}}_h, \dot{\mathbf{w}}_h \right))(0) \\
& +  \mathcal{M}_T(T_h,T_h)(0) + \mathcal{A}_{uw,h}((\mathbf{u}_h, \mathbf{w}_h), (\mathbf{u}_h, \mathbf{w}_h))(0).
\end{aligned}
\end{equation}
We use the fundamental theorem of calculus on the $[0,t]$ together with the Cauchy--Schwarz inequality to infer $\mathcal{B}(\mathbf{w}, \mathbf{w})(t) \lesssim \mathcal{B}(\mathbf{w}, \mathbf{w})(0) + \int_0^t \mathcal{B}(\dot{\mathbf{w}}, \dot{\mathbf{w}})(s)ds$. 
Then, applying Lemma~\ref{lemma:cont_coer_stab} and recalling the definition of the energy norms gives
\begin{equation}
\label{eq:stab_est_3}
\begin{aligned}
& \| \left(\mathbf{u}_h, \mathbf{w}_h, T_h \right)(t)\|_{\mathcal{E},*}^2 \lesssim \| \left(\mathbf{u}_h, \mathbf{w}_h, T_h \right)(0)\|_{\mathcal{E}}^2 + 2 \int_0^t ((\mathbf{f}, \mathbf{g}, H), (\dot{\mathbf{u}}_h, \dot{\mathbf{w}}_h, T_h))(s)ds.
\end{aligned}
\end{equation}
Focusing now on the integral term on the right-hand side of \eqref{eq:stab_est_3}, we apply the Cauchy-Schwarz inequality to get
\begin{equation}
    \label{eq:rhs_wp_semidiscr_timeint}
    \begin{aligned}
   \int_0^t ((\mathbf{f}, \mathbf{g}, H), (\dot{\mathbf{u}}_h, \dot{\mathbf{w}}_h, T_h))(s)ds & \leq \int_0^t \| (\mathbf{f}, \mathbf{g}, H)(s) \| \|(\dot{\mathbf{u}}_h, \dot{\mathbf{w}}_h, T_h)(s) \| ds \\
   & \lesssim \int_0^t \| (\mathbf{f}, \mathbf{g}, H)(s) \| \|(\mathbf{u}_h, \mathbf{w}_h, T_h)(s) \|_{\mathcal{E},*} ds 
    \end{aligned}
\end{equation}
Therefore, it is inferred that
\begin{equation}
\label{eq:stab_est_4}
\begin{aligned}
& \| \left(\mathbf{u}_h, \mathbf{w}_h, T_h \right)(t)\|_{\mathcal{E},*}^2 \lesssim \| \left(\mathbf{u}_h, \mathbf{w}_h, T_h \right)(0)\|_{\mathcal{E}}^2 + \int_0^t \| (\mathbf{f}, \mathbf{g}, H)(s) \| \|(\mathbf{u}_h, \mathbf{w}_h, T_h)(s) \|_{\mathcal{E},*} ds 
\end{aligned}
\end{equation}
Finally, we can apply Gromwall's lemma \cite{Quarteroni2017} to \eqref{eq:stab_est_4}. Since \eqref{eq:stab_est_4} holds for an  arbitrary $t \in (0, T_f]$, this concludes the proof.

\subsection{Error estimate (proof of Theorem~\ref{thm:error_est})}
\label{sec:error_analysis_proof}
We start by introducing the following auxiliary result. Owing to Lemma~\ref{lemma:interp} and \cite[Lemma 22, Lemma 33]{CangianiDong:17} we observe that:
\begin{equation}
\label{eq:interp}
\begin{aligned}
|||\left(\mathbf{v} - \mathbf{v}_I , \mathbf{z} - \mathbf{z}_I, S - S_I \right)|||_{dG}^2 \lesssim & \sum_{\kappa \in \mathcal{T}_h} \bigg(  \frac{h_{\kappa}^{2q_{\kappa} - 2}}{\ell_{\kappa}^{2l_{\kappa}-3}} \ ||\mathcal{E}\mathbf{v}||_{\mathbf{H}^{l_{\kappa}}(\mathcal{K}_\kappa)}^2 + \frac{h_{\kappa}^{2r_{\kappa} - 2}}{\ell_{\kappa}^{2m_{\kappa}-3}} \ ||\mathcal{E}\mathbf{z}||_{\mathbf{H}^{m_{\kappa}}(\mathcal{K}_\kappa)}^2 \\
& + \frac{h_{\kappa}^{2s_{\kappa} - 2}}{\ell_{\kappa}^{2n_{\kappa}-3}} \ ||\mathcal{E}S||_{H^{n_{\kappa}}(\mathcal{K}_\kappa)}^2  \bigg),\\
|||\left(\mathbf{v} - \mathbf{v}_I , \mathbf{z} - \mathbf{z}_I, S - S_I \right)|||_{\mathcal{E}}^2 \lesssim & \sum_{\kappa \in \mathcal{T}_h} \bigg( \frac{h_{\kappa}^{2q_{\kappa}}}{\ell_{\kappa}^{2l_{\kappa}}} \ ||\mathcal{E} \dot{\mathbf{v}}||_{\mathbf{H}^{l_{\kappa}}(\mathcal{K}_\kappa)}^2 + \frac{h_{\kappa}^{2r_{\kappa}}}{\ell_{\kappa}^{2m_{\kappa}}} \ ||\mathcal{E} \dot{\mathbf{z}}||_{\mathbf{H}^{m_{\kappa}}(\mathcal{K}_\kappa)}^2 \\
& +\frac{h_{\kappa}^{2s_{\kappa}}}{\ell_{\kappa}^{2n_{\kappa}}} \ ||\mathcal{E}S||_{H^{n_\kappa}(\mathcal{K}_\kappa)}^2  \bigg) + \sum_{\kappa \in \mathcal{T}_h} \bigg( \frac{h_{\kappa}^{2q_{\kappa} - 2}}{\ell_{\kappa}^{2l_{\kappa}-3}} \ ||\mathcal{E}\mathbf{v}||_{\mathbf{H}^{l_{\kappa}}(\mathcal{K}_\kappa)}^2 \\
& + \frac{h_{\kappa}^{2r_{\kappa} - 2}}{\ell_{\kappa}^{2m_{\kappa}-3}} \ ||\mathcal{E}\mathbf{z}||_{\mathbf{H}^{m_{\kappa}}(\mathcal{K}_\kappa)}^2 + \frac{h_{\kappa}^{2s_{\kappa} - 2}}{\ell_{\kappa}^{2n_{\kappa}-3}} \ ||\mathcal{E}S||_{H^{n_{\kappa}}(\mathcal{K}_\kappa)}^2  \bigg).\\
\end{aligned}
\end{equation}

To derive the error equation for our problem we need to extend the bilinear forms \eqref{eq:bil_form_discr} to the space of continuous solutions. Thus, we need further regularity assumptions on the exact solutions $\textbf{X}$. Indeed, we consider the solid displacement, filtration displacement, and temperature to have at least local $H^2$-regularity, as reported in the statement of Theorem \ref{thm:error_est}. Moreover, without any loss of generality, we assume the continuity of the normal stress, of the heat flux, and of the normal traces of the two velocities $\dot{\mathbf{u}}$ and $\dot{\mathbf{w}}$ across the interfaces $F\in\mathcal{F}_I$ for all time $t\in (0,T_f]$. Under these assumptions, we can insert the exact solutions into \eqref{eq:semidiscr_form} obtaining a formulation equivalent to \eqref{eq:DYN_weak_form}. Now, we can subtract the resulting equation from \eqref{eq:semidiscr_form} to infer the error equation
\begin{equation}
    \label{eq:error_eq}
    \begin{aligned}
    & \mathcal{M}_{uw}((\ddot{\mathbf{e}}^u, \ddot{\mathbf{e}}^w), \left( \mathbf{v}_h, \mathbf{z}_h \right)) + \mathcal{B}( \dot{\mathbf{e}}^w, \mathbf{z}_h) +  \mathcal{M}_T(\dot{e}^T,S_h) + \mathcal{C}_h(\left(\dot{\mathbf{e}}^u, \dot{\mathbf{e}}^w \right), S_h) + \mathcal{A}_{uw,h}((\mathbf{e}^u, \mathbf{e}^u), \left( \mathbf{v}_h, \mathbf{z}_h \right)) \\
    & + \mathcal{A}_{T,h}(e^T,S_h) - \mathcal{C}_h(\left( \mathbf{v}_h, \mathbf{z}_h \right), e^T) = 0,
\end{aligned}
\end{equation}
for all $(\mathbf{v}_h, \mathbf{w}_h, S_h) \in \mathbf{V}_h\times \mathbf{V}_h \times V_h$. 
We assume that the semi-discrete problem \eqref{eq:semidiscr_form} is completed by initial conditions $\textbf{X}_h(0) = (\mathbf{u}_I(0), \mathbf{w}_I(0), T_I(0))$ and $(\dot{\mathbf{u}}_h(0), \dot{\mathbf{w}}_h(0)) = (\dot{\mathbf{u}}_I(0), \dot{\mathbf{w}}_I(0))$  where $\mathbf{u}_I, \mathbf{w}_I, T_I$ are the interpolants of the exact solutions given by Lemma~\ref{lemma:interp}, so that the error equation \eqref{eq:error_eq} is supplemented by the condition
\begin{equation}
\label{eq:IC_error}
E_h(0) = \mathbf{0}, \quad \dot{\mathbf{e}}_h^u = \mathbf{0}, \quad \dot{\mathbf{e}}_h^w = \mathbf{0}.
\end{equation}
We now take $\left( \mathbf{v}_h, \mathbf{z}_h, S_h \right) = \left( \dot{\mathbf{e}}_h^u, \dot{\mathbf{e}}_h^w, e_h^T \right)$ in \eqref{eq:error_eq}, bringing to the right-hand side all the terms that involve the interpolation errors. Thus, we get
\begin{equation}
    \label{eq:error_eq1}
    \begin{aligned}
    & \mathcal{M}_{uw}((\ddot{\mathbf{e}}_h^u, \ddot{\mathbf{e}}_h^w), (\dot{\mathbf{e}}_h^u, \dot{\mathbf{e}}_h^w)) + \mathcal{B}( \dot{\mathbf{e}}_h^w, \dot{\mathbf{e}}_h^w) +  \mathcal{M}_T(\dot{e}_h^T,e_h^T) + \mathcal{A}_{uw,h}((\mathbf{e}_h^u, \mathbf{e}_h^w), (\dot{\mathbf{e}}_h^u, \dot{\mathbf{e}}_h^w)) + \mathcal{A}_{T,h}(e_h^T,e_h^T) \\
    & =  \mathcal{M}_{uw}((\ddot{\mathbf{e}}_I^u, \ddot{\mathbf{e}}_I^w), (\dot{\mathbf{e}}_h^u, \dot{\mathbf{e}}_h^w)) + \mathcal{B}( \dot{\mathbf{e}}_I^w, \dot{\mathbf{e}}_h^w) +  \mathcal{M}_T(\dot{e}_I^T,e_h^T) + \mathcal{C}_h(\left(\dot{\mathbf{e}}_I^u, \dot{\mathbf{e}}_I^w \right), e_h^T) + \mathcal{A}_{uw,h}((\mathbf{e}_I^u, \mathbf{e}_I^w), (\dot{\mathbf{e}}_h^u, \dot{\mathbf{e}}_h^w)) \\
    & + \mathcal{A}_{T,h}(e_I^T,e_h^T) - \mathcal{C}_h((\dot{\mathbf{e}}_h^u, \dot{\mathbf{e}}_h^w), e_I^T).
\end{aligned}
\end{equation}
For treating the left-hand side of \eqref{eq:error_eq1} we follow the same arguments of \ref{sec:stability_analysis_proof}. For the right-hand side, we move the time derivatives from the discretization errors to the interpolation ones in the fifth and seventh bilinear forms via Leibniz's formula. By doing so, we infer that
\begin{equation}
    \label{eq:error_eq2}
    \begin{aligned}
    & \frac{1}{2}  \frac{d}{dt} \bigg[ \mathcal{M}_{uw}((\dot{\mathbf{e}}_h^u, \dot{\mathbf{e}}_h^w), (\dot{\mathbf{e}}_h^u, \dot{\mathbf{e}}_h^w)) +  \mathcal{M}_T(e_h^T,e_h^T) + \mathcal{A}_{uw,h}((\mathbf{e}_h^u, \mathbf{e}_h^w), (\mathbf{e}_h^u, \mathbf{e}_h^w)) \bigg] + \mathcal{B}( \dot{\mathbf{e}}_h^w, \dot{\mathbf{e}}_h^w) + \mathcal{A}_{T,h}(e_h^T,e_h^T) \\
    & =  \mathcal{M}_{uw}((\ddot{\mathbf{e}}_I^u, \ddot{\mathbf{e}}_I^w), (\dot{\mathbf{e}}_h^u, \dot{\mathbf{e}}_h^w)) \hspace{-0.75pt} + \hspace{-0.75pt} \mathcal{B}( \dot{\mathbf{e}}_I^w, \dot{\mathbf{e}}_h^w) +  \mathcal{M}_T(\dot{e}_I^T,e_h^T) + \mathcal{C}_h(\left(\dot{\mathbf{e}}_I^u, \dot{\mathbf{e}}_I^w \right), e_h^T) \hspace{-1.5pt} + \hspace{-1.5pt} \frac{d}{dt} \mathcal{A}_{uw,h}((\mathbf{e}_I^u, \mathbf{e}_I^w), (\mathbf{e}_h^u, \mathbf{e}_h^w)) \\
    & - \mathcal{A}_{uw,h}((\dot{\mathbf{e}}_I^u, \dot{\mathbf{e}}_I^w), (\mathbf{e}_h^u, \mathbf{e}_h^w)) + \mathcal{A}_{T,h}(e_I^T,e_h^T) - 
    \frac{d}{dt} \mathcal{C}_h((\mathbf{e}_h^u, \mathbf{e}_h^w), e_I^T) + \mathcal{C}_h((\mathbf{e}_h^u, \mathbf{e}_h^w), \dot{e}_I^T).
\end{aligned}
\end{equation}
Now, we integrate with respect to time between $0$ and $t \leq T_f$, recalling \eqref{eq:IC_error}, and owing on the same arguments as in the proof of Theorem~\ref{thm:stability_est} we obtain
\begin{equation}
    \label{eq:error_eq3}
    \begin{aligned}
    & \|E_h(t)\|_{\mathcal{E}}^2 + \int_0^t \|e_h^T(s)\|_{dG,T}^2 \, ds \lesssim \mathcal{R}_1(t) + \int_0^t \bigg( \mathcal{R}_2(s) + \mathcal{R}_3(s) \bigg) \, ds,
    \end{aligned}
\end{equation}
where
\begin{equation}
    \label{eq:error_rhs}
    \begin{aligned}
    & \mathcal{R}_1 = \mathcal{A}_{uw,h}((\mathbf{e}_I^u, \mathbf{e}_I^w), (\mathbf{e}_h^u, \mathbf{e}_h^w))  - \mathcal{C}_h((\mathbf{e}_h^u, \mathbf{e}_h^w), e_I^T), \\
    & \mathcal{R}_2 = \mathcal{M}_{uw}((\ddot{\mathbf{e}}_I^u, \ddot{\mathbf{e}}_I^w), (\dot{\mathbf{e}}_h^u, \dot{\mathbf{e}}_h^w)) + \mathcal{B}( \dot{\mathbf{e}}_I^w, \dot{\mathbf{e}}_h^w) +  \mathcal{M}_T(\dot{e}_I^T,e_h^T), \\
    & \mathcal{R}_3 =  \mathcal{C}_h(\left(\dot{\mathbf{e}}_I^u, \dot{\mathbf{e}}_I^w \right), e_h^T) - \mathcal{A}_{uw,h}((\dot{\mathbf{e}}_I^u, \dot{\mathbf{e}}_I^w), (\mathbf{e}_h^u, \mathbf{e}_h^w)) + \mathcal{A}_{T,h}(e_I^T,e_h^T) + \mathcal{C}_h((\mathbf{e}_h^u, \mathbf{e}_h^w), \dot{e}_I^T).
    \end{aligned}
\end{equation}
We bound the terms $\mathcal{R}_1, \mathcal{R}_2, \mathcal{R}_3$ by the repeated use of Cauchy-Schwarz, Young and triangles inequalities, and Lemma~\ref{lemma:bound_bilinearform_error}:
\begin{equation}
    \label{eq:error_rhs_bound}
    \begin{aligned}
    \mathcal{R}_1 \lesssim & \ \left( |||\mathbf{e}_I^u|||_{dG,e}^2 + |||\mathbf{e}_I^u|||_{dG,p}^2 + |||\mathbf{e}_I^w|||_{dG,p}^2 + |||e_I^T|||_{dG,T}^2 \right) + \left( \|\mathbf{e}_h^u\|_{dG,e}^2 + |\mathbf{e}_h^u|_{dG,p}^2 + |\mathbf{e}_h^w|_{dG,p}^2 \right) \\[5pt]
    \mathcal{R}_2 \lesssim & \ \mathcal{M}_{uw}((\ddot{\mathbf{e}}_I^u, \ddot{\mathbf{e}}_I^w), (\ddot{\mathbf{e}}_I^u, \ddot{\mathbf{e}}_I^w)) + \mathcal{B}( \dot{\mathbf{e}}_I^w, \dot{\mathbf{e}}_I^w) +  \mathcal{M}_T(\dot{e}_I^T,\dot{e}_I^T) + \mathcal{M}_{uw}((\dot{\mathbf{e}}_h^u, \dot{\mathbf{e}}_h^w), (\dot{\mathbf{e}}_h^u, \dot{\mathbf{e}}_h^w))
     + \mathcal{B}( \dot{\mathbf{e}}_h^w, \dot{\mathbf{e}}_h^w)  \\
     & +  \mathcal{M}_T(e_h^T,e_h^T) \\[5pt]
    \mathcal{R}_3 \lesssim & \ \big( ||| \dot{\mathbf{e}}_I^u |||_{dG,e}^2 + ||| \dot{\mathbf{e}}_I^u |||_{dG,p}^2 + ||| \dot{\mathbf{e}}_I^w |||_{dG,p}^2 + ||| \dot{e}_I^T |||_{dG,T}^2 + ||| e_I^T |||_{dG,T}^2 \big) + \big( \| \mathbf{e}_h^u \|_{dG,e}^2 + | \mathbf{e}_h^u |_{dG,p}^2 \\
    \ & + | \mathbf{e}_h^w |_{dG,p}^2 + \|e_h^T\|_{dG,T}^2 \big) 
    \end{aligned}
\end{equation}
By plugging \eqref{eq:error_rhs_bound} into \eqref{eq:error_eq3} we get
\begin{equation}
    \label{eq:error_eq4}
    \begin{aligned}
     \|E_h(t)\|_{\mathcal{E}}^2 + \int_0^t \|e_h^T(s)\|_{dG,T}^2 \, ds \lesssim & 
     |||E_I(t)|||_{dG}^2 + \int_0^t \|E_h(s)\|_{\mathcal{E}}^2 ds \\
    & + \int_0^t \left( |||\dot{E}_I(s)|||_{\mathcal{E}}^2 + |||e_I^T(s)|||_{dG,T}^2 \right) ds 
    \end{aligned}
\end{equation}
and by using Gronwall's lemma we obtain
\begin{equation}
    \label{eq:error_eq5}
    \begin{aligned}
    & \|E_h(t)\|_{\mathcal{E}}^2 + \int_0^t \|e_h^T(s)\|_{dG,T}^2 \, ds \lesssim |||E_I(t)|||_{dG}^2 + \int_0^t \left( |||\dot{E}_I(s)|||_{\mathcal{E}}^2 + |||e_I^T(s)|||_{dG,T}^2 \right) ds.
    \end{aligned}
\end{equation}
The thesis follows by bounding the right-hand side of equation \eqref{eq:error_eq5} via the interpolation estimates of Lemma~\ref{lemma:interp} and \eqref{eq:interp}.

\bigskip

\textbf{CRediT authorship contribution statement}

\textbf{S. Bonetti:} Conceptualization, Data curation, Formal analysis, Investigation, Methodology, Software, Validation, Visualization, Writing – original draft. \textbf{M. Botti:} Conceptualization, Formal analysis, Methodology, Writing – review \& editing. \textbf{I. Mazzieri:} Formal analysis, Methodology, Software, Validation, Visualization, Writing – review \& editing. \textbf{P.F. Antonietti:} Conceptualization, Funding acquisition, Methodology, Project administration, Supervision, Writing – review \& editing.

\bigskip
\textbf{Declaration of competing interest}

The authors declare that they have no known competing financial interests or personal relationships that could have appeared to influence the work reported in this paper.

\bigskip
\textbf{Acknowledgements}

This work has received funding from the European Union’s Horizon 2020 research and innovation programme under the Marie Skłodowska-Curie grant agreement No. 896616 (project PDGeoFF).
P.F.A. has been partially funded by the research grants PRIN2017 n. 201744KLJL and PRIN2020 n. 20204LN5N5 funded by the Italian Ministry of Universities and Research (MUR). P.F.A. and I.M. have been partially funded by European Union - Next Generation EU. S.B., M.B., I.M., and P.F.A. are members of INdAM-GNCS. The work of M.B. has been partially supported by the INdAM-GNCS project CUP E55F22000270001.

\bibliography{bibliography.bib}

\end{document}